\documentclass{gtart_h}  

\def\ifplaintex{\expandafter\ifx\csname documentclass\endcsname\relax}

\def\ifplaintex{\expandafter\ifx\csname documentclass\endcsname\relax}

\ifplaintex 
\else
\fi

\expandafter\ifx\csname epsfbox\endcsname\relax\input epsf\fi

\def\gt{{\mathsurround=0pt\it $\cal G\mskip-2mu$eometry \&\ 
$\cal T\!\!$opology}}        %  journal title in recommended style

\def\gtp{{\mathsurround=0pt\it $\cal G\mskip-2mu$eometry \&\ 
$\cal T\!\!$opology $\cal P\!$ublications}}  % GT publications

\def\lognumber#1{\def\thelognumber{#1}}
\def\volumenumber#1{\def\thevolumenumber{#1}}
\def\papernumber#1{\def\thepapernumber{#1}}
\def\volumeyear#1{\def\thevolumeyear{#1}}

\def\pagenumbers#1#2{\def\startpage{#1}\def\finishpage{#2}}
\def\published#1{\def\publishdate{#1}}
\def\proposed#1{\def\theproposer{#1}}
\def\seconded#1{\def\theseconders{#1}}
\def\received#1{\def\receiveddate{#1}}

\def\accepted#1{\def\accepteddate{#1}}
\def\asciititle#1{\def\theasciititle{#1}}

\def\asciiaddress#1{\def\theasciiaddress{#1}}
\def\asciiemail#1{\def\theasciiemail{#1}}

\long\def\asciiabstract#1{\long\def\theasciiabstract{#1}}
\def\asciikeywords#1{\def\theasciikeywords{#1}}

\let\\\par\let\thelognumber\relax
\let\thevolumenumber\relax\let\thepapernumber\relax
\let\thevolumeyear\relax\let\thesamplenumber\relax\let\startpage\relax
\let\finishpage\relax\let\publishdate\relax\let\receiveddate\relax
\let\reviseddate\relax\let\accepteddate\relax\let\theasciititle\relax
\let\theasciiauthors\relax\let\theasciiaddress\relax
\let\theasciiabstract\relax\let\theasciikeywords\relax
\let\theasciiemail\relax\let\theshortauthors\relax\let\theshorttitle\relax

\long\def\maketitlep{   % start of definition of \maketitlep

\count0=\startpage

\gt\hfill      %   Journal title (top left) 
\hbox to 77pt{\vbox to 0pt{\vglue -15pt\epsfbox{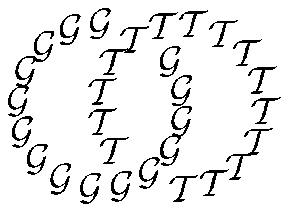}\vss}\hss}
\break
{\small\ifx\thesamplenumber\relax % sample?  
Volume \else Sample
\fi\thevolumenumber\ (\thevolumeyear)
\startpage--\finishpage\nl
Published: \publishdate}
\vglue 0.5truein plus 0.4fil minus 0.1truein

{\parskip=0pt\leftskip 0pt plus 1fil\def\\{\par\smallskip}{\ifplaintex\large
\else\Large\fi\bf\thetitle}\par\medskip}   

\vglue 0pt plus 0.1fil 

{\parskip=0pt\leftskip 0pt plus 1fil\def\\{\par}{\sc\theauthors}
\par\medskip}

\vglue 0pt plus 0.1fil 

{\small\parskip=0pt\let\newline\\
{\leftskip 0pt plus 1fil\def\\{\par}{\sl\theaddress}\par}
\expandafter\ifx\theemail\relax    % email address?
\relax\else\vglue 5pt plus 0.02fil minus 2pt\def\\{\stdspace{\rm 
and}\stdspace} 
\cl{Email:\stdspace\tt\theemail}\fi
\ifx\theurl\relax                  % URL given?
\relax\else\vglue 5pt plus 0.02fil minus 2pt\def\\{\stdspace{\rm 
and}\stdspace}
\cl{URL:\stdspace\tt\theurl}\fi\par}

\vglue 7pt plus 0.3fil minus 3pt

{\bf Abstract}
\vglue 5pt plus 0.1fil minus 2pt

\theabstract

\vglue 7pt plus 0.3fil minus 3pt

{\bf AMS Classification numbers}\quad Primary:\quad \theprimaryclass

Secondary:\quad \thesecondaryclass

\vglue 5pt plus 0.3fil minus 2pt

{\bf Keywords:}\quad \thekeywords

\vglue 10pt plus 0.5fil minus 5pt

{\small  Proposed: \theproposer\hfill Received: \receiveddate\nl
Seconded: \theseconders\hfill 
\ifx\reviseddate\relax                         % paper revised?
Accepted: \accepteddate                        % no
\else
Revised: \reviseddate                          % yes
\fi}
\eject
}       %  end of definition of \maketitlep

\font\phead=cmsl9 scaled 950
\font=cmsl9 scaled 1050
\font\pnum=cmbx10 scaled 913
\font\lnum=cmbx10 
\font\pfoot=cmsl9 scaled 950
\font=cmsl9 scaled 1050
\ifplaintex
\headline{\vbox to 0pt{\vskip -4.5mm\line{\small\phead\ifnum
\count0=\startpage ISSN 1364-0380 (on line)
1465-3060 (printed) \hfill {\pnum\folio}\else\ifodd\count0\def\\{ }% 
\ifx\theshorttitle\relax\thetitle\else\theshorttitle\fi\hfill{\pnum\folio}
\else\def\\{ and }{\pnum\folio}\hfill\ifx\theshortauthors\relax\theauthors
\else\theshortauthors\fi\fi\fi}\vss}}
\footline{\vbox to 0pt{\vglue 0mm\line{\small\pfoot\ifnum\count0=\startpage
\copyright\ \gtp\hfill\else
\gt, Volume \thevolumenumber\ (\thevolumeyear)\hfill\fi}\vss
}}
\else
\makeatletter
\def\@oddhead{{\small\lhead\ifnum\count0=\startpage ISSN 1364-0380 (on line)
1465-3060 (printed) \hfill {\lnum\number\count0}\else\ifodd\count0
\def\\{ }\ifx\theshorttitle\relax \thetitle \else\theshorttitle\fi\hfill
{\lnum\number\count0}\else\def\\{ and }{\lnum\number\count0}
\hfill\ifx\theshortauthors\relax 
\theauthors\else\theshortauthors\fi\fi\fi}}\def\@evenhead{\@oddhead}
\def\@oddfoot{\small\lfoot\ifnum\count0=\startpage\copyright\ \gtp\hfill\else
\gt, Volume \thevolumenumber\ (\thevolumeyear)\hfill\fi}
\def\@evenfoot{\@oddfoot}
\makeatother
\fi

\newwrite\gtoutfile
\long\gdef\makeheadfile{  %%% start of definition of \makeheadfile
{\def\\{, }\def\s{ }
\immediate\openout\gtoutfile head.xxx
\immediate\write\gtoutfile{Proxy-for: \ifx\theasciiauthors\relax
\theauthors\else\theasciiauthors\fi\s<\ifx\theasciiemail\relax\theemail\else\theasciiemail\fi>}
\immediate\write\gtoutfile{\noexpand\\}
\immediate\write\gtoutfile{Authors: \ifx\theasciiauthors\relax
\theauthors\else\theasciiauthors\fi}
{\def\\{ }\immediate\write\gtoutfile{Title: \ifx\theasciititle\relax
\thetitle\else\theasciititle\fi}}
\immediate\write\gtoutfile{Subj-class: GT or SG or MG etc}
\immediate\write\gtoutfile{MSC-class: \theprimaryclass\ifx\thesecondaryclass\relax\else, \thesecondaryclass\fi}
\immediate\write\gtoutfile{Journal-ref: Geom. Topol. \thevolumenumber
(\thevolumeyear) \startpage-\finishpage}
\immediate\write\gtoutfile{Comments: Published by Geometry and Topology at}
\immediate\write\gtoutfile{\s\s http://www.maths.warwick.ac.uk/gt/GTVol\thevolumenumber/paper\thepapernumber.abs.html}
\immediate\write\gtoutfile{\noexpand\\}
\immediate\write\gtoutfile{}
\ifx\theasciiabstract\relax
\immediate\write\gtoutfile{\theabstract}\else
\immediate\write\gtoutfile{\theasciiabstract}\fi
\immediate\write\gtoutfile{}
\immediate\write\gtoutfile{\noexpand\\}
\immediate\write\gtoutfile{}
\immediate\closeout\gtoutfile}}  %%% end of definition of \makeheadfile

\def\maketitlepage{\maketitlep\makeheadfile}
\let\maketitle\maketitlepage

\lognumber{412}

\volumenumber{8}\papernumber{23}\volumeyear{2004}
\pagenumbers{877}{924}
\received{26 January 2004}
\published{5 June 2004}
\accepted{28 May 2004}
\proposed{Robion Kirby}
\seconded{Martin Bridson, Steven Ferry}

\usepackage{amssymb, amsmath, latexsym}
\newcount\localcountcount
\newcount\localboxcount
\def\newlocalcount#1{
	\advance\localcountcount by -1
	\allocationnumber= \localcountcount
	\ifnum\count10<\allocationnumber
	  \else
	    \errmessage{No room for a new local count.}
	  \fi
	\countdef#1=\allocationnumber
    }

\def\newlocalbox#1{
	\advance\localboxcount by -1
	\allocationnumber= \localboxcount
	\ifnum\count14<\allocationnumber
	  \else
	    \errmessage{No room for a new local box.}
	  \fi
	\chardef#1=\allocationnumber
    }

\newcount\coordinatescalex
\newcount\coordinatescaley
\def\bigger(#1/#2){
	\multiply\coordinatescalex by #1
	\divide\coordinatescalex by #2
	\multiply\coordinatescaley by #1
	\divide\coordinatescaley by #2
    }
\def\wider(#1/#2){
	\multiply\coordinatescalex by #1
	\divide\coordinatescalex by #2
    }
\def\taller(#1/#2){
	\multiply\coordinatescaley by #1
	\divide\coordinatescaley by #2
    }

\newlength{\objectmargin}
\newlength{\arrowmargin}
\newlength{\cornersize}
\newlength{\arrowheadlength}
\newlength{\nudgesize}
\def\enlarge #1 by #2/#3;{
	\multiply\csname #1\endcsname by #2
	\divide\csname #1\endcsname by #3
    }

\newif\ifflippinglabels
\def\fliplabel{\global\flippinglabelstrue}

\newif\ifepimorphic

\newif\ifmonomorphic

\newcount\labelposition
\newcount\nudgex
\newcount\nudgey
\def\nudge(#1,#2){
	\global\advance\nudgex by #1
	\global\advance\nudgey by #2
    }
\def\nudgeright{
	\global\multiply\nudgex by \unitlength
	\global\advance\nudgex by \nudgesize
	\global\divide\nudgex by \unitlength
    }
\def\nudgeleft{
	\global\multiply\nudgex by \unitlength
	\global\advance\nudgex by -\nudgesize
	\global\divide\nudgex by \unitlength
    }
\def\nudgeup{
	\global\multiply\nudgey by \unitlength
	\global\advance\nudgey by \nudgesize
	\global\divide\nudgey by \unitlength
    }
\def\nudgedown{
	\global\multiply\nudgey by \unitlength
	\global\advance\nudgey by -\nudgesize
	\global\divide\nudgey by \unitlength
    }

\def\negate(#1){
	\csname #1\endcsname= -\csname #1\endcsname
    }

\def\makepositive(#1){
	\ifnum\csname #1\endcsname<0
	    \negate(#1)
	  \fi
    }

\def\calcgdivisor(#1,#2) into #3 {{
	\newlocalcount\left
	\newlocalcount\right
	\newlocalcount\temp

	\left=  #1
	\right= #2
	\makepositive(left)
	\makepositive(right)

	\loop
	    \ifnum\left<\right
		\temp=  \left
		\left=  \right
		\right= \temp
	      \fi
	  \ifnum\right>0 
	    \advance\left by -\right
	  \repeat
	\global\csname #3\endcsname= \left
    }}

\def\newpoint #1 {
	\expandafter\newlocalcount\csname #1x\endcsname
	\expandafter\newlocalcount\csname #1y\endcsname
    }

\def\setpoint #1 to (#2,#3){
	\csname #1x\endcsname = #2
	\csname #1y\endcsname = #3
    }

\def\copypoint #1 to #2 {
	\csname #2x\endcsname = \csname #1x\endcsname
	\csname #2y\endcsname = \csname #1y\endcsname
    }

\def\addpoint #1 to #2 {
	\advance\csname #2x\endcsname by \csname #1x\endcsname
	\advance\csname #2y\endcsname by \csname #1y\endcsname
    }

\def\subtractpoint #1 from #2 {
	\advance\csname #2x\endcsname by -\csname #1x\endcsname
	\advance\csname #2y\endcsname by -\csname #1y\endcsname
    }

\def\multiplypoint #1 by #2 {
	\multiply\csname #1x\endcsname by #2
	\multiply\csname #1y\endcsname by #2
    }

\def\dividepoint #1 by #2 {
	\divide\csname #1x\endcsname by #2
	\divide\csname #1y\endcsname by #2
    }

\def\printpoint #1 {
     \message{(\number\csname #1x\endcsname,\number\csname #1y\endcsname)}
    }

\def\lookuptrig #1:#2 {
	\ifnum1=#1
	    \ifnum0=#2		% 1:0
		\multiply\trigpointx by 10000
		\multiply\trigpointy by 0
	      \else		% 1:1
		\multiply\trigpointx by 7071
		\multiply\trigpointy by 7071
	      \fi
	  \fi
	\ifnum2=#1		% 2:1
	    \multiply\trigpointx by 8944
	    \multiply\trigpointy by 4472
	  \fi	
	\ifnum3=#1
	    \ifnum1=#2		% 3:1
		\multiply\trigpointx by 9487
		\multiply\trigpointy by 3162
	      \else		% 3:2
		\multiply\trigpointx by 8321
		\multiply\trigpointy by 5547
	      \fi
	  \fi	
	\ifnum4=#1
	    \ifnum1=#2		% 4:1
		\multiply\trigpointx by 9701
		\multiply\trigpointy by 2425
	      \else		% 4:3
		\multiply\trigpointx by 8000
		\multiply\trigpointy by 6000
	      \fi
	  \fi	
    }

\def\calcunsignedtrig #1:#2 {
	\ifnum#1<#2
	    \lookuptrig #2:#1
	    \advance\trigpointx by -\trigpointy
	    \advance\trigpointy by \trigpointx
	    \advance\trigpointx by -\trigpointy
	    \trigpointx= -\trigpointx
	  \else
	    \lookuptrig #1:#2
	  \fi
    }

\def\calctrig #1:#2 {
	\ifnum0>#2
	    \ifnum0>#1
		\calcunsignedtrig -#1:-#2
	        \trigpointx= -\trigpointx
	      \else
		\calcunsignedtrig #1:-#2
	      \fi
	    \trigpointy= -\trigpointy
	  \else
	    \ifnum0>#1
		\calcunsignedtrig -#1:#2
	        \trigpointx= -\trigpointx
	      \else
		\calcunsignedtrig #1:#2
	      \fi
	  \fi
    }

\def\calctrigoffset #1 along #2:#3 {
	\trigpointx= #1
	\divide\trigpointx by \unitlength
	\trigpointy= \trigpointx
	\calctrig #2:#3
	\dividepoint trigpoint by 10000
    }

\def\trigoffset #1 by #2 along #3:#4 {
	\calctrigoffset #2 along #3:#4
	\addpoint trigpoint to #1
    }

\def\newvector #1 {
	\newpoint {#1head}
	\newpoint {#1tail}
	\newpoint {#1slope}
	\expandafter\newlocalcount\csname #1length\endcsname
    }

\def\shortenhead #1 by #2 {
	\subtractpoint {#2} from {#1head}
	\newlocalcount\changeinlength
	\ifnum\csname #1slopex\endcsname=0
	    \changeinlength= \csname #2y\endcsname
	  \else
	    \changeinlength= \csname #2x\endcsname
	  \fi
	\makepositive(changeinlength)
	\advance\csname #1length\endcsname by -\changeinlength
    }

\def\shortentail #1 by #2 {
	\addpoint {#2} to {#1tail}
	\newlocalcount\changeinlength
	\ifnum\csname #1slopex\endcsname=0
	    \changeinlength= \csname #2y\endcsname
	  \else
	    \changeinlength= \csname #2x\endcsname
	  \fi
	\makepositive(changeinlength)
	\advance\csname #1length\endcsname by -\changeinlength
    }

\def\setvector #1 from #2 to #3 {
	\copypoint #2 to {#1tail}
	\copypoint #3 to {#1head}

	\copypoint #3 to {#1slope}
	\subtractpoint #2 from {#1slope}

	\ifnum\csname #1slopex\endcsname=0
	    \csname #1length\endcsname= \csname #1slopey\endcsname
	  \else
	    \csname #1length\endcsname= \csname #1slopex\endcsname
	  \fi
	\makepositive(#1length)

	\newlocalcount\gdivisor
	\calcgdivisor(
		\csname #1slopex\endcsname,
		\csname #1slopey\endcsname
	    ) into gdivisor
	\dividepoint {#1slope} by {\gdivisor}
    }

\def\getvectormidpoint #1 into #2
    {
	\copypoint {#1head} to #2
	\subtractpoint  {#1tail} from #2
	\multiplypoint #2 by {\labelposition}
	\dividepoint #2 by 100
	\addpoint {#1tail} to #2
	\global\labelposition=50
    }

\def\drawline #1 {
	\put(\csname #1tailx\endcsname,\csname #1taily\endcsname){\line( \csname #1slopex\endcsname,\csname #1slopey\endcsname){\csname #1length\endcsname}}
    }

\def\drawsecondhead #1 {{
	\newpoint {}
	\newpoint trigpoint
	\copypoint {#1head} to {}
	\trigoffset {}
		by {-\arrowheadlength}
		along {\csname #1slopex\endcsname}:{\csname #1slopey\endcsname}
\put(\x,\y){\vector(\csname #1slopex\endcsname,\csname #1slopey\endcsname){0}}
    }}

\def\drawarrowtail #1 {{
	\newpoint {}
	\newpoint trigpoint
	\copypoint {#1tail} to {}
	\trigoffset {}
		by {\arrowheadlength}
		along {\csname #1slopex\endcsname}:{\csname #1slopey\endcsname}
\put(\x,\y){\vector(\csname #1slopex\endcsname,\csname #1slopey\endcsname){0}}
    }}

\def\drawvector #1 {
	\put(\csname #1tailx\endcsname,\csname #1taily\endcsname){\vector( \csname #1slopex\endcsname,\csname #1slopey\endcsname){\csname #1length\endcsname}}
	\ifepimorphic
	    \drawsecondhead {#1}
	    \global\epimorphicfalse
	  \fi	
	\ifmonomorphic
	    \drawarrowtail {#1}
	    \global\monomorphicfalse
	  \fi	
    }

\def\drawdoublevector #1 {
	\drawvector {#1}
	\put(\csname #1tailx\endcsname,\csname #1taily\endcsname){\vector(-\csname #1slopex\endcsname,-\csname #1slopey\endcsname){ 0 }}
    }

\def\expandxbounds(#1){
	\ifnum\leastobjectx>#1
	    \global\leastobjectx=#1
	  \fi
	\ifnum\greatestobjectx<#1
	    \global\greatestobjectx=#1
	  \fi
    }

\def\expandybounds(#1,#2){
	\ifnum\leastobjecty>#2
	    \global\leastobjecty=#2
	  \fi
	\ifnum\greatestobjecty<#1
	    \global\greatestobjecty=#1
	  \fi
    }

\def\object #1 #2,#3 #4;{
	\expandafter\newlocalcount\csname objectx#1\endcsname
	\expandafter\newlocalcount\csname objecty#1\endcsname
	\csname objectx#1\endcsname=#2
	\csname objecty#1\endcsname=#3
	\multiply\csname objectx#1\endcsname by \coordinatescalex
	\multiply\csname objecty#1\endcsname by \coordinatescaley

	\expandafter\newlocalcount\csname halfwidth#1\endcsname
	\expandafter\newlocalcount\csname halfheight#1\endcsname
	{
	    \newlocalbox\objectbox

	    \savebox{\objectbox}{$\displaystyle{#4}$}
	    \ifdim\ht\objectbox=0pt 
		\savebox{\objectbox}{$\bullet$}
	      \fi

	    \newpoint {}
	    \setpoint {} to (#2,#3)
	    \multiply\x by \coordinatescalex
	    \multiply\y by \coordinatescaley

	    \put(\x,\y){\makebox(0,0){\usebox{\objectbox}}}
	
	    \newlocalcount\halfwidth
	    \newlocalcount\halfheight

	    \expandxbounds(\x)

	    \halfwidth=  \wd\objectbox
	    \halfheight= \ht\objectbox
	    \advance\halfheight by \dp\objectbox

	    \divide\halfwidth  by 2
	    \divide\halfheight by 2

	    \divide\halfwidth  by \unitlength
	    \divide\halfheight by \unitlength

	    \newlocalcount\top
	    \newlocalcount\bottom
	    \top= \y
	    \bottom= \y
	    \advance\top    by \halfheight
	    \advance\bottom by -\halfheight
	    \expandybounds(\top,\bottom)

	    \newlocalcount\margin
	    \margin= \objectmargin
	    \divide\margin by \unitlength
	    \advance\halfwidth  by \margin
	    \advance\halfheight by \margin

	    \global\csname halfwidth#1\endcsname=  \halfwidth
	    \global\csname halfheight#1\endcsname= \halfheight
	}
    }

\newlength{\starsize}
\def\starzero #1;{{
	\newpoint center
	\newpoint foo
	\newpoint trigpoint
	\centerx= \csname objectx#1\endcsname
	\centery= \csname objecty#1\endcsname
	\copypoint center to foo
	\trigoffset foo by {\starsize} along 0:1
	\put(\foox,\fooy){\makebox(0,0){$\bullet$}}
	\copypoint center to foo
	\trigoffset foo by {\starsize} along 0:-1
	\put(\foox,\fooy){\makebox(0,0){$\bullet$}}
	\copypoint center to foo
	\trigoffset foo by {\starsize} along 1:0
	\put(\foox,\fooy){\makebox(0,0){$\bullet$}}
	\copypoint center to foo
	\trigoffset foo by {\starsize} along -1:0
	\put(\foox,\fooy){\makebox(0,0){$\bullet$}}
    }}

\def\starone #1;{{
	\newpoint center
	\newpoint foo
	\newpoint trigpoint
	\centerx= \csname objectx#1\endcsname
	\centery= \csname objecty#1\endcsname
	\copypoint center to foo
	\trigoffset foo by {\starsize} along 1:1
	\put(\foox,\fooy){\makebox(0,0){$\bullet$}}
	\copypoint center to foo
	\trigoffset foo by {\starsize} along 1:-1
	\put(\foox,\fooy){\makebox(0,0){$\bullet$}}
	\copypoint center to foo
	\trigoffset foo by {\starsize} along -1:1
	\put(\foox,\fooy){\makebox(0,0){$\bullet$}}
	\copypoint center to foo
	\trigoffset foo by {\starsize} along -1:-1
	\put(\foox,\fooy){\makebox(0,0){$\bullet$}}
    }}

\def\startwo #1;{{
	\newpoint center
	\newpoint foo
	\newpoint trigpoint
	\centerx= \csname objectx#1\endcsname
	\centery= \csname objecty#1\endcsname
	\copypoint center to foo
	\trigoffset foo by {\starsize} along 1:2
	\put(\foox,\fooy){\makebox(0,0){$\bullet$}}
	\copypoint center to foo
	\trigoffset foo by {\starsize} along 1:-2
	\put(\foox,\fooy){\makebox(0,0){$\bullet$}}
	\copypoint center to foo
	\trigoffset foo by {\starsize} along -1:2
	\put(\foox,\fooy){\makebox(0,0){$\bullet$}}
	\copypoint center to foo
	\trigoffset foo by {\starsize} along -1:-2
	\put(\foox,\fooy){\makebox(0,0){$\bullet$}}
	\copypoint center to foo
	\trigoffset foo by {\starsize} along 2:1
	\put(\foox,\fooy){\makebox(0,0){$\bullet$}}
	\copypoint center to foo
	\trigoffset foo by {\starsize} along -2:1
	\put(\foox,\fooy){\makebox(0,0){$\bullet$}}
	\copypoint center to foo
	\trigoffset foo by {\starsize} along 2:-1
	\put(\foox,\fooy){\makebox(0,0){$\bullet$}}
	\copypoint center to foo
	\trigoffset foo by {\starsize} along -2:-1
	\put(\foox,\fooy){\makebox(0,0){$\bullet$}}
    }}

\def\starthree #1;{{
	\newpoint center
	\newpoint foo
	\newpoint trigpoint
	\centerx= \csname objectx#1\endcsname
	\centery= \csname objecty#1\endcsname
	\copypoint center to foo
	\trigoffset foo by {\starsize} along 1:3
	\put(\foox,\fooy){\makebox(0,0){$\bullet$}}
	\copypoint center to foo
	\trigoffset foo by {\starsize} along 1:-3
	\put(\foox,\fooy){\makebox(0,0){$\bullet$}}
	\copypoint center to foo
	\trigoffset foo by {\starsize} along -1:3
	\put(\foox,\fooy){\makebox(0,0){$\bullet$}}
	\copypoint center to foo
	\trigoffset foo by {\starsize} along -1:-3
	\put(\foox,\fooy){\makebox(0,0){$\bullet$}}
	\copypoint center to foo
	\trigoffset foo by {\starsize} along 3:1
	\put(\foox,\fooy){\makebox(0,0){$\bullet$}}
	\copypoint center to foo
	\trigoffset foo by {\starsize} along -3:1
	\put(\foox,\fooy){\makebox(0,0){$\bullet$}}
	\copypoint center to foo
	\trigoffset foo by {\starsize} along 3:-1
	\put(\foox,\fooy){\makebox(0,0){$\bullet$}}
	\copypoint center to foo
	\trigoffset foo by {\starsize} along -3:-1
	\put(\foox,\fooy){\makebox(0,0){$\bullet$}}
	\copypoint center to foo
	\trigoffset foo by {\starsize} along 2:3
	\put(\foox,\fooy){\makebox(0,0){$\bullet$}}
	\copypoint center to foo
	\trigoffset foo by {\starsize} along 2:-3
	\put(\foox,\fooy){\makebox(0,0){$\bullet$}}
	\copypoint center to foo
	\trigoffset foo by {\starsize} along -2:3
	\put(\foox,\fooy){\makebox(0,0){$\bullet$}}
	\copypoint center to foo
	\trigoffset foo by {\starsize} along -2:-3
	\put(\foox,\fooy){\makebox(0,0){$\bullet$}}
	\copypoint center to foo
	\trigoffset foo by {\starsize} along 3:2
	\put(\foox,\fooy){\makebox(0,0){$\bullet$}}
	\copypoint center to foo
	\trigoffset foo by {\starsize} along -3:2
	\put(\foox,\fooy){\makebox(0,0){$\bullet$}}
	\copypoint center to foo
	\trigoffset foo by {\starsize} along 3:-2
	\put(\foox,\fooy){\makebox(0,0){$\bullet$}}
	\copypoint center to foo
	\trigoffset foo by {\starsize} along -3:-2
	\put(\foox,\fooy){\makebox(0,0){$\bullet$}}
    }}

\def\starfour #1;{{
	\newpoint center
	\newpoint foo
	\newpoint trigpoint
	\centerx= \csname objectx#1\endcsname
	\centery= \csname objecty#1\endcsname
	\copypoint center to foo
	\trigoffset foo by {\starsize} along 1:4
	\put(\foox,\fooy){\makebox(0,0){$\bullet$}}
	\copypoint center to foo
	\trigoffset foo by {\starsize} along 1:-4
	\put(\foox,\fooy){\makebox(0,0){$\bullet$}}
	\copypoint center to foo
	\trigoffset foo by {\starsize} along -1:4
	\put(\foox,\fooy){\makebox(0,0){$\bullet$}}
	\copypoint center to foo
	\trigoffset foo by {\starsize} along -1:-4
	\put(\foox,\fooy){\makebox(0,0){$\bullet$}}
	\copypoint center to foo
	\trigoffset foo by {\starsize} along 4:1
	\put(\foox,\fooy){\makebox(0,0){$\bullet$}}
	\copypoint center to foo
	\trigoffset foo by {\starsize} along -4:1
	\put(\foox,\fooy){\makebox(0,0){$\bullet$}}
	\copypoint center to foo
	\trigoffset foo by {\starsize} along 4:-1
	\put(\foox,\fooy){\makebox(0,0){$\bullet$}}
	\copypoint center to foo
	\trigoffset foo by {\starsize} along -4:-1
	\put(\foox,\fooy){\makebox(0,0){$\bullet$}}
	\copypoint center to foo
	\trigoffset foo by {\starsize} along 3:4
	\put(\foox,\fooy){\makebox(0,0){$\bullet$}}
	\copypoint center to foo
	\trigoffset foo by {\starsize} along 3:-4
	\put(\foox,\fooy){\makebox(0,0){$\bullet$}}
	\copypoint center to foo
	\trigoffset foo by {\starsize} along -3:4
	\put(\foox,\fooy){\makebox(0,0){$\bullet$}}
	\copypoint center to foo
	\trigoffset foo by {\starsize} along -3:-4
	\put(\foox,\fooy){\makebox(0,0){$\bullet$}}
	\copypoint center to foo
	\trigoffset foo by {\starsize} along 4:3
	\put(\foox,\fooy){\makebox(0,0){$\bullet$}}
	\copypoint center to foo
	\trigoffset foo by {\starsize} along -4:3
	\put(\foox,\fooy){\makebox(0,0){$\bullet$}}
	\copypoint center to foo
	\trigoffset foo by {\starsize} along 4:-3
	\put(\foox,\fooy){\makebox(0,0){$\bullet$}}
	\copypoint center to foo
	\trigoffset foo by {\starsize} along -4:-3
	\put(\foox,\fooy){\makebox(0,0){$\bullet$}}
    }}

\def\getobjectposition #1 into #2 {
	\csname #2x\endcsname = \csname objectx#1\endcsname
	\csname #2y\endcsname = \csname objecty#1\endcsname
    }

\def\getobjectsize #1 into #2 {
	\csname #2x\endcsname = \csname halfwidth#1\endcsname
	\csname #2y\endcsname = \csname halfheight#1\endcsname
    }

\def\labelbox(#1,#2)(#3,#4)#5{{
	\newlocalbox\labelbox
	\savebox{\labelbox}{${#5}$}

	\newpoint {}
	\setpoint {} to (#1,#2)

	\newlocalcount\margin
	\margin= \arrowmargin
	\divide\margin by \unitlength
	
	\newlocalcount\height
	\newlocalcount\width
	\newlocalcount\depth
	\height= \ht\labelbox
	\width=  \wd\labelbox
	\depth=  \dp\labelbox
	\divide\height by \unitlength
	\divide\width  by \unitlength
	\divide\depth  by \unitlength

	\ifnum#3=0
	    \ifflippinglabels
	    	\advance\x by \margin
	      \else
	    	\advance\x by -\width
	    	\advance\x by -\margin
	      \fi
	     \advance\y by \y
	     \advance\y by -\height
	     \advance\y by \depth
	     \divide\y by 2
	  \else
	    \newlocalcount\slopesign
	    \slopesign= #3
	    \multiply\slopesign by #4

	    \ifflippinglabels
		\advance\y by -\height
		\advance\y by -\margin
		\multiply\slopesign by -1
	      \else
		\advance\y by \depth
		\advance\y by \margin
	      \fi
	    
	    \ifnum\slopesign>0
		\advance\x by -\width
		\advance\x by -\margin
	      \else
		\ifnum\slopesign<0
		    \advance\x by \margin
		  \else
		    \divide\width by 2
		    \advance\x by -\width
		  \fi
	      \fi
	  \fi

	\put(\x,\y){\usebox{\labelbox}}

	\newlocalcount\top
	\newlocalcount\bottom
	\top=    \y
	\bottom= \y
	\advance\top    by \height
	\advance\bottom by -\depth
	\expandybounds(\top,\bottom)
	\global\flippinglabelsfalse
    }}

\def\avoidobject #1 (#2,#3){{
	\getobjectsize {#1} into delta

	\newlocalcount\hweight
	\newlocalcount\vweight

	\hweight=\deltay
	\vweight=\deltax
	\multiply\hweight by #2
	\multiply\vweight by #3
	\makepositive(hweight)
	\makepositive(vweight)

	\ifnum\hweight<\vweight
		\ifnum#3<0 \deltay= -\deltay \fi
		\deltax=\deltay
		\multiply\deltax by #2
		\divide\deltax by #3
	    \else
		\ifnum#2<0 \deltax= -\deltax \fi
		\deltay=\deltax
		\multiply\deltay by #3
		\divide\deltay by #2
	    \fi

	\global\deltax= \deltax
	\global\deltay= \deltay
    }}

\def\calcarrow #1 #2 {
	\newpoint start
	\newpoint end
	\getobjectposition {#1} into start
	\getobjectposition {#2} into end

	\newvector arrow
	\setvector arrow from start to end

	\newpoint delta
	\avoidobject #1 (\arrowslopex,\arrowslopey)
	\shortentail arrow by delta

	\avoidobject #2 (\arrowslopex,\arrowslopey)
	\shortenhead arrow by delta

	\addpoint nudge to arrowhead
	\addpoint nudge to arrowtail
	
	\global\nudgex=0
	\global\nudgey=0
    }

\def\checkarrowslope #1 from #2 to #3 {{
	\ifnum\arrowslopex>4
	    \errmessage{Bad slope for arrow <#1> from #2 to #3:
			(\number\arrowslopex,\number\arrowslopey)}
	  \fi
	\ifnum\arrowslopex<-4
	    \errmessage{Bad slope for arrow <#1> from #2 to #3:
			(\number\arrowslopex,\number\arrowslopey)}
	  \fi
	\ifnum\arrowslopey>4
	    \errmessage{Bad slope for arrow <#1> from #2 to #3:
			(\number\arrowslopex,\number\arrowslopey)}
	  \fi
	\ifnum\arrowslopey<-4
	    \errmessage{Bad slope for arrow <#1> from #2 to #3:
			(\number\arrowslopex,\number\arrowslopey)}
	  \fi
    }}

\def\checksegmentslope #1 from #2 to #3 {{
	\ifnum\arrowslopex>6
	    \errmessage{Bad slope for segment <#1> from #2 to #3:
			(\number\arrowslopex,\number\arrowslopey)}
	  \fi
	\ifnum\arrowslopex<-6
	    \errmessage{Bad slope for segment <#1> from #2 to #3:
			(\number\arrowslopex,\number\arrowslopey)}
	  \fi
	\ifnum\arrowslopey>6
	    \errmessage{Bad slope for segment <#1> from #2 to #3:
			(\number\arrowslopex,\number\arrowslopey)}
	  \fi
	\ifnum\arrowslopey<-6
	    \errmessage{Bad slope for segment <#1> from #2 to #3:
			(\number\arrowslopex,\number\arrowslopey)}
	  \fi
    }}

\def\arrowlabel #1{{
	 \newpoint mid
	 \getvectormidpoint arrow into mid
	 \labelbox(\midx,\midy)(\arrowslopex,\arrowslopey){#1}
    }}

\def\arrow #1 #2 #3;{{
	\calcarrow {#1} {#2}
	\checkarrowslope {#3} from #1 to #2
	\drawvector arrow
	\arrowlabel {#3}
    }}

\def\segment #1 #2 #3;{{
	\calcarrow {#1} {#2}
	\checksegmentslope {#3} from #1 to #2
	\drawline arrow
	\arrowlabel {#3}
    }}

\def\dblarrow #1 #2 #3;{{
	\calcarrow {#1} {#2}
	\checkarrowslope {#3} from #1 to #2
	\drawdoublevector arrow
	\arrowlabel {#3}
    }}

\def\rightcell #1,#2 #3;{{
	\newpoint {}
	\copypoint coordinatescale to {}
	\multiply\x by #1
	\multiply\y by #2
	\put(\x,\y){\makebox(0,0){${\displaystyle{#3}}\atop\Longrightarrow$}}
    }}

\def\leftcell #1,#2 #3;{{
	\newpoint {}
	\copypoint coordinatescale to {}
	\multiply\x by #1
	\multiply\y by #2
	\put(\x,\y){\makebox(0,0){${\displaystyle{#3}}\atop\Longleftarrow$}}
    }}

\def\isoarrow #1 #2;{{
	\calcarrow {#1} {#2}
	\checkarrowslope {isomorphism} from #1 to #2
	\drawdoublevector arrow
	\ifnum0=\arrowslopex
	    \arrowlabel {\wr}
	  \else
	    \arrowlabel {\sim}
	  \fi
    }}

\def\calcrectlength {
	\newlocalcount\cornerlength
	\cornerlength= \cornersize
	\divide\cornerlength by \unitlength
    }

\def\calcdiaglength {
	\calcrectlength
	\multiply\cornerlength by 1000
	\divide\cornerlength by 1414
    }

\def\westcorner #1;{{
	\calcdiaglength

	\newpoint {}
	\getobjectposition {#1} into {}
	\advance\x by -\cornerlength
	\advance\x by -\cornerlength

	\addpoint nudge to {}
	\setpoint nudge to (0,0)

	\put(\x,\y){\line(1, 1){\cornerlength}}
	\put(\x,\y){\line(1,-1){\cornerlength}}
    }}

\def\eastcorner #1;{{
	\calcdiaglength

	\newpoint {}
	\getobjectposition {#1} into {}
	\advance\x by \cornerlength
	\advance\x by \cornerlength

	\addpoint nudge to {}
	\setpoint nudge to (0,0)

	\put(\x,\y){\line(-1, 1){\cornerlength}}
	\put(\x,\y){\line(-1,-1){\cornerlength}}
    }}

\def\northcorner #1;{{
	\calcdiaglength

	\newpoint {}
	\getobjectposition {#1} into {}
	\advance\y by \cornerlength
	\advance\y by \cornerlength

	\addpoint nudge to {}
	\setpoint nudge to (0,0)

	\put(\x,\y){\line( 1,-1){\cornerlength}}
	\put(\x,\y){\line(-1,-1){\cornerlength}}
    }}

\def\southcorner #1;{{
	\calcdiaglength

	\newpoint {}
	\getobjectposition {#1} into {}
	\advance\y by -\cornerlength
	\advance\y by -\cornerlength

	\addpoint nudge to {}
	\setpoint nudge to (0,0)

	\put(\x,\y){\line( 1, 1){\cornerlength}}
	\put(\x,\y){\line(-1, 1){\cornerlength}}
    }}

\def\northeastcorner #1;{{
	\calcrectlength

	\newpoint {}
	\getobjectposition {#1} into {}
	\advance\x by \cornerlength
	\advance\y by \cornerlength

	\addpoint nudge to {}
	\setpoint nudge to (0,0)

	\put(\x,\y){\line( 0,-1){\cornerlength}}
	\put(\x,\y){\line(-1, 0){\cornerlength}}
    }}

\def\northwestcorner #1;{{
	\calcrectlength

	\newpoint {}
	\getobjectposition {#1} into {}
	\advance\x by -\cornerlength
	\advance\y by \cornerlength

	\addpoint nudge to {}
	\setpoint nudge to (0,0)

	\put(\x,\y){\line( 0,-1){\cornerlength}}
	\put(\x,\y){\line( 1, 0){\cornerlength}}
    }}

\def\southeastcorner #1;{{
	\calcrectlength

	\newpoint {}
	\getobjectposition {#1} into {}
	\advance\x by \cornerlength
	\advance\y by -\cornerlength

	\addpoint nudge to {}
	\setpoint nudge to (0,0)

	\put(\x,\y){\line( 0, 1){\cornerlength}}
	\put(\x,\y){\line(-1, 0){\cornerlength}}
    }}

\def\southwestcorner #1;{{
	\calcrectlength

	\newpoint {}
	\getobjectposition {#1} into {}
	\advance\x by -\cornerlength
	\advance\y by -\cornerlength

	\addpoint nudge to {}
	\setpoint nudge to (0,0)

	\put(\x,\y){\line( 0, 1){\cornerlength}}
	\put(\x,\y){\line( 1, 0){\cornerlength}}
    }}

\long\def\fig#1{{
\newlocalcount\leastobjectx
\newlocalcount\leastobjecty
\newlocalcount\greatestobjectx
\newlocalcount\greatestobjecty
\leastobjectx= 2147483647
\leastobjecty= 2147483647
\greatestobjectx= -2147483647
\greatestobjecty= -2147483647
\newlocalbox\thefigurebox
\savebox{\thefigurebox}{#1}
\advance\greatestobjectx by -\leastobjectx
\advance\greatestobjecty by -\leastobjecty
\begin{center}
\begin{picture}(\greatestobjectx,\greatestobjecty)
	       (\leastobjectx,\leastobjecty)
\put(0,0){\usebox{\thefigurebox}}
\end{picture}
\end{center}
}}

\newtheorem{theorem}{Theorem}[section]
\newtheorem{prop}[theorem]{Proposition}
\newtheorem{lemma}[theorem]{Lemma}

\newtheorem{cor}[theorem]{Corollary}
\def\example{\par\medskip{\bf Example}\qua}
\def\examples{\par\medskip{\bf Examples}\qua}

\def\dfinition{\par\medskip{\bf Definition}\qua}
\def\remark{\par\medskip{\bf Remark}\qua}

\let\eproof\endproof
\def\eps{\varepsilon}
\def\Diff{{\rm Diff}}
\def\Aff{{\rm Aff}}
\def\Alpha{A}
\def\BS{{\rm BS}}

\newcommand\Q{{\mathbb Q}}
\newcommand\Z{{\mathbb Z}} 
\newcommand\R{{\mathbb R}} 
 
\newcommand\N{{\mathbb N}} 

\newcommand\C{{\mathbb C}}

\newcommand\RP{{\R P^1}}
\newcommand\CP{{\C P^1}}

\newcommand\s{{\bf s}}

\begin{document}
\title{Global rigidity of solvable group actions on $S^1$}
\asciititle{Global rigidity of solvable group actions on S^1}
\author{Lizzie Burslem\\Amie Wilkinson}

\address{Department of Mathematics, University of Michigan\\
2074 East Hall, Ann Arbor, MI 48109-1109 USA}

\secondaddress{Department of Mathematics, Northwestern University\\
2033 Sheridan Road, Evanston, IL 60208-2730 USA}

\asciiaddress{Department of Mathematics, University of Michigan
2074 East Hall, Ann Arbor, MI 48109-1109 USA
and
Department of Mathematics, Northwestern University
2033 Sheridan Road, Evanston, IL 60208-2730 USA}

\gtemail{\mailto{burslem@umich.edu}, \mailto{wilkinso@math.northwestern.edu}}
\asciiemail{burslem@umich.edu, wilkinso@math.northwestern.edu}

\begin{abstract}
In this paper we find all solvable subgroups of $\Diff^\omega(S^1)$
and classify their actions.  We also investigate the $C^r$ local
rigidity of actions of the solvable Baumslag--Solitar groups on the
circle.

The investigation leads to two novel phenomena in the study of 
infinite group actions on compact manifolds.  We exhibit
a finitely generated group $\Gamma$ and a manifold $M$ such that:
\begin{itemize}
\item $\Gamma$ has exactly countably infinitely many 
effective real-analytic 
actions on $M$, up to conjugacy in $\Diff^\omega(M)$;
\item every effective, real analytic action of $\Gamma$ on $M$ 
is $C^r$ locally rigid, for some $r\geq 3$, and for every
such $r$, there are infinitely many nonconjugate, effective real-analytic
actions of $\Gamma$ on $M$ that are $C^r$ locally rigid, but not
$C^{r-1}$ locally rigid.
\end{itemize}
\end{abstract}

\asciiabstract{%
In this paper we find all solvable subgroups of Diff^omega(S^1) and
classify their actions.  We also investigate the C^r local rigidity of
actions of the solvable Baumslag-Solitar groups on the circle.

The investigation leads to two novel phenomena in the study of infinite
group actions on compact manifolds. We exhibit a finitely generated
group Gamma and a manifold M such that: * Gamma has exactly countably
infinitely many effective real-analytic actions on M, up to conjugacy
in Diff^omega(M); * every effective, real analytic action of Gamma on
M is C^r locally rigid, for some r>=3, and for every such r, there are
infinitely many nonconjugate, effective real-analytic actions of Gamma
on M that are C^r locally rigid, but not C^(r-1) locally rigid.}

\asciikeywords{Group action, solvable group, rigidity, Diff^omega(S^1)}



\maketitlepage

\addcontentsline{toc}{section}{Introduction}
\section*{Introduction}

This paper describes two novel phenomena in the study of 
infinite group actions on compact manifolds.  We exhibit
a finitely generated group $\Gamma$ and a manifold $M$ such that:
\begin{itemize}
\item $\Gamma$ has exactly countably infinitely many 
effective real-analytic 
actions on $M$, up to conjugacy in $\Diff^\omega(M)$;
\item every effective, real analytic action of $\Gamma$ on $M$ 
is $C^r$ locally rigid, for some $r\geq 3$, and for every
such $r$, there are infinitely many nonconjugate, effective real-analytic
actions of $\Gamma$ on $M$ that are $C^r$ locally rigid, but not
$C^{r-1}$ locally rigid.
\end{itemize}
In the cases we know of where an infinite group $\Gamma$ has exactly countably
many smooth effective actions on a manifold $M$, that countable number is
finite, and indeed usually $0$. While many actions 
have been shown to
to be $C^r$ locally rigid, in the cases where a precise cutoff 
in rigidity has been established, it 
occurs between $r=1$ and $r=2$.  
For a survey of some of the recent results on smooth group actions,
see the paper of Labourie \cite{La}.

Our manifold $M$ is the circle $S^1$ and our group $\Gamma$
is the solvable Baumslag--Solitar group:
$$\BS(1,n) = \langle a,b,\,\mid\, aba^{-1} = b^n\rangle ,$$
where $n\geq 2$.

As a natural by-product of our techniques, we obtain a classification of
all solvable subgroups of $\Diff^\omega(S^1)$.  We show that every such subgroup 
${\mathcal G}$ is either 
%{\em virtually abelian}, meaning that $G$
%for some $m \in \Z$,
%the group ${\mathcal G}^m := \{g^m : g \in {\mathcal G} \}$ is abelian, or
conjugate in $\Diff^\omega(S^1)$
to a subgroup of a ramified affine group $\Aff^{\bf s}(\R)$, or, 
for
some $m \in \Z$, the group 
${\mathcal G}^m := \{g^m : g \in {\mathcal G} \}$ is abelian.
The ramified affine groups are defined and their properties
discussed in Section~\ref{s=ramaff}. Each 
ramified affine group is abstractly isomorphic to a direct
product $\Aff_+(\R)\times H$, where $\Aff_+(\R)$ is the group of 
orientation-preserving affine transformations of $\R$, and
$H$ is a subgroup of a finite dihedral group.

\section{Statement of results}
\subsection{Notation and preliminary definitions}
In places, we shall use two different
analytic coordinatizations of the circle $S^1$.  To denote an element
of the additive group, $\R/\Z$,
we will use $u$, and for an element of
the real projective line $\RP$, we will use  
$x$.  These coordinate systems are identified by:
$u\in \R/\Z \mapsto x = \tan(\pi u)\in \RP$.
When we are not specifying a coordinate system, we will use $p$ or $q$
to denote an element of $S^1$.  We fix an orientation on $S^1$
and  use ``$<$'' to denote the counterclockwise cyclic ordering
on $S^1$.

If $G$ is a group, then we denote by ${\mathcal R}^r(G)$ 
the set of all representations 
$\rho_0\co G \to \Diff^r(S^1)$,
and we denote by ${\mathcal R}^r_+(G)$ the set of all orientation-preserving 
representations in ${\mathcal R}^r(G)$.
Two representations $\rho_1, \rho_2\in{\mathcal R}^r(G)$ are
{\em conjugate (in $\Diff^r(S^1)$)} 
if there exists $h\in \Diff^r(S^1)$  such that,
for every $\gamma\in G$, $h\rho_1(\gamma) h^{-1} = \rho_2(\gamma)$.
%(equivalently, if this holds on a generating set for $G$).

We use the standard $C^k$ topology on representations of a
finitely-generated group into $\Diff^r(S^1)$, 
$r\in \{1,\ldots,\infty, \omega\}$ and $k\leq r$.
If $\Gamma$ is a finitely-generated group, then the $C^k$--open
sets in ${\mathcal R}^r(\Gamma)$  take generators
in a fixed generating 
set for $\Gamma$ into $C^k$--open sets.
A representation $\rho_0\in{\mathcal R}^r(\Gamma)$ is {\em ($C^r$) locally
rigid} if there exists a $C^1$ neighborhood ${\mathcal U}$ of $\rho_0$ in
${\mathcal R}^r(\Gamma)$
such that every $\rho\in {\mathcal U}$ is conjugate in
$\Diff^r(S^1)$ to $\rho_0$.
Finally, we say that $\Gamma$ is {\em globally rigid in $\Diff^r(S^1)$} 
if there
exists a countable set of locally rigid representations in
${\mathcal R}^r(\Gamma)$ such that every faithful representation
in  ${\mathcal R}^r(\Gamma)$ is conjugate to an element of this set.

To construct the subgroups and representations in this paper, we use
a procedure we call real ramified lifting.  

\dfinition A real analytic surjection $\pi\co S^1\to S^1$ is called
a {\em ramified covering map over $p \in S^1$} if
the restriction of $\pi$ to $\pi^{-1}(S^1\setminus\{p\})$ is a regular
analytic covering map  onto $S^1\setminus\{p\}$ of degree $d \geq 1$.
The {\em degree} of of $\pi$ is defined to be this integer $d$.

\medskip

Examples and properties of ramified covering maps and ramified lifts 
are described in Section~\ref{s=ramlift}. 

Let  $\pi\co S^1\to S^1$ be a ramified covering map over $p\in S^1$, and let
$f\co S^1\to S^1$ be a real analytic diffeomorphism that fixes $p$.  
We say that $\hat{f}\in \Diff^\omega(S^1)$ is a {\em $\pi$--ramified lift} 
of $f$ if the following diagram commutes:

\fig{
\object {a}  0,4 {S^1};
\object {b}  8,4 {S^1};
\object {c}  0,0 {S^1};
\object {d}  8,0  {S^1};
\arrow  {a}  {b}  {\hat f};
\arrow  {c}  {d}  {f};
\arrow  {a}  {c}  {\pi};
\arrow  {b}  {d}  {\pi};}

More generally, let $\rho\co  \Gamma \to \Diff^\omega(S^1)$ be a representation with global
fixed point $p$.  A representation $\hat\rho \co  \Gamma \to \Diff^\omega(S^1)$
is a {\em $\pi$--ramified lift of $\rho$}   if $\hat\rho(\gamma)$ is
a $\pi$--ramified lift $\rho(\gamma)$, for every $\gamma\in\Gamma$.  We will show in
Proposition~\ref{p=lift} that
a representation can have more than one $\pi$--ramified lift.

For $G$  a subgroup
of $\Diff^\omega(S^1)$ with a global fixed point $p$,
we define ${\hat G}^\pi$, the
{\em $\pi$--ramified lift of
$G$} to be the collection of all $\pi$--ramified lifts of elements of $G$.
By Proposition~\ref{p=lift} and Propostion~\ref{p=extension}, ${\hat G}^\pi$ 
is a subgroup of $\Diff^\omega(S^1)$,
abstractly isomorphic to an $H$--extension of $G_+$, where
$G_+ = \Diff_+^\omega(S^1)\cap G$, and $H$ is a subgroup of
a dihedral group determined by $\pi$.

\subsection{Rigidity of solvable Baumslag--Solitar groups}

In real projective coordinates on $\RP$, 
the {\em standard representation} $\rho_n$ of $\BS(1,n)$ into
$\Diff_+^\omega(S^1)$
takes the generators $a$ and $b$ to the affine maps
$$x\mapsto nx,\qquad\hbox{and}\qquad x\mapsto x+1.$$
This representation has a global fixed point $\infty\in\RP$. Our
first result states that $\BS(1,n)$ is globally rigid in $\Diff^\omega(S^1)$:

\begin{theorem}
\label{t=main}
For each $n\geq 2$, there are exactly countably infinitely many
faithful representations of $\BS(1,n)$ into $\Diff^\omega(S^1)$, up
to conjugacy in $\Diff^\omega(S^1)$.  Each conjugacy class
contains a $\pi$--ramified lift of $\rho_n$, where $\pi\co \RP\to\RP$
is a rational map that is ramified over $\infty$.
Furthermore, if $\rho\co \BS(1,n)\to \Diff^\omega(S^1)$ is not faithful,
then
there exists a $k\geq 1$ such that $\rho(b)^k=id$.
\end{theorem}
 
We give an explicit description of these conjugacy classes in
Section~\ref{s=bslifts}.

The conclusion of Theorem~\ref{t=main} 
does not hold when $C^\omega$ is replaced 
by a lower differentiability class such as $C^\infty$, 
even when analytic
conjugacy is replaced by topological conjugacy in the statement.
Nonetheless, as $r$ increases, there is a sort of ``quantum rigidity''
phenomenon.
Let $\rho\co \BS(1,n) \to \Diff^2(S^1)$ 
be a representation,
and let $f = \rho(a)$. We make a preliminary observation:

\begin{lemma} If the rotation number of $f$ is irrational,
then $g^k = id$, for some $k\leq n+1$,
where $g=\rho(b)$.
\end{lemma}
(See the beginning of Section~\ref{p=mainproofs} for a proof).
Hence, if $\rho\in {\mathcal R}^2(\BS(1,n))$ is faithful,
then $f$ must have periodic points.
For $\rho\in {\mathcal R}^2(\BS(1,n))$ a faithful representation, 
we define the {\em inner spectral radius}
$\sigma(\rho)$  by:
$$\sigma(\rho) = \sup\{|(f^k)^\prime(p)|^{\frac1k}\,\mid\, 
p\in {\rm Fix}(f^k)\hbox{ and } 
|(f^k)^\prime(p)|\leq  1\}.$$
For the standard representation, $\sigma(\rho_n) = \frac{1}{n}$, 
and if $\hat\rho_n$ is a ramified lift of $\rho_n$, then
$\sigma(\hat\rho_n)= \left(\frac1n\right)^{\frac1s}$,
for some $s\in \N_{\geq 1}$.

\begin{theorem}\label{t=maincr} Let 
$\rho\co  \BS(1,n) \to \Diff^{r}(S^1)$ be a faithful
representation, where $r\in [2,\infty]$.  
If either $r < \infty$ 
and $\sigma(\rho) \leq \left(\frac{1}{n}\right)^{\frac{1}{r-1}}$, 
or $r = \infty$ and $\sigma(\rho) < 1$,
then $\rho$ is conjugated by an element of 
$\Diff^{r}(S^1)$ into a unique conjugacy class in 
${\mathcal R}^\omega(\BS(1,n))$.

If $\rho$ takes values in $\Diff^{r}_+(S^1)$,
then  $\rho$ is conjugated by an element of $\Diff^{r}_+(S^1)$
into a unique conjugacy class in 
${\mathcal R}^\omega_+(\BS(1,n))$.
\end{theorem}

Theorem~\ref{t=maincr} has the following corollary:

\begin{cor}  Every  representation $\rho\co \BS(1,n)\to \Diff^\omega(S^1)$
is $C^\infty$ locally rigid.   Further, if
$\sigma(\rho) < \left(\frac{1}{n}\right)^{\frac{1}{r-1}}$,
then $\rho$ is $C^{r}$ locally rigid. 
\end{cor}

This corollary implies that the standard representation is $C^r$
locally rigid, for all $r\geq 3$, 
and every representation in ${\mathcal R}^\omega(\BS(1,n))$ is
locally rigid in some finite differentiability classes.
This local rigidity breaks down, however, if the differentiability class is
lowered.

\begin{prop}\label{p=deform}
For every representation $\rho\co \BS(1,n)\to \Diff^\omega(S^1)$,
if $\sigma(\rho) =  \left(\frac{1}{n}\right)^{\frac{1}{r-1}}$, for
some $r\geq 2$ then
there exists a
family of representations $\rho_t\in {\mathcal R}^{r}(\BS(1,n))$,
$t\in(-1,1)$,
with the following properties:
\begin{enumerate}
\item $\rho_0 = \rho$,
\item $t\mapsto \rho_t$ is continuous in the $C^{r-1}$ topology on ${\mathcal R}^{r}(\BS(1,n))$,
\item for every  $t_1, t_2\in (-1,1)$, if $\rho_{t_1}$ is conjugate
to $\rho_{t_2}$ in $\Diff^1(S^1)$ then $t_1 = t_2$.
\end{enumerate}
\end{prop}

It follows from our characterization of the conjugacy classes
in ${\mathcal R}^\omega (\BS(1,n))$ in the
next section that, for each value of $r\in [1,\infty)$, 
there are infinitely many nonconjugate representations
$\rho\in {\mathcal R}^\omega(\BS(1,n))$  satisfying  
$\sigma(\rho) =  \left(\frac{1}{n}\right)^{\frac{1}{r}}.$
Hence, for each $r\geq 3$ 
there are infinitely many distinct (nonconjugate) representations 
in ${\mathcal R}^\omega(\BS(1,n))$ that are $C^r$ locally rigid, but not 
$C^{r-1}$ locally rigid.

A. Navas has given a complete classification of
$C^2$ solvable group actions, up to finite index subgroups
and topological semiconjugacy.  
%As in \cite{FaFr}, a main tool in this classification is  Kopell's
%lemma (Lemma~\ref{l=kopell}), which classifies the centralizer of a $C^2$
%contraction on $[0,\infty)$.
One corollary of his result is that every
faithful $C^2$ representation $\rho$ of $\BS(1,n)$ into $\Diff^2(S^1)$ is
is virtually topologically semiconjugate to the standard representation:

\begin{theorem}\label{t=navas}{\rm\cite{Nav}}\qua
Let $\rho\co \BS(1,n)\to \Diff^r(S^1)$ be a representation, where
$r\geq 2$.  Then either
$\rho$ is  unfaithful, in which case $\rho(b)^m = id$, for some $m$, or
there exists an integer $m\geq 1$, 
a finite collection of closed, connected sets $I_1,\ldots,
I_k$, and a surjective continuous map 
$\varphi\co  S^1 \to \RP$ with the following properties:
\begin{enumerate}
\item $\rho(b)^m$ is the identity on each set $I_k$;
\item $\varphi$ sends each set $I_k$ to $\infty$;
\item the restriction of $\varphi$ to $S^1\setminus \bigcup_{i=1}^k I_i$
is a $C^r$ covering map of $\R$;
\item For every $\gamma\in \BS(1,n)$, the following diagram commutes:

\fig{
\object {a}  0,4 {S^1};
\object {b}  9,4 {S^1};
\object {c}  0,0 {S^1};
\object {d}  9,0  {S^1};
\arrow  {a}  {b}  {\rho(\gamma^m)};
\arrow  {c}  {d}  {\rho_n(\gamma^m)};
\arrow  {a}  {c}  {\varphi};
\arrow  {b}  {d}  {\varphi};}
where $\rho_n\co \BS(1,n)\to \Diff^\omega(S^1)$ is the standard representation.
\end{enumerate}
\end{theorem}

The map $\varphi$ in Theorem~\ref{t=navas} is a sort of ``broken $C^r$
ramified cover.''  The regularity of $\varphi$ at the preimages of the
point $\infty$ can be poor, and the map can be infinite-to-one on the
sets $I_1,\ldots, I_k$, but a map with these features is nothing more
than a deformation of a ramified covering map. 
Combining Theorem~\ref{t=navas} with Theorem~\ref{t=maincr} and
the proof Proposition~\ref{p=deform}, we obtain:

\begin{cor} Let $\rho\co \BS(1,n)\to \Diff^{r}(S^1)$,  be any
representation, with $r\geq 2$.  Then either:
\begin{enumerate}
\item $\rho$ is not faithful, and there exists an $m\geq 1$ such that
$\rho(b)^m=id$;
\item  $\rho$ admits $C^{r-1}$
deformations as in Proposition~\ref{p=deform}; or
\item $\rho$ is $C^r$ conjugated into a unique conjugacy class in
${\mathcal R}^\omega(\BS(1,n))$.
\end{enumerate}
\end{cor}

Since the statement of Theorem~\ref{t=navas} does not appear explicitly in 
Navas's paper, and we don't use this result elsewhere in the paper,
we sketch the proof at the end of Section~\ref{p=mainproofs}.

Finally, note that the 
trivial representation $\rho_0(a) = \rho_0(b) = id$ is not rigid in any
topology;  it can
be approximated  by the representation $\rho(b) = id, \rho(a)= f$, 
where $f$ is an any diffeomorphism close to the identity. 
Another nice consequence of Navas's theorem is that this is the
{\em only} way to $C^2$ deform the
trivial representation.

\begin{cor}\label{t=deformtriv} There is a 
$C^2$ neighborhood ${\mathcal U}\subset {\mathcal R}^2(\BS(1,n))$
of the trivial representation such that, for all $\rho\in {\mathcal U}$,
$\rho(b) = id$.
\end{cor}

\proof Let $\rho$ be a $C^2$ representation.  Since $\rho(b)$ 
is conjugate by $\rho(a)$ to $\rho(b)^n$, it will have rotation number
of the form $\frac{k}{n-1}$ if $\rho(a)$ is orientation-preserving, and
of the form $\frac{k}{n+1}$ if $\rho(a)$ is orientation-reversing.  Therefore,
if $\rho$ is sufficiently $C^0$--close to $\rho_0$ and if
$\rho(b)^m=id$, for some $m\geq 1$, then $m=1$.  So we may assume
that there exists a map $\varphi$ as in Theorem~\ref{t=navas} and
that $m=1$.  On a
component of $S^1\setminus \bigcup I_i$, $\varphi$ is a diffeomorphism
conjugating the action of $\rho$ to the restriction
of the standard representation $\rho_n$ to $\R$ (in general
$\varphi$ fails to extend to a diffeomorphism at either endpoint
of $\R$).
But in the standard action, the element $\rho_n(a)$ has a fixed point
in $\R$ of derivative $n$.  If $\rho$ is sufficiently $C^1$ close
to $\rho_0$, this can't happen.\eproof

We remark that, in contrast to the results in this paper, 
there are uncountably many topologically distinct 
faithful representations of $\BS(1,n)$ into
$\Diff^\omega(\R)$ (see \cite{FaFr}, Proposition 5.1).  The proof of our results uses the existence
of a global fixed point on $S^1$ for a finite index subgroup of
$\BS(1,n)$; such a point need not exist when $\BS(1,n)$ acts on 
${\bf\R}$.  Farb and Franks \cite{FaFr} studied actions of
Baumslag--Solitar groups on the line and circle.  Among their results,
they prove that if $m>1$, the (nonsolvable) Baumslag--Solitar group:
$$\BS(m,n) = \langle  a,b\,|\, ab^ma^{-1} = b^n\rangle ,$$
has no faithful $C^2$ actions on $S^1$ if $m$ does not
divide $n$.  They ask whether the actions of $B(1,n)$ on the circle can
be classified.  This question inspired the present paper.

\subsection{Classification of solvable subgroups of $\Diff^\omega(S^1)$}

Several works address the properties of solvable subgroups of 
$\Diff^r(S^1)$; we mention a few here.  Building on work
of Kopell \cite{Ko}, Plante and Thurston
\cite{PT} showed that any nilpotent subgroup of $\Diff^2(S^1)$
is in fact abelian.  Ghys \cite{Gh1} proved that
any solvable subgroup of $\Diff^\omega(S^1)$ is metabelian, ie,
two-step solvable.  In the same work, 
he remarks that there are
solvable subgroups of $\Diff^\infty(S^1)$ that are not metabelian.
The subgroups he constructs contain {\em infinitely flat elements} ---
nontrivial diffeomorphisms $g\in \Diff^\infty(S^1)$ with the property that
for some $p \in S^1$,
$g(p)=p$,
$g'(p)  =  1$, and $g^{(k)}(p) = 0$ for all $k \geq 2$.

Navas \cite{Nav} constructed further examples of solvable subgroups 
of $\Diff^\infty(S^1)$
with arbitrary degree of solvability, again using infinitely flat elements.
As mentioned above, Navas's work also contains a topological
classification of solvable subgroups 
of $\Diff^2(S^1)$.
As part of a study of ergodicity of actions of discrete groups on $S^1$, 
Rebelo and Silva \cite{RS} also study the solvable subgroups of
$\Diff^\omega(S^1)$.

Our main result in this part of the paper, Theorem~\ref{t=solv}, implies
that any solvable 
subgroup of $\Diff^\infty(S^1)$ that does not contain infinitely
flat elements is either virtually abelian
or conjugate to a subgroup of a ramified lift of the affine group:
$${\rm Aff}(\R) =  \{x\mapsto c x + d \co  \RP\to \RP\,\mid\, c,d\in \R,\,
c\neq 0\}.$$

\begin{theorem}\label{t=solv} Let $G < \Diff^r(S^1)$ be a solvable group, where $r \in
\{\infty, \omega\}$.  Then either:
\begin{enumerate}
\item for some $m \in \Z$, the group $G^m: =\{g^m: g\in G\}$ is abelian,
\item $G$ contains infinitely flat elements (which can't happen if 
$r=\omega$), or 
\item $G$ is conjugate in $\Diff^r(S^1)$ to 
a  subgroup of a $\pi$--ramified lift of 
$\Aff(\R)$, where $\pi\co \RP\to \RP$ is a ramified cover over $\infty$. 
\end{enumerate}
Further, if $G < \Diff^r_+(S^1)$ and {\rm(3)} holds, then the conjugacy 
can be taken in $\Diff^r_+(S^1)$.
\end{theorem}

In Section~\ref{s=ramlift}, we characterize the ramified lifts
of $\Aff(\R)$.  To summarize the results there, we have:
 
\begin{theorem}\label{t=ramchar} There exists a collection
$${\mathcal{RAFF}}:= \{ \widehat{{\rm Aff}}^{\bf s}(\R)
<\Diff^\omega(S^1)\,\mid\, {\bf s}\in
\overline{\mathcal S}\},$$ 
where $\overline{\mathcal S}$ 
is a countably infinite index set, with the following
properties:
\begin{enumerate} 
\item if ${\bf s}_1,{\bf s}_2\in \overline{\mathcal S}$ and
$\widehat{{\rm Aff}}^{{\bf s}_1}(\R)$ is conjugate to 
$\widehat{{\rm Aff}}^{{\bf s}_2}(\R)$ in $\Diff^1(S^1)$, then
${\bf s}_1={\bf s}_2$;
\item  for each ${\bf s}\in \overline{\mathcal S}$, there exists a subgroup
$H$ of a dihedral group such that $\widehat{{\rm Aff}}^{\bf s}(\R)
\simeq {\rm Aff}_+(\R)\times H$,
\item for each finite dihedral or cyclic group $H$, there exist
infinitely many ${\bf s}\in\overline{\mathcal S}$ so that
$\widehat{{\rm Aff}}^{\bf s}(\R) \simeq {\rm Aff}_+(\R)\times H$,
\item each element of ${\mathcal{RAFF}}$ is the $\pi$--ramified lift of
$\Aff(\R)$, for some rational ramified cover $\pi\co \RP\to\RP$ over $\infty$, and
every $\pi$--ramified lift of $\Aff(\R)$
is conjugate in $\Diff^\omega(S^1)$ to an element of ${\mathcal{RAFF}}$.
\end{enumerate}

There also exists a collection 
$${\mathcal{RAFF}}_+:= \{ \widehat{{\rm Aff}}^{\bf s}_+(\R)
<\Diff_+^\omega(S^1)\,\mid\, {\bf s}\in
\overline{\mathcal S}_+\},$$
with the same properties, except that in {\rm(1)} and {\rm(4)}, the conjugacy is
orientation-preserving, and in {\rm(2)} and {\rm(3)},  $H$ is cyclic.
\end{theorem}

Hence we have found all solvable groups that act effectively
on the circle as real-analytic diffeomorphisms.
%, and we have found all of the actions of those that are not 
%virtually abelian.

\section{Introduction to ramified lifts}\label{s=ramlift}

Let $G$ be a group and let 
$\rho\co G\to {\rm Diff }^\omega(S^1)$ be a representation with
a global fixed point $p$.  Restricting each element of this
representation to a suitably small neighborhood of $p$,
we obtain a representation $\tilde{\rho}\co G\to {\mathcal G}^\omega$,
where ${\mathcal G}^\omega$ is the group of analytic germs of 
diffeomorphisms. It is known \cite{CM, Na, EISV} that if 
if $G$ is solvable, then for some $k\geq 1$,
$\tilde\rho$ is conjugate in ${\mathcal
G}^\omega$ to a representation taking values in the {\em ramified affine
group} ${\rm Aff }^k(\R)$:
$${\rm Aff }^k(\R) = \{\frac{x}{(ax^k + b)^{\frac{1}{k}}}\,\mid\, a,b\in
\R, b>0\}$$
(see \cite{Gh1} for a proof in
the context of circle diffeomorphisms).  The name ``ramified affine
group'' is explained by the fact that the elements of ${\rm Aff }^k(\R)$
are lifts of the elements of the affine group under the branched (or
ramified) cover $z\mapsto z^k$.  These lifts are well-defined as 
holomophic germs, but do not extend to diffeomorphisms of $\CP$.

The key observation of this paper is that the elements of ${\rm Aff}(\R)$
do admit global ramified lifts as diffeomorphisms of $\RP$.
The reason is that,
in contrast to a ramified cover $\pi\co \CP\to \CP$, which must be ramified over
$2$ points,  a ramified cover $\pi\co \RP\to \RP$ is ramified
over one point, which can be chosen to coincide with the global fixed
point of ${\rm Aff}(\R)$.  

\medskip

\noindent{\bf Examples of real ramified covers}\qua
The  map $\pi_1\co  \R/\Z \rightarrow \R/\Z$ given by
$\pi_1(u) = \sin^2(\pi u)$
is a ramified covering map over $0$, with critical points of order
$2$ at $\pi_1^{-1}(0) = \{ 0, \frac{1}{2} \}$.

The rational map $\pi_2\co \RP\to \RP$ given by:
$$\pi_2(x) = \frac{(x+1)^2(x-1)^2}{x (x^2+1)}$$
is also a ramified covering map over $0$, with critical points of order
$2$ at $\pm 1$. It is clear that $\pm 1$ are critical points of $\pi_2$,
and one verifies directly that the other critical points of
$\pi_2$ in $\CP$ occur off of $\RP$, at $\pm i \sqrt{3 \pm\sqrt{8}}$. 

We will define an equivalence relation on ramified covering maps in
which $\pi_1$ and $\pi_2$ are equivalent, and show
 that, under this notion of
equivalence, all possible ramified covering maps occur as rational maps.

If $\pi\co S^1\to S^1$ is a ramified
covering map over $p$ then 
for each $q\in \pi^{-1}(p)$, there exists an integer
$s(q)\geq 1$ such that the leading (nonconstant) term in the
Taylor expansion of $\pi$ at $q$ is of order $s(q)$.
A regular covering map is a ramified covering map; in this case, 
$d$ is the topological degree of the map, 
and $s(q)=1$, for each $q\in \pi^{-1}(p)$.  As the examples show,
a ramified covering map need not be a regular covering map 
(even topologically), as it is
possible to have $s(q) > 1$.

Let $\pi$ be a ramified covering map over $p$, and let 
$q_1, \ldots, q_d$ be the elements
of $\pi^{-1}(p)$, ordered so that
$p\leq q_1 < q_2  < \cdots < q_d < p$.
For each $i\in \{1,\ldots , d\}$ we define $o_i\in \{\pm 1\}$
by:
$$o_i =  \begin{cases}1 & \hbox{if }\pi\vert_{(q_i, q_{i+1})} \hbox{ is
orientation-preserving,}\cr
-1 &  \hbox{if }\pi\vert_{(q_i, q_{i+1})} \hbox{ is
orientation-reversing.}\end{cases}$$

We call the vector ${\bf s}(\pi) = (s(q_1), s(q_2), \ldots , s(q_d), 
o_1, \ldots o_d)\in
\N^d\times \{\pm 1\}^d$ the
{\em signature} of $\pi$. 
Geometrically, we think of a signature as a regular $d$--gon in $\R^2$
with vertices labelled by $s_1, \ldots, s_d$ and edges labelled
by $o_1, \ldots , o_d$.   
Every signature vector ${\bf s}=(s_1,\ldots, s_d, o_1, \ldots, o_d)$ 
has the following two properties:
\begin{enumerate}
\item The number of vertices with an even label is  even:
$\#\{ 1\leq i\leq d\,\mid\, s_i \in 2\N\}\in 2\N.$
\item If a vertex has an odd label, then both edges connected to that
vertex have the same label, and if a vertex has an even label, then the
edges have opposite labels:
$(-1)^{s_i +1} = o_{i-1}\, o_{i},$
where addition is$\mod d$.
\end{enumerate}
We will call any  vector ${\bf s}\in \N^d\times\{\pm 1\}^d$ 
with these properties
a {\em signature vector}.  Note that a signature vector of length $2d$ 
is determined by its first $d+1$ entries.
Let ${\mathcal S}_d$ be the set of all
signature vectors with length $2d$,
and let ${\mathcal S}$ be the set of
all signature vectors.

\begin{prop}\label{p=ramcovexists}
Given any ${\bf s}\in {\mathcal S}$ and $p\in S^1$, there
is a ramified covering map $\pi\co  S^1 \rightarrow S^1$ over $p$
with signature ${\bf s}$.
\end{prop}

\proof Let ${\bf s} = (s_1,\ldots, s_d, o_1, \ldots, o_d)$ be a signature, 
and let $p\in \R/\Z$.  Choose points
$u_1< \cdots < u_d$ evenly spaced in $\R/\Z$, and let
$F \co  \R/\Z \to \R/Z$  be the piecewise affine map
that sends the $u_i$ to $p$, and which sends each component of 
$\R/\Z\setminus \{u_1, \ldots, u_d\}$ onto
$\R/Z\setminus  \{p\}$, with orientation determined by $o_i$.

Put a new analytic structure on $\R/\Z$ as follows.
In the intervals $I_j =
(u_j,u_{j+1})$ use the standard analytic charts, but in each
interval  $J_j = (u_j-\epsilon,
u_j+\epsilon)$  compose the standard chart (that identifies $u_j$ with
$0$) with the homeomorphism
$\sigma_j\co  \R\to \R$ defined by:
$$\sigma_j(x) = \begin{cases} x^{1/s_j} & \hbox{if } x > 0\\
- |x|^{1/s_j} & \hbox{if } x \leq 0.
\end{cases}$$
Since the overlaps are analytic, this defines a real anaytic atlas on
$\R/\Z$.                                                                              

Note that the map $F\co  (\R/\Z, \hbox{new structure})\to (\R/\Z,
\hbox{standard structure})$ is analytic: in charts around $u_j$
and $p=F(u_j)$, 
the map $F$ takes the form $x\mapsto x^{s_j}$.  
Since there
is a unique real analytic structure on the circle, there is an analytic homeomorphism
of the circle $h\co (\R/\Z, \hbox{standard structure})\to (\R/\Z,
\hbox{new structure})$. Let $\pi = F\circ h$.  Then $\pi$ is a
a ramified covering map over $0$ with signature ${\bf s}$.\eproof

In fact, ramified covers exist in the purely algebraic
category; every signature can be realized by a rational map. We have:

\begin{prop}\label{p=ratramcovexists}
Given any ${\bf s}\in {\mathcal S}$ and $p\in \RP$, there
is a rational map $\pi\co  \RP\rightarrow \RP$ that is a ramified cover 
over $p$
with signature ${\bf s}$.
\end{prop}

\proof Since the proof of Proposition~\ref{p=ratramcovexists} is somewhat lengthy, we omit the details.
The construction proceeds as follows.  Let ${\bf s} = (s_1,\ldots, s_d, o_1, \ldots, o_d)$ be a signature, and assume that $p=0\in\RP$ and
$o_1=1$.  Choose a sequence of real numbers 
$a_0 < a_1 < \ldots < a_{2d-2}$, 
let $P(x) = (x-a_0)^{s_1} (x-a_2)^{s_2} \ldots (x-a_{2d-2})^{s_d}$
and let $Q(x) = (x-a_1) (x-a_3) \ldots (x-a_{2d-3})$.  The desired rational
function $\pi$ will be a modification of $P/Q$.

Let $h(x)$ be a polynomial of even degree with no zeros, with critical 
points of even degree at $a_i$, $0 \leq i \leq 2d-2$, and with no other
critical points. One first shows that, for $N$ sufficiently large, 
the rational function:
$${\pi_0} = \frac{P h^N}{Q}$$
has zeroes of order $s_1, \ldots, s_d$ at $a_0, a_2, \ldots, a_{2d-2}$,
simple poles at $a_1, a_3, \ldots,$ $a_{2d-3}$, a pole of odd order at $\infty$, and no other zeroes, poles or critical points.  Hence $\pi_0$ is
a ramified covering map over $0$ with signature ${\bf s}$, 
except at $\infty$, where it may fail to be a diffeomorphism.

Choose such an $N$, and let $2m+1$ be the order of the pole $\infty$ for 
$\pi_0$. One then shows 
that for $\eps$ sufficiently small, the rational function:
$$\pi(x) =  \frac{{\pi_0}(x)}{1 + \eps x^{2m}}$$
has the same properties as $\pi_0$, except that $\infty$ is now a simple pole; it
is the desired ramified cover.\eproof

If $\pi$ is a ramified covering map, then the cyclic and dihedral groups:
$$C_d =  \langle b: b^d = id \rangle, \; \; \hbox{and} \; \;
D_d = \langle a,b: b^d = id, a^2 = id, 
aba^{-1} = b^{-1} \rangle,$$
respectively, act on $\pi^{-1}(p)$ and on the set ${\mathcal E}(\pi)$
of oriented components
of $S^1\setminus \pi^{-1}(p)$ in a natural way.  By an
orientation-preserving
homeomorphism, we identify the circle with a regular oriented $d$--gon,
sending the elements of $\pi^{-1}(p)$ to the vertices and the elements
of ${\mathcal E}(\pi)$ to the edges.  The groups $C_d \lhd D_d$
act by symmetries of the $d$--gon, inducing actions on  $\pi^{-1}(p)$ and
${\mathcal E}(\pi)$ that are clearly independent of choice of
homeomorphism.  For $q\in \pi^{-1}(p)$, $e\in{\mathcal E}(\pi)$, and
$\zeta\in D_d$, we write $\zeta(q)$ and $\zeta(e)$ for 
their images under this action.

These symmetry groups also act on the signature vectors in ${\mathcal
S}_d$  in the 
natural way, permuting both vertex labels and edge labels.
For $\zeta\in D_d$, we will write
$\zeta({\bf s})$ for the image of  ${\bf s}\in{\mathcal S}_d$
under this action.
In this notation, the  action is generated by:
$$b(s_1,  \ldots , s_d, o_1, \ldots, o_d) = (s_2, s_3, \ldots , s_d,
s_1, o_2, o_3, \ldots, o_d, o_1),$$
and
$$a(s_1, \ldots , s_d, o_1, \ldots, o_d) = (s_1, s_d, s_{d-1}
\ldots , s_3
,s_2, -o_d, -o_{d-1}, \ldots, -o_2, -o_1).$$
Denote by ${\rm Stab}_{C_d}({\bf s})$ and
${\rm Stab}_{D_d}({\bf s})$ the
stabilizer of ${\bf s}$ in $C_d$ and $D_d$, respectively, under this
action:
$${\rm Stab}_H({\bf s}) = \{\zeta\in H\,\mid\, \zeta({\bf s}) = {\bf s}\},$$
for $H= C_d$ or $D_d$.

\examples The signature vector of $\pi_1(u)=\sin^2(\pi u)$
is ${\bf s}_1=(2,2,1,-1)$. The stabilizer of ${\bf s}_1$
in $D_d$ is ${\rm Stab}_{D_d}({\bf s}_1) = \langle a\rangle $, and the 
stabilizer of ${\bf s}_1$
in $C_d$ is trivial. The signature vector of $\pi_2(x) =
((x-1)^2(x+1)^2)/(x (x^2 + 1))$ is ${\bf s}_2 = (2,2,-1, 1)$.  Note that
${\bf s}_2$ lies in the $C_d$--orbit of ${\bf s}_1$, and so 
${\rm Stab}_{C_d}({\bf s}_2)$ and ${\rm Stab}_{D_d}({\bf s}_2)$ 
must be conjugate to ${\rm Stab}_{C_d}({\bf s}_1)$ and
${\rm Stab}_{D_d}({\bf s}_1)$, respectively, by an element of $C_d$.
In this simple case, the stabilizers are equal.

%The signature of the rational map $\pi\co \RP\to \RP$ given by:
%$$\pi(x) = \frac{x^2(x-2)^3(x-4)(x-6)^2(x-8)^3(x-10)}
%{(x-1)(x-3)(x-5)(x-7)(x-9)(1+\eps x^6)},$$
%for $\eps$ sufficiently small, is 
For another example, consider the signature vector 
$$(2,3,1,2,3,1,-1,-1, -1, 1, 1, 1),$$ which geometrically is 
represented by the following labelled graph:                                                                               
\fig{
\object {a}  4,3 {2};
\object {b}  8,6 {1};
\object {c}  8,10 {3};
\object {d}  4,13  {2};
\object {e}  0,10 {1};
\object {f}  0,6 {3};
\arrow  {b}  {a}  \fliplabel{-};
\arrow  {c}  {b}  \fliplabel{-};
\arrow  {d}  {c}  {-};
\arrow  {d}  {e}  {+};
\arrow  {e}  {f}  {+};
\arrow  {f}  {a}  \fliplabel{+};
}
\noindent
This labelling has no symmetries, despite the fact that the
edge labels have a flip symmetry and the vertex labels have a rotational
symmetry. By contrast, the signature
$(2, 1, 4, 2, 1, 4, -1, -1, 1, -1, -1, 1)$
has a $180$ degree rotational symmetry corresponding to the element 
$b^3\in C_6$, and so both stabilizer subgroups are $\langle b^3\rangle $.

\subsection{Characterization of ramified lifts of the standard representation
$\rho_n$ of $\BS(1,n)$}\label{s=bslifts}

The next proposition gives the key tool for lifting representations
under ramified covering maps.

\begin{prop}\label{p=lift}
Let $G$ be a group, and let
$\rho\co  G \rightarrow {\rm Diff }_+^\omega(S^1)$ be a representation
with global fixed point $p$.  
Let $\pi\co  S^1 \rightarrow S^1$ be a ramified covering map over
$p$ with signature ${\bf s}\in {\mathcal S}_d$, for some $d\geq1$.

Then for every homomorphism 
$h\co G \rightarrow {\rm Stab}_{D_d}({\bf s})$, there is a unique representation
$$\hat{\rho} = \hat{\rho}(\pi, h)\co  G \rightarrow {\rm Diff }^\omega(S^1)$$
such that, for all $\gamma\in G$,
\begin{enumerate}
\item $\hat{\rho}$ is a $\pi$--ramified lift of $\rho$; 
\item  $\hat{\rho}(\gamma)(q) = h(\gamma)(q)$, for each 
$q\in \pi^{-1}(p)$;
\item  $\hat{\rho}(\gamma)(e) =  h(\gamma)(e)$, for each oriented component
$e\in{\mathcal E}(\pi)$; 
\end{enumerate}

\noindent Furthermore, if $h$ takes values
in ${\rm Stab}_{C_d}({\bf s})$, then $\hat\rho$ takes values in 
$\Diff_+^\omega(S^1)$.
\end{prop}

Proposition~\ref{p=lift} is a special case of Proposition~\ref{p=cr}, 
which is proved in Section~\ref{s=morelifts}. Note that the representation $\rho$ in Proposition~\ref{p=lift} must be
orientation preserving, although the lift $\hat\rho$ might not be,
depending on where the image of $h$ lies.
There is also a criterion for lifting representations into
$\Diff^\omega(S^1)$ that are not necessarily orientation-preserving.
We discuss this issue in the next subsection.

\begin{lemma}\label{l=conj}  
Suppose that $\pi_1$ and $\pi_2$ are two ramified covering maps over $p$
such that ${\bf s}(\pi_2)$ lies in the $D_d$--orbit of ${\bf s}(\pi_1)$;
that is, suppose
there exists $\zeta\in D_d$ such that
${\bf s}(\pi_2) = \zeta ({\bf s}(\pi_1)).$  Then given any representation
$\rho\co  G \rightarrow {\rm Diff }_+^\omega(S^1)$ with global fixed point $p$
and homomorphism $h\co  G \rightarrow {\rm Stab}_{D_d}({\bf s}_1)$,
the representations $\hat{\rho}(\pi_1, h)$ and
$\hat{\rho}(\pi_2, \zeta h \zeta^{-1})$ are
conjugate in $\Diff^\omega(S^1)$, where
$$(\zeta h \zeta^{-1})(\gamma) := \zeta h(\gamma) \zeta^{-1}.$$

Furthermore, if $\zeta\in C_d$ and $h$ takes values
in ${\rm Stab}_{C_d}({\bf s})$, then $\hat{\rho}(\pi_1, h)$ and
$\hat{\rho}(\pi_2, \zeta h \zeta ^{-1})$ are
conjugate in $\Diff^\omega_+(S^1)$.
\end{lemma}

Lemma~\ref{l=conj} follows from Lemma~\ref{l=crconj}, which is proved in
Section~\ref{s=morelifts}. We now characterize the countably many conjugacy classes in ${\mathcal R}^\omega(\BS(1,n))$.
Note that the elements of ${\mathcal S}_d$ are totally 
ordered by the lexicographical order on $\R^n$.  
Hence we can write ${\mathcal S}_d$ as a disjoint union of 
$C_d$--orbits:
$${\mathcal S}_d = \bigsqcup_{\alpha\in A_+} C_d({\bf s}_\alpha),$$
where for each $\alpha\in A_+$, ${\bf s}_\alpha$ is the smallest element in
its $C_d$--orbit.  Similarly, there is an index set $A\supset A_+$
such that:
$${\mathcal S}_d = \bigsqcup_{\alpha\in A} D_d({\bf s}_\alpha).$$
Let $\overline {\mathcal S}_d = \{{\bf s}_\alpha\,\mid\,\alpha\in A\}$,
and let $\overline {\mathcal S}_d^+ = \{{\bf s}_\alpha\,\mid\,\alpha\in A_+\}$.
Finally, let $\overline{\mathcal S} = \bigcup_d \overline{\mathcal S}_d$
and let  $\overline{\mathcal S}^+ = \bigcup_d \overline{\mathcal S}_d^+$.

\dfinition 
Let $\rho_n\co  \BS(1,n) \rightarrow\Diff^\omega_+(S^1)$ denote the
standard projective action,
with global fixed point at $\infty\in\RP$.  Then we define:
$$
{\mathcal V} = \{ \hat{\rho}_n(\pi_{\bf s}, h)
\,\mid\, {\bf s}\in \overline{\mathcal
S}_d,\qua h\in {\rm Hom }(\BS(1,n), \, {\rm Stab}_{D_d}({\bf s}))/\equiv,
\qua d\in\N,\, d\geq 1\},$$
and let
$$
{\mathcal V}_+ = \{ \hat{\rho}_n(\pi_{\bf s}, h)
\,\mid\, {\bf s}\in \overline{\mathcal S}_d^+,\qua h\in 
{\rm Hom }(\BS(1,n), \, 
{\rm Stab}_{C_d}({\bf s})),
\qua d\in\N,\, d\geq 1\},$$
where, for ${\bf s}\in {\mathcal S}_d$,
$\pi_{{\bf s}}\co S^1\to S^1$ is the rational ramified cover over
$\infty$ with signature
${\bf s}$ given by Proposition~\ref{p=ratramcovexists}, and
$\equiv$ denotes conjugacy in ${\rm Stab}_{D_d}({\bf s})$.

\begin{prop}\label{p=c1conj} Each element of ${\mathcal V}$ and ${\mathcal V}_+$ represents a 
distinct conjugacy class of faithful representations. 

That is, if  $\hat{\rho}_n(\pi_{{\bf s}_1}, h_1), 
\hat{\rho}_n(\pi_{{\bf s}_2}, h_2)\in {\mathcal V}$ 
(resp.\ $\in {\mathcal V}_+$) are conjugate
in $\Diff^1(S^1)$ (resp.\ in $\Diff^1_+(S^1)$), then
${\bf s}_1 = {\bf s}_2$ and $h_1 = h_2$.  
\end{prop}

Proposition~\ref{p=c1conj} is proved at the end of Section~\ref{s=morelifts}.
Our main result, Theorem~\ref{t=main}, states that
the elements of ${\mathcal V}$ and ${\mathcal V}_+$ are the only faithful
representations of $\BS(1,n)$, up to conjugacy in $\Diff^\omega(S^1)$
and  $\Diff^\omega_+(S^1)$, respectively.

\subsection{Proof of Theorem~\ref{t=ramchar}}\label{s=ramaff}
To characterize the ramified lifts of $\Aff(\R)$, we 
need to define ramified lifts of orientation-reversing
diffeomorphisms.  In the end, our description is complicated by the
following fact: in contrast to lifts by regular covering maps,
ramified lifts of orientation-preserving diffeomorphisms can 
be orientation-reversing, and vice versa.

To deal with this issue, we introduce another action of the dihedral group 
$D_d=\langle a,b,\,\mid\, a^2=1, b^d=1, aba^{-1} = b^{-1}\rangle $ on ${\mathcal S}_d$
that ignores the edge labels completely.   To distinguish from
the action of $D_d$ on ${\mathcal S}_d$ already defined, we will write
$\zeta^\#\co {\mathcal S}_d\to {\mathcal S}_d$ for the action
of an element $\zeta\in D_d$. In this notation, the action is generated
by:
$$b^\#(s_1,  \ldots , s_d, o_1, \ldots, o_d) = (s_2, s_3, \ldots , s_d,
s_1, o_1, o_2, \ldots, o_d),$$
and
$$a^\#(s_1, \ldots , s_d, o_1, \ldots, o_d) = (s_1, s_d, s_{d-1}
\ldots , s_3
,s_2, o_1, o_2, \ldots, o_d).$$
For ${\bf s}\in{\mathcal S}_d$, we denote by 
${\rm Stab}^{\#}_{D_d}({\bf s})$ and ${\rm Stab}^{\#}_{C_d}({\bf s})$
the stabilizers of ${\bf s}$ in $D_d$ and $C_d$, respectively under this
action.  

\begin{lemma} For each ${\bf s}\in{\mathcal S}_d$, 
there exists a 
homomorphism $$\Delta_{\bf s}\co  {\rm Stab}^{\#}_{D_d}({\bf
s}) \to \Z/ 2\Z,$$  
such that ${\rm Stab}_{D_d}({\bf s}) = \ker(\Delta_{\bf s})$.
\end{lemma}

\proof Let ${\bf s}\in {\mathcal S}_d$ be given.
Clearly ${\rm Stab}_{D_d}({\bf s})$ is a subgroup of
${\rm Stab}^{\#}_{D_d}({\bf s})$. 
Let  $I\co {\mathcal S}_d\to
{\mathcal S}_d$ be the involution:
$$I(s_1,\ldots, s_d, o_1, \ldots , o_d) = (s_1, \ldots, s_d, -o_1,
\ldots , -o_d).$$
We show that for every $\zeta\in
{\rm Stab}^{\#}_{D_d}({\bf s})$,
either $\zeta({\bf s}) = {\bf s}$ (so that 
$\zeta\in{\rm Stab}_{D_d}({\bf s})$) or 
$\zeta({\bf s}) = I({\bf s})$.  This follows from the
property (2) in the definition of signature vector, which implies that
every element of ${\mathcal S}_d$ is determined by its first $d+1$ entries. 
Hence
we may define
$\Delta_{\bf s}(\zeta)$ to be $0$ if  $\zeta({\bf s}) = {\bf s}$
and $1$ otherwise.  Since $I$ is an involution, $\Delta_{\bf s}$ is
a homomorphism.  \eproof

\example  Consider the signature $${\bf
s}=(2,1,2,1,2,1,2,1,1,1,-1,-1,1,1,-1,-1).$$  
For this example we have ${\rm Stab}_{C_8}({\bf s}) = \langle b^4\rangle $,
${\rm Stab}_{D_8}({\bf s}) = \langle a, b^4\rangle $, ${\rm Stab}^\#_{C_8}({\bf s}) =
\langle b^2\rangle $,
and ${\rm Stab}^\#_{D_8}({\bf s}) = \langle a,b^2\rangle $.  In this example the
homomorphism $\Delta_{\bf s}$ is surjective, with nontrivial kernel.  
For 
${\bf s}=(2,1,4,1,2,1,4,1,1,1,-1,-1,1,1,-1,-1)$, on the other hand,
the image of  $\Delta_{\bf s}$ is trivial, and  
${\rm Stab}_{D_8}({\bf s}) = {\rm Stab}^\#_{D_8}({\bf s}) =
\langle a, b^4\rangle $, ${\rm Stab}_{C_8}({\bf s}) = {\rm Stab}^\#_{C_8}({\bf s}) =
\langle b^4\rangle $.

For a third example, recall that the stabilizer ${\rm Stab}_{D_6}({\bf s})$
of the signature vector ${\bf s} = (2,3,1,2,3,1,-1,-1, -1, 1, 1, 1)$ is
trivial.  Because of the rotational symmetry of the vertex labels,
however,  
${\rm Stab}^\#_{C_6}({\bf s}) =  {\rm Stab}_{D_6}^\#({\bf s}) = \langle b^3\rangle 
\simeq \Z/2\Z$.
In this example, $\Delta_{\bf s}$ is an isomorphism.

Let $G < \Diff^\omega(S^1)$ be a subgroup with global fixed point $p\in
S^1$:
$f(p) = p$, for all $f\in G.$
We now show how to assign, to each ${\bf s}\in{\mathcal S}$,
a subgroup $\hat G^{\bf s}$ consisting of ramified lifts of elements
of $G$. We first write 
%$\Diff^\omega(S^1) = \Diff^\omega_+(S^1)\sqcup \Diff^\omega_-(\R)$ and
$G = G_+\sqcup G_-$, where $G_+ = G\cap \Diff^{\omega}_+(S^1)$ is the 
kernel of the homomorphism $O\co G\to {\Z/2\Z}$ given by:
$$O(f) =  \begin{cases}0 & \hbox{if }f \hbox{ is
orientation-preserving,}\cr
1 & \hbox{ otherwise.}\end{cases}$$

Suppose that $\pi\co  S^1 \rightarrow S^1$ is a ramified
covering map over $p$.  Then, for every
$f\in G_+$,  Proposition~\ref{p=lift}
implies that for every $\zeta\in {\rm Stab}_{D_d}({\bf s}(\pi))$,
there is a unique lift $\hat{f}(\pi,\zeta)\in \Diff^\omega(S^1)$
satisfying: 

\begin{enumerate}
\item $\hat{f}(\pi, \zeta)$ is a $\pi$--ramified lift of $f$,
\item$\hat{f}(\pi,\zeta)(q) = \zeta(q)$, for all $q\in\pi^{-1}(p)$,
and
\item $\hat{f}(\pi,\zeta)(e) = \zeta(e)$, for all $e\in {\mathcal E}(p)$
\end{enumerate}
(Further, this lift is orientation-preserving 
if $\zeta\in {\rm Stab}_{C_d}(\s)$.)
Suppose, on the other hand, that $f\in G_-$.  In Section~\ref{s=morelifts},
we prove Lemma~\ref{l=liftdiff}, which 
implies that if $\zeta\in {\rm Stab}^\#_{D_d}(\s)$ satisfies:
\begin{eqnarray}\label{e=coset}
\zeta({\bf s}) = I({\bf s}),
\end{eqnarray}
then there exists a unique lift $\hat{f}(\pi,\zeta)\in
\Diff^\omega(S^1)$ satisfying {\rm(1)}--{\rm(3)} (and, further,
$\hat{f}(\pi,\zeta)\in
\Diff^\omega_+(S^1)$ if $\zeta\in {\rm Stab}^\#_{C_d}(\s)$.)
We can rephrase condition (\ref{e=coset}) as:
\begin{eqnarray*}
\Delta_{\bf s}(\zeta) = 1.
\end{eqnarray*}
To summarize this discussion, we have proved the following:
\begin{lemma}\label{l=alllifts}
If $f\in G$, and $\zeta\in{\rm Stab}_{D_d}^\#(\s)$,
then there exists a lift $\hat f(\pi,\zeta)$ satisfying {\rm(1)}--{\rm(3)}
if and only if: $$O(f) = \Delta_{\bf s}(\zeta).$$
\end{lemma}

For ${\bf s}\in {\mathcal S}$, let $\pi_{\bf s}$ be the
ramified covering map over $p$ with signature ${\bf s}$ given
by Proposition~\ref{p=ramcovexists}.  If $G < \Diff^\omega(S^1)$
has global fixed point $p$,
we define $\hat G^{\bf s} \subset \Diff^\omega(S^1)$ 
to be the fibered product of 
$G$ and ${\rm Stab}_{D_d}^\#(\s)$ with respect to
$O$ and $\Delta_{\bf s}$:
$$\hat G^{\bf s}:= \{\hat{f}(\pi_{\bf s},\zeta)\,\vert\, (f,\zeta)\in 
G\times {\rm Stab}_{D_d}^\#(\s),\, 
O(f) = \Delta_{\bf s}(\zeta)\}.$$
Similarly, we define:
$$\hat G^{\bf s}_+:= \{ \hat{f}(\pi_{\bf s},\zeta)\,\vert\,
(f,\zeta)\in G\times 
{\rm Stab}_{C_d}^\#(\s),\, 
O(f) = \Delta_{\bf s}(\zeta)\}.$$
Lemma~\ref{l=alllifts} tells us that $\hat G^{\bf s}$ coincides
with $\hat G^{\pi_{\bf s}}$, the set of all $\pi_{\bf s}$--ramified 
lifts of $G$, and, similarly, that $\hat G^{\bf s}_+ = \hat G^{\pi_{\bf s}}
\cap \Diff_+(S^1)$.
It follows from Lemma~\ref{l=liftdiff2} that $\hat G^{\bf s}$ and
$\hat G^{\bf s}_+$ are subgroups of $\Diff^\omega(S^1)$ and
$\Diff^\omega_+(S^1)$, respectively, with:
$$ \hat{f_1}(\pi_{\bf s},\zeta_1)\circ \hat{f_2}(\pi_{\bf s},\zeta_2) 
=  \widehat{f_1\circ f_2}(\pi_{\bf s},\zeta_1\zeta_2).$$
  Further, we have:

\begin{prop}\label{p=extension} Assume that
$G_-$ is nonempty.  Then $\hat G^{\bf s}$ and $\hat G^{\bf s}_+$ are 
both finite extensions of
$G_+$; there exist exact sequences:
\begin{eqnarray}\label{e=firstexact0}
1\to G_+ \to \hat G^{\bf s} \to {\rm Stab}_{D_d}^\#(\s) \to 1,
\end{eqnarray}
and
\begin{eqnarray}\label{e=secondexact0}
1\to G_+ \to \hat G^{\bf s}_+ \to {\rm Stab}_{C_d}^\#(\s) \to 1.
\end{eqnarray}
Furthermore, if the sequence 
$$1\to G_+ \to G \stackrel{O}{\to} O(G)\to 1$$
splits, where $O\co G\to \Z/2Z$ is  the orientation homomorphism, then
the sequences (\ref{e=firstexact0}) and (\ref{e=secondexact0}) are split, and
so is the sequence
\begin{eqnarray}\label{e=thirdexact}
1\to \hat G_+^{\bf s} \to \hat G^{\bf s} \to O(\hat G^{\bf s})\to 1.
\end{eqnarray}
\end{prop}

\proof The maps in the first sequence (\ref{e=firstexact0}) are given by:
$$\iota\co  G_+\to \hat G^{\bf s}\qquad f\mapsto \hat f(\pi_{\bf s}, id)$$
$$\sigma\co  \hat G^{\bf s} \to {\rm Stab}_{D_d}^\#(\s) \qquad \hat f(\pi_{\bf s},
\zeta)\mapsto \zeta.$$
It is easy to see that $\iota$ is injective and $\sigma$ is surjective.
Moreover, $\hat f(\pi_{\bf s},\zeta)$ is in the kernel of $\sigma$ if and only
if $\zeta = id$, if and only if $f$ is orientation-preserving, if and
only if $\hat f(\pi_{\bf s}, id)$ is in the image of $\iota$. Hence the first
sequence is exact. Similarly, the second 
sequence (\ref{e=secondexact0}) is exact.

Now suppose that $1\to G_+ \to G \to O(G)\to 1$ is split exact.  If 
$O(G)$ is trivial then $G_+ = G$, and there is nothing to prove.  If
$O(G) = \Z/2\Z$, then $G$ contains an involution $g\in G_-$ with
$g^2 = id$, namely the image of $1$ under the homomorphism $O(G)\to G$.
We use $g$ to  define a
homomorphism $j\co {\rm Stab}_{D_d}^\#(\s) \to \hat G^{\bf s}$ as follows:
$$j(\zeta) = \begin{cases}\hat{id}(\pi_{\bf s}, \zeta)&\hbox{ if
} \Delta_{\bf s}(\zeta) = 0 \cr  \hat{g}(\pi_{\bf s},\zeta)&\hbox{ if }
\Delta_{\bf s}(\zeta) = 1. \end{cases}$$
Hence the sequence (\ref{e=firstexact0}) is split.
The restriction of $j$ to ${\rm Stab}_{C_d}^\#(\s)$ splits the
sequence (\ref{e=secondexact0}).

If $O(\hat G^{\bf s})$ is trivial, then the last sequence
(\ref{e=thirdexact})
is trivially split.  If
$O(\hat G^{\bf s}) = \Z/2\Z$, then there exists 
a $\zeta\in {\rm Stab}_{D_d}^\#(\s)$ such that $\Delta_{\bf s}(\zeta) = 1$.
We then define $k\co  O(\hat G^{\bf s})\to G^{\bf s}$ by $k(0) = id$,
$k(1) = \hat{g}(\pi_{\bf s},\zeta)$, which implies that (\ref{e=thirdexact})
is split.
\eproof

Setting $G={\rm Aff}(\R)$, which has the global fixed point $\infty\in\RP$,
we thereby define
$\widehat{{\rm Aff}}^{\bf s}(\R)$ and $\widehat{{\rm Aff}}^{\bf
s}_+(\R)$, for ${\bf s}\in {\mathcal S}$. 
Let $\overline{\mathcal S}$ and $\overline{\mathcal S}_+$ be the
sets of signatures defined at the end of the previous subsection.  The
elements of  $\overline{\mathcal S}$ and $\overline{\mathcal S}_+$
are representatives of distinct orbits in ${\mathcal S}$ under the
dihedral and cyclic actions, respectively.  We now define:
$${\mathcal{RAFF}}:=
\{\widehat{{\rm Aff}}^{{\bf s}}(\R)\, \mid\,{\bf s}\in
\overline{\mathcal S}\},\, \hbox{  and  }\,{\mathcal{RAFF}}_+:=
\{\widehat{{\rm Aff}}^{{\bf s}}_+(\R)\, \mid\,{\bf
s}\in\overline{\mathcal S}_+\}.$$  Since ${\rm Aff}(\R)$
contains the involution $x\mapsto -x$, the sequence
$$1\to {\rm Aff}_+(\R)\to  {\rm Aff}(\R) \stackrel{O}{\to} \Z/2\Z \to 1$$
splits. Proposition~\ref{p=extension} implies that 
$$\widehat{{\rm Aff}}^{\bf s}(\R) \simeq  {\rm Aff}_+(\R)\times
{\rm Stab}_{D_d}^\#(\s),\qquad \widehat{{\rm Aff}}^{\bf s}_+(\R) \simeq  
{\rm Aff}_+(\R)\times
{\rm Stab}_{C_d}^\#(\s),$$
and either $\widehat{{\rm Aff}}^{\bf s}(\R) = 
\widehat{{\rm Aff}}^{\bf s}_+(\R)$ or
 $\widehat{{\rm Aff}}^{\bf s}(\R) \simeq  
\widehat{{\rm Aff}}^{\bf s}_+(\R)\times \Z/2\Z$,
depending on whether $\Delta_{\bf s}$ is surjective.  This proves (2)
of Theorem~\ref{t=ramchar}.

 Corollary~\ref{c=sigdet} implies that
  $\widehat{{\rm Aff}}^{{\bf s}_1}(\R)$
and $\widehat{{\rm Aff}}^{{\bf s}_2}(\R)$ are conjugate subgroups
in $\Diff^1(S^1)$, only if ${\bf s}_1\in D_d {\bf s}_2$. 
Similarly, $\widehat{{\rm Aff}}^{{\bf s}_1}_+(\R)$
and $\widehat{{\rm Aff}}^{{\bf s}_2}_+(\R)$ are conjugate subgroups
in $\Diff^1_+(S^1)$, if and only if ${\bf s}_1\in C_d {\bf s}_2$. 
This proves (1) of the theorem.

Finally, every finite dihedral or cyclic group $H$ is the stabilizer of
infinitely many ${\bf s}\in\overline{\mathcal S}$, and hence there are
infinitely many ${\bf s}\in\overline{\mathcal S}$ so that
$\widehat{{\rm Aff}}^{\bf s}(\R) \simeq {\rm Aff}_+(\R)\times H$.
This proves part (3).  
Property (4) follows from the definition of
${\mathcal{RAFF}}$.
This completes the proof of Theorem~\ref{t=ramchar}.

\section{The relation $fgf^{-1} = g^\lambda$: the central technical result}\label{s=schwarz}

In this section we analyze the relation $fgf^{-1} = g^\lambda$ 
near a common fixed point of $f$ and $g$.
If $f$ and $g$ are real-analytic,  then they can be locally conjugated into one 
of the ramified affine groups described
at the beginning of Section~\ref{s=ramlift} \cite{CM, Na, EISV}.
This gives a local characterization of
diffeomorphisms of $S^1$ satisfying this relation about each common
fixed point, and to obtain a global characterization, it is a matter of gluing
together these local ones. This was the way Ghys \cite{Gh1} proved that every
solvable subgroup of $\Diff^\omega(S^1)$ is metabelian.
To prove Theorems~\ref{t=main} and \ref{t=maincr}, we adapt the
arguments in \cite{Na} to the $C^r$ setting, 
where additional hypotheses on $f$ and $g$ are required.

The initial draft of this paper contained a completely different proof
of the local characterization of \cite{CM, Na, EISV} 
that does not rely on vector fields and works for
$C^r$ diffeomorphisms as well, under the right assumptions.
In the $C^\omega$ and $C^\infty$ case, this original proof gives identical 
results as the vector fields proof, but in the 
general $C^r$ setting, the proof using vector fields gives sharper results. 
At the end of this section,
we outline the alternate proof method. The main idea behind this method 
is to study the implications of the relation
$fgf^{-1} = g^\lambda$ for the Schwarzian derivative of $g$ near a
common fixed point at $0$.

We now state the main technical result of this section.
Let $[q, q_1)$ be a half open interval, let 
$r\in [2,\infty]\cup \{\omega\}$,
and let $f,g\in\Diff^r([q, q_1))$
be diffeomorphisms.
Assume that $g$ has no fixed points in $(q, q_1)$.

\noindent {\bf Standing Assumptions}\qua We assume that
either (A), (B), (C) or (D) holds:

(A)\qua $r=\omega$,  and there exists an {\em integer} $\lambda > 1$ such that
$f g f^{-1} = g^\lambda.$

(B)\qua $r \in [2,\infty)$, and there is an {\em integer} $\lambda > 1$
such that 
$f g f^{-1} = g^\lambda,$ and 
$f'(q) \leq \left(\frac{1}{\lambda}\right)^{\frac 1 {r-1}}$.
% or $f^\prime(q) > 1$.

(C)\qua $r = \infty$, $f^\prime(q) < 1$, and for some {\em integer} $\lambda > 1$,
$f g f^{-1} = g^\lambda.$

(D)\qua $r \in \{\infty, \omega\}$, $g$ is not infinitely flat, and 
there is a $C^\infty$ flow $g^t\co [q, q_1)\to [q, q_1)$ such that: 
\begin{enumerate}
\item $g^1 = g$, and
\item $f g f^{-1} = g^\lambda$
for some positive {\em real number} $\lambda \neq 1$.
\end{enumerate}

Assumptions (A), (B) and (C) will arise in the proof of Theorems~\ref{t=main}
and \ref{t=maincr}, and assumption (D) will arise in the proof of
Theorem~\ref{t=solv}.  The main technical result that we will use
in these proofs is the following.

\begin{prop}\label{p=schwarz}
Assume that either {\rm(A)}, {\rm(B)}, {\rm(C)} or {\rm(D)} holds.  Then 
there is a $C^r$ diffeomorphism 
$\alpha\co (q,q_1)\to (-\infty,\infty)\subset \RP$ such that
for all $p \in (q,\alpha^{-1}(0))$:
\begin{enumerate}
\item 
$\alpha(p) = \epsilon h(p)^s,$ where  
$h\co [q,\alpha^{-1}(0))\to [-\infty,0)$ is a $C^r$
diffeomorphism, $s$ is an integer satisfying $1 \leq s< r$, 
and $\epsilon\in\{\pm 1\}$,
\item 
$\alpha g (p) = \alpha(p) + 1 \; \hbox{  and  } \;
\alpha f (p) = \lambda \alpha(p).$
\end{enumerate}
\end{prop}

We start with a lemma describing which values of $f'(q)$ and
$g'(q)$ can occur. 

\begin{lemma}\label{l=derivs}  Assume that one of assumptions 
{\rm(A)}--{\rm(D)} holds. Then $g^\prime(q) = 1$, and either
\begin{enumerate}
\item $g^{(i)}(q) = 0\hbox{ for }2 \leq i \leq r$
(in particular, neither assumption {\rm(A)} nor {\rm(D)} can hold in this case),
or

\item $f^\prime(q) = (\frac{1}{\lambda})^{\frac{1}{s}}$ 
for some integer $1 \leq s < r$, and
$$g^{(i)}(q) = 0 \hbox{ for }2 \leq i \leq s, \; \; \hbox{and} \; \;
g^{(s+1)}(q) \neq 0.$$ 
\end{enumerate}

\end{lemma}

\proof  Since $fg = g^\lambda f$,
$$f^\prime(g(p)) \: g^\prime(p) = (g^\lambda )^\prime (f(p)) \: f^\prime(p).$$
When $p = q$, we thus have $g^\prime(q) = (g^\lambda )^\prime(q).$
But $(g^\lambda )^\prime (q) = (g^\prime(q))^\lambda $, and so
$g^\prime(q) = 1$.

Suppose that $f^\prime(q) = \kappa \neq 1$.  Then there is an interval $[q, p)$ on
which $f$ is $C^r$ conjugate to the linear map
$x \mapsto \kappa x$ (\cite{St}, Theorem 2).  So in local coordinates, 
identifying $q$ with $0$, 
\begin{eqnarray*}
f(x) &=& \kappa x, \; \; \hbox{and}\\
g(x) &=& x + ax^{s+1} + o(x^{s+1}) \; \; \hbox{for some } s \geq 1.
\end{eqnarray*}
Then 
\begin{eqnarray*}
fgf^{-1}(x) &=& x + (\frac{a}{\kappa^{s}})x^{s+1} + o(x^{s+1}), \; \; \hbox{and}\\
g^\lambda(x) &=& x + \lambda ax^{s+1} + o(x^{s+1}).
\end{eqnarray*}
So either 
\begin{enumerate}
\item $a = 0$, and therefore $g^{(i)}(q) = 0$ for $2 \leq i \leq r$, or
\item $a \neq 0$, in which case $\kappa = (\frac{1}{\lambda})^\frac{1}{s}$, and 
\begin{eqnarray*}
g^{(i)}(q) &=& 0 \hbox{ for } 2 \leq i \leq s \\
	&\neq& 0 \hbox{ for } i = s+1.
\end{eqnarray*}
\end{enumerate}

Now suppose that $f^\prime(q) = 1$. Then in a neighborhood of $q$, we can write
$$f(x) = x + bx^{k+1} + o(x^{k+1}) \; \;  \hbox{ and } \; \;  g(x) = x + ax^{s+1} + o(x^{s+1})$$
for some $k$, $s \geq 1$.
If $k = s$, then 
\begin{eqnarray*}
[f,g](x) &:=& fgf^{-1}g^{-1}(x) \; = \; x + o(x^{s+1}) \\
	&=& g^{\lambda-1}(x) \; = \; x + (\lambda-1)a x^{s+1} + o(x^{s+1}).
\end{eqnarray*}
So $a= 0$, and hence $g^{(i)}(0) = 0$ for $2 \leq i \leq r$.

If $k \neq s$, then we use the following well known result (see, eg, \cite{RS}):

\begin{lemma}\label{l=wkr}
If $f(x) = x + bx^{k+1} + o(x^{k+1})$ and $g(x) = x + ax^{s+1} + o(x^{s+1})$,
and if $s > k \geq 1$, then 
$$[f,g](x) = x + (s-k)ab x^{s+k} + o(x^{s+k}).$$
\end{lemma}

Assume that $s>k$.  (If $s<k$, then the proof is similar).
It follows from Lemma~\ref{l=wkr} that
\begin{eqnarray*}
x + (s-k)ab x^{s+k} + o(x^{s+k}) &=& [f,g](x) \\
	&=& g^{\lambda -1}(x) \; = \; x + (\lambda-1)ax^{s+1} + o(x^{s+1}).
\end{eqnarray*}
So either
\begin{enumerate}
\item  $k \geq 2$, and therefore $a = 0$ and $g^{(i)}(q) = 0$, for $2 \leq i \leq r$, or
\item  $k=1$, and therefore $b = \frac{\lambda -1}{s-1}$.
\end{enumerate}

But if $k = 1$, then
\begin{eqnarray*}
x + (s-1)2ab x^{s+1} + o(x^{s+1}) &=&  [f^2, g](x)\\
	&=& f^2 g f^{-2} g^{-1}(x) \; = \; g^{\lambda^2 -1}(x) \\
	&=& x + (\lambda^2 -1) ax^{s+1} + o(x^{s+1}).
\end{eqnarray*}

So $b = \frac{\lambda^2 -1 }{2(s-1)} = \frac{\lambda -1}{s-1}$, which is impossible,
since $\lambda \neq 1$.  Therefore $g^{(i)}(q) = 0$, for $2 \leq i \leq r$.
\eproof

\begin{lemma}\label{l=schwarz} 
Assume that {\rm(A)}, {\rm(B)}, {\rm(C)} or {\rm(D)} holds.
Then there is a neighborhood $[q,p_1) \subset [q,q_1)$ and
a $C^r$ map $\Alpha\co  [q, p_1) \rightarrow \R $
such that for all $p \in [q, p_1)$,
\begin{enumerate}
\item  $\Alpha(p) = o H(p)^s$, where $H\co  [q, p_1) \rightarrow [0, \infty)$ is a $C^r$ 
diffeomorphism, $1 \leq s < r$ is an integer, and $o\in \{\pm 1\}$;
\item  $$\Alpha f(p) = \frac{1}\lambda \Alpha(p);$$
\item  $$\Alpha g(p) = \frac{\Alpha(p)}{1-\Alpha(p)}.$$
\end{enumerate}
Consequently, $f^\prime(q) = (\frac{1}\lambda)^\frac{1}{s}$ for some integer 
$1 \leq s <r$, 
and $g^{(s+1)}(q) \neq 0$.
\end{lemma}

Before giving a proof of Lemma~\ref{l=schwarz}, we will show how this lemma implies 
Proposition~\ref{p=schwarz}.

\begin{lemma} Assume that {\rm(A)}, {\rm(B)}, {\rm(C)} or {\rm(D)} holds.  Let $A, H, s$, and $o$ be
given by Lemma~\ref{l=schwarz}.  For $p\in[q,p_1)$, let 
$$h(p) = \frac{-1}{H(p)};\qquad\qquad  \alpha_0(p) =
\frac{-1}{\Alpha(p)}.$$
Then $\alpha_0$ extends to a $C^r$ map $\alpha\co (q,q_1)\to
(-\infty,\infty)$ satisfying
the conclusions of Proposition~\ref{p=schwarz}.
\end{lemma}

\proof Lemma~\ref{l=schwarz} implies that for
all $p\in (q,p_1)$,  $\alpha_0 g(p) = \alpha_0(p) + 1$ and $\alpha_0 f (p) = \lambda \alpha_0(p)$. Since $\alpha_0$ has been defined in a fundamental domain for $g$,
we can now extend this map to a $C^r$ diffeomorphism $\alpha$ from $(q, q_1)$ to $(-\infty,\infty)$ as
follows.  Since $g$ has no fixed points in 
$(q, q_1)$, given
any $p \in (q, q_1)$, there is some $j \in \Z$ such that $g^{j}(p) \in (q,p_1)$.
Let $\alpha(p) = \alpha_0 (g^{j}(p)) - j$ (which is easily seen to be
independent of choice of $j$).  By construction, $\alpha g(p) = \alpha(p) + 1$ for all
$p \in (q, q_{1})$.  Since $f(p) = (g^{-j})^{\lambda} f g^{j}(p)$, we also have:
$$\alpha f(p) \; = \;
(\alpha (g^{-j})^{\lambda} \alpha_0^{-1})(\alpha_0 f \alpha_0^{-1}) (\alpha_0 g^{j}(p)) 
\; = \; \lambda (\alpha(p) + j) - j \lambda \; = \; \lambda \alpha(p).$$
Hence the conclusions of Proposition~\ref{p=schwarz} hold. \eproof

\medskip
\noindent
{\bf Proof of Lemma~\ref{l=schwarz}}\qua
We say that a $C^2$ function $c\co [a,b)\to [a,b)$ is a {\em $C^2$ contraction} if
$c'$ is positive on $[a,b)$ and $c(x) < x$, for all $x\in (a,b)$. Since $g$ has no fixed
points in $(q, q_1)$, either $g$ or $g^{-1}$ is a $C^2$ contraction.
We will assume until the end of the proof that $g$ is a $C^2$
contraction. Replacing $g$ by $g^{-1}$ does not change the relation
$fgf^{-1} = g^{\lambda}$.

Since $g$ has no fixed points in $(q,q_1)$, there is a  a unique $C^1$ vector field $X_0$ on
$[q,q_1)$ that generates a $C^1$ flow $g^t$ such that
$g\vert_{[q,q_1)} = g^1$ (Szekeres, see \cite{Nav} for a discussion).

\begin{lemma}
For all $j \in \N$ and $x\in[q,q_1)$, $f^{-j}gf^j(x) = g^{\frac{1}{\lambda^j}}(x)$.
\end{lemma}
 
\proof  
We will use the following result of Kopell:

\begin{lemma}[\cite{Ko} Lemma 1]\label{l=kopell}   Let $g\in\Diff^2[q, q_1)$ be
a $C^2$ contraction that embeds in a  $C^1$ flow $g^t$, so that
$g=g^1$.  
If $h\in \Diff^1[q, q_1)$
satisfies $hg=gh$, then $h = g^t$ for some $t \in \R$.  
\end{lemma}
It follows from the relation $f^jgf^{-j}=g^{\lambda^j}$ that
$f^{-j}gf^j$ commutes with $g$, and therefore
Lemma~\ref{l=kopell} implies that $f^{-j}gf^j = g^t$ for some $t \in \R$.
This relation also implies that
$(f^{-j}gf^j)^{\lambda^j} = g$.  So $f^{-j}gf^j = g^{\frac{1}{\lambda^j}}$.
\eproof

Let $\kappa = f^\prime(q)$.  We may assume, by Lemma~\ref{l=derivs}, that
$\kappa \neq 1$, and therefore
there is an interval $(q,p_1)$ on which $f$ has no fixed points, and
a $C^r$ diffeomorphism $H\co  [q, p_1)
\rightarrow [0, \infty)$ such that
$H\: f \: H^{-1}(x) =  \kappa x$ (\cite{St}, Theorem 2).
The diffeomorphism $H$ is unique up to multiplication by a constant.
Let $F = HfH^{-1}$ and let $G = HgH^{-1}$.
Since we have assumed that $g$ is a contraction, we have $g([q, p_1)) \subseteq [q, p_1)$.

Let $X$ be the push-forward of the vector field $X_0$ to $[0,\infty)$
under $H$, and let $G^t$ be the semiflow generated by $X_0$, so that
$G=G^1$.

\begin{lemma}\label{l=flat}  If $F^\prime(0) \leq (\frac{1}{\lambda})^\frac{1}{r-1}$ and $G$ is $r$--flat at $0$,
then $X(x) = 0$ on $[0, \infty)$.
\end{lemma}

\proof  We will show that for all $x \in [0, \infty)$,
$$\lim_{t \rightarrow 0} \frac{G^t(x) -x}{t} = 0.$$
Since the limit exists, it is enough to show that it converges to $0$ for a
subsequence $t_i \rightarrow 0$.  We will use the subsequence $t_i = \frac{1}{\lambda^i}$.
Writing $\kappa = F^\prime(0)$ as before, we have
 $$F(x) = \kappa x \; \; \hbox{and } \; G(x) = x + R(x),$$
where $R(x)/x^r\to 0$ as $x\to 0$, and therefore:
$$G^{t_i}(x) = F^{-i}GF^i(x) = x + \frac{1}{\kappa^i} R((\kappa^i x)).$$
So
\begin{eqnarray*}
0 \; \leq \; \lim_{i \rightarrow \infty} \frac{|G^{t_i}(x) - x|}{t_i} &=& 
	\lim_{i \rightarrow \infty}  \lambda^{i}
\left(\frac{1}{\kappa^{i}} |R( \kappa^{i}x)|\right) \\
&=& \lim_{i \rightarrow \infty}  (\lambda \kappa^{r-1})^i x^r \;
	\frac{|R(( \kappa^{i}x))|}{( \kappa^{i}x)^r} \\
&\leq&  \lim_{i \rightarrow \infty}  x^r \;
	\frac{|R(( \kappa^{i}x))|}{( \kappa^{i}x)^r}  \; = \;0,
\end{eqnarray*} since $\kappa^{r-1} \leq \frac{1}{\lambda}$.
\eproof

\begin{cor}\label{c=notflat} Under any assumption {\rm(A)}--{\rm(D)}, $g$ is not $r$--flat at $q$, and therefore $f^\prime(q)
= (\frac{1}{\lambda})^\frac{1}{s}$, for some integer $1 \leq s <r$.
\end{cor}

\proof  Clearly $g$ cannot be infinitely flat  if (A) or (D) holds. Under assumption
(C),  $f'(q) < (\frac{1}{\lambda})^\frac{1}{k}$, for some $k>0$ and $f,g$ are
$C^k$, so (C) reduces to (B).
By Lemma~\ref{l=flat}, under assumption (B), if $g$ is $r$--flat at $q$, then
the semiflow $G^t$ is tangent to the trivial vector field, 
$X(x) = 0$.  But then $G = id$, and therefore $g = id$ on $[q, p_1)$,
contradicting the assumption that $g$ has no fixed points in $(q, q_1)$.
\eproof

\begin{lemma}\label{l=notflat} If $F^\prime(0) 
= (\frac{1}{\lambda})^\frac{1}{s}$
for some integer $1 \leq s <r$, then for some $a < 0$,
$X(x) = ax^{s+1}$ on $[0,\infty)$.
\end{lemma}

\proof  As in the proof of Lemma~\ref{l=flat}, it is enough to show that
for all $x \in [0, \infty)$, and for $t_i = \frac{1}{\lambda^i}$,
$$\lim_{i \rightarrow \infty} \frac{G^{t_i}(x) -x}{t_i} = ax^{s+1}$$
for some $a \in \R$.  If $F^\prime(0) 
= (\frac{1}{\lambda})^\frac{1}{s}$
for some integer $1 \leq s <r$, then by Lemma~\ref{l=derivs}, 
$$F(x) = (\frac{1}{\lambda})^\frac{1}{s} x \; \; \hbox{ and } \; \;
G(x) = x + ax^{s+1} + R(x)$$
for some $a \in \R$, 
where $R(x)/x^{s+1} \rightarrow 0$ as $x \rightarrow 0$.  The value of $a$ depends
on the choice of linearizing map $h$ for $f|_{[q,p_1)}$.  
For all $i \in \N$,
\begin{eqnarray*}
G^{t_i}(x) \; = \;  F^{-i}GF^i(x) &=& \lambda^\frac{i}{s} \;
G\left(\frac{x}{\lambda^\frac{i}{s}}\right) \\
	&=& x + a  x^{s+1} \frac{1}{\lambda^i} + \lambda^\frac{i}{s} \;
R\left(\frac{x}{\lambda^\frac{i}{s}}\right).
\end{eqnarray*}

So
\begin{eqnarray*}
\lim_{i \rightarrow \infty} \frac{G^{t_i}(x) - x}{t_i} &=& 
\lim_{i \rightarrow \infty} \lambda^i \left(a  x^{s+1} \frac{1}{\lambda^i} 
+ \lambda^\frac{i}{s} \;
R\left(\frac{x}{\lambda^{\frac{i}{s}}}\right) 
\right) \\
&=& \lim_{i \rightarrow \infty} ax^{s+1} + 
R\left(\frac{ x }{ \lambda^{ \frac{i}{s} } } \right) \;
\frac{\lambda^{\frac{i(s+1)}{s}}}{{x^{s+1}}} \; x^{s+1} \; \; = \; \; ax^{s+1}
\end{eqnarray*}
Since $G$ is
a contraction on $[0,\infty)$, it follows that $a \leq 0$,
and since $G \neq id$, we must have $a < 0$.
\eproof

\begin{cor}\label{c=H}  If one of {\rm(A)}--{\rm(D)} holds, then $f^\prime(q) = (\frac{1}{\lambda})^\frac{1}{s}$ for some
integer $1 \leq s <r$, and, after a suitable
rescaling of the linearizing map $H$, $$G(x) = \frac{x}{(1+x^s)^\frac{1}{s}}.$$ 
\end{cor}

\proof Choose $H$ so that $a= -1/s$.  Solving the differential equation $\frac{\delta}{\delta t}
G^t(x)$ $= aG^t(x)^{s+1}$ with initial condition $G^0(x) = x$, we obtain
$G^t(x) = x/(1 + t x^{s})^\frac{1}{s}$.
Since  $G(x) = G^1(x)$, the conclusion follows.
\eproof

To complete the proof of Lemma~\ref{l=schwarz} assuming that $g$ is a contraction, 
let $s, H$ be given by
Corollary~\ref{c=H}, and let $o = - 1$.
 Then Corollary~\ref{c=H} implies that $A(p) =
o (H(p)^s)$ satisfies
the desired conditions. If $g$ is not a contraction, we replace $g$ by
$g^{-1}$ in the proof.  Setting $o = 1$, we obtain the desired
conclusions. \eproof

\subsection{Idea behind an alternate proof of Proposition~\ref{p=schwarz}}
Suppose that $f$ and $g\neq id$ are $C^r$ diffeomorphisms, defined in a 
neighborhood of $0$ in $\R$, both fixing the origin,
and satisfying the relation:
 $$fgf^{-1} = g^\lambda,$$
for some $\lambda>1$.  In this context, the conclusion of Proposition~\ref{p=schwarz} 
can be reformulated as follows:  $f$ and $g$ are
conjugate, via $-1/h$, to the maps 
$$ x \mapsto (\frac{1}{\lambda})^\frac{1}{s} x \; \; \hbox{and} \; \;
x \mapsto \frac{x}{(1 - ox^s)^\frac{1}{s}}$$
for some integer $1 \leq s \leq r$ and some $o \in \{\pm 1\}$. 
The proof of Proposition~\ref{p=schwarz} uses vector fields; here we sketch
an alternative proof of this reformulation, 
using the Schwarzian derivative.  This 
sketch can be made into a
complete proof of Proposition~\ref{p=schwarz} under assumptions
(A), (C) and (D), but gives a weaker result
in case (B): for this proof we will need both $r \geq 2s + 1$ and $f^\prime(q) \leq
(\frac{1}{\lambda})^{1/(s-1)}$, for some $s \geq 1$.

For simplicity, assume that $r = \omega$ and that $\lambda = 2$.  
First note that, since $g$ is not infinitely flat, Lemma~\ref{l=derivs} 
implies that 
$f'(0)\in \{\left(\frac{1}{2}\right)^{\frac{1}{s}}\,\mid\, s\geq 1\}$.  After
conjugating $f$ and $g$ by an analytic diffeomorphism, we 
may assume, then, that:
$$f(x) = \frac{x}{2^{\frac{1}{s}}},$$
for some $s\geq 1$.

Let $F(x) = f(x^{\frac{1}{s}})^s = x/2$ and let 
$G(x) = g(x^{\frac{1}{s}})^s$.  Rewriting the relation
$FGF^{-1} = G^2$, we obtain:
$$\frac{1}{2} G(2x) = G^2(x);$$
rearranging and iterating this relation, we obtain:
\begin{eqnarray}\label{e=iterate}
G(x) = 2^k G^{2^k} \left(\frac{x}{2^k}\right),
\end{eqnarray}
for all $k\geq 1$.

Recall that the Schwarzian derivative of a $C^3$ function $H$ is
defined by:
$$S(H)(x) = \frac{H^{\prime \prime \prime}(x)}{H^\prime(x)} - \frac{3}{2}
\left( \frac{H^{\prime \prime}(x)}{H^\prime(x)} \right)^2,$$
and has the following properties:

\begin{enumerate}
\item $S(H)(x) = 0$ for all $x$ iff G is M{\"o}bius, and
\item for any $C^3$ function $K$, $S(H \circ K)(x) =  K^\prime(x)^2 S(H)(K(x)) + S(K)(x)$.
\end{enumerate}
Combining these properties with (\ref{e=iterate}), we will show that
$S(G)=0$, which implies that $G$ is M\"obius.  Lemma~\ref{l=derivs}
implies that $G^\prime(0)=1$ and
$G^{\prime \prime}(0) \neq 0$, 
so we have $G(x) = \frac{x}{1-ox}$ for $o \in \{\pm1\}$.   Writing
$g(x) = G(x^s)^{1/s}$, we obtain the desired result.

The first thing to check is that $G$ is $C^3$. To obtain this, we use 
a slightly stronger version of Lemma~\ref{l=derivs} (whose proof is left as an exercise), which states that,
if $g$ is not infinitely flat, then 
$$g(x) = x + ax^{s+1} + bx^{2s+1} + \cdots$$
Performing the substitution $G = g(x^{1/s})^s$ in this series,
one finds that $G$ is $C^3$. (This requires that $g$ be at least
$C^{2s+1}$, in contrast to the proof of Proposition~\ref{p=schwarz}, which requires 
only $C^{s+1}$).

Equation (\ref{e=iterate}) implies that 
\begin{eqnarray*}
S(G)(x) & = & \frac{1}{2^{2k}} S(G^{2^k})(\frac{x}{2^k}),
\end{eqnarray*}
for all $k\geq 1$.  Thus, by the cocycle condition (2) of the Schwarzian, 
we have:
\begin{eqnarray}\label{e=series1}
S(G)(x) & = & \frac{1}{2^{2k}} \sum_{i=1}^{2^k} S(G)(G^{i-1}(\frac{x}{2^k})) 
\,\left((G^{i-1})'(\frac{x}{2^k})\right)^2\\
&=& \frac{1}{2^{2k}} \sum_{i=1}^{2^k} S(G)( x_i ) 
\,\left(\Pi_{j=1}^{i-1} G'(x_j)\right)^2
\end{eqnarray}
where $x_i:= G^{i-1}(\frac{x}{2^k})$. 

Fix $x$, and assume without loss of generality that $G^j(x)\to 0$ as
$j\to\infty$.  Since 
$G$ is $C^3$ and $G'(0) = 1$, there is a constant $C>0$ such that
$|G'(x_i)| \leq 1 + \frac{C}{2^k}$, and $|S(G)(x_i)| \leq C$,
for all $i$ between $1$ and $2^k$ and all $k\geq 1$.  
Combined with (\ref{e=series1}), this
gives us a bound on the Schwarzian of $G$ at $x$:
\begin{eqnarray*}
|S(G)(x)| & \leq & \frac{C}{2^{2k}} \sum_{i=1}^{2^k}  
\,\left(1 + \frac{C}{2^k}  \right)^{2(i-1)}\\
&\leq &  \frac{C}{2^{2k}} \left(\frac{1 - (1+ \frac{C}{2^k})^{2^{k+1}}}{1 -  (1+ \frac{C}{2^k})^2}\right)\\
&\leq & \frac{1}{2^{k}}\left(\frac{e^{2C}-1}{2 + \frac{C}{2^k}}
\right),
\end{eqnarray*}
for all $k\geq 1$.  Hence 
$S(G)(x) = 0$, for all $x$, which implies that $G$ is M\"obius.

\section{Further properties of ramified covers: proofs of Proposition~\ref{p=lift},
Lemma~\ref{l=conj} and Proposition~\ref{p=c1conj}}\label{s=morelifts}

The next lemma describes a useful normal form for ramified covering maps.

\begin{lemma}\label{l=normal}
Let $\pi: \RP \rightarrow \RP$ be a ramified covering map over $0$, where
$\pi^{-1}(0) = \{x_1, \ldots, x_d\}$.  Let ${\bf s} = (s(x_1), \ldots, s(x_d),
o_1,\ldots,o_d)$
be the signature of $\pi$.  Then given any $x_i \in \pi^{-1}(0)$, there is
a neighborhood $U$ of $x_i$ and an analytic diffeomorphism $h\co  U \rightarrow
\R$ such that for all $x \in U$,
$$\pi(x) = h(x)^{s(x_i)}.$$
\end{lemma}
\proof  In local coordinates at $x_i$, identifying $x_i$ with $0$, we can write
\begin{eqnarray*}
\pi(x) &=& a x^s (1 + O(x))\\
&=& a x^sg(x) 
\end{eqnarray*}
where $a > 0$, $s = s(x_i)$ and $g(x) = (1 + O(x))$.  Let $h(x) =
a^\frac{1}{s} x g(x)^\frac{1}{s}$.  Then $h(x)$ is analytic in a neighborhood
of $0$, and $\pi(x) = h(x)^s$. \eproof

This lemma motivates the following definition.

\dfinition A {\em $C^r$ ramified cover over $p\in S^1$} is a map
$\pi\co  S^1 \rightarrow S^1$ satisfying:
\begin{enumerate}
\item  $\pi^{-1}(p) = \{q_1, q_2, \ldots, q_d\}$, where 
$q_1 <  q_2 < \ldots< q_d$;
\item  the restriction of $\pi$ to 
$\pi^{-1}(S^1\setminus\{p\})$ is a regular
$C^r$ covering map  onto $S^1\setminus\{p\}$ of degree $d \geq 1$;
\item  for all $1 \leq i \leq d$, there are neighborhoods $U_i$ of $q_i$
and $V$ of $p$, and $C^r$ charts $h_i\co  U_i \rightarrow \R$ and 
$k_i\co  V \rightarrow \R$ with
$h_i(q_i) = 0$ and $k_i(p) = 0$, such that
$$k_i \: \pi \: h_i^{-1}(x) = x^{s_i}$$ for some integer $s_i > 0$.
\end{enumerate}

\remark By Lemma~\ref{l=normal}, a ramified cover is a
$C^\omega$ ramified cover.

We define the signature of a $C^r$ ramified cover in the obvious way.

\dfinition  Let $\pi_1$ and $\pi_2$ be $C^r$ ramified covering maps of degree $d$  over $p_1$ and $p_2$, respectively.  Fix an orientation 
preserving identification
between $\pi_1^{-1}(p)$ and $\pi_2^{-1}(p)$ and between ${\mathcal E}(\pi_1)$
and ${\mathcal E}(\pi_2)$.
Suppose that $f\in \Diff^r(S^1)$   satisfies $f(p_1) = p_2$,
and let  $\zeta \in D_d$. We say that  $\hat f\in \Diff^r(S^1)$ is
a {\em $(\pi_1, \pi_2, \zeta)$--ramified lift of $f$} if:
\begin{enumerate}
\item  $\hat{f}(q) = \zeta(q)$, for all $q\in \pi_1^{-1}(p_1)$,
\item $\hat{f}(e) = \zeta(e)$, for all $e\in {\mathcal E}(\pi_1)$,
and
\item the following diagram commutes:

\fig{
\object {a}  0,8 {S^1};
\object {b}  8,8 {S^1};
\object {c}  0,4  {S^1};
\object {d}  8,4  {S^1};
\arrow  {a}  {b}  {\hat{f}};
\arrow  {c}  {d}  {f};
\arrow  {a}  {c}  {\pi_1};
\arrow  {b}  {d}  {\pi_2};}
\end{enumerate}

\begin{lemma}\label{l=liftdiff}
Let $\pi_1$, $\pi_2$ and $f$ be as above.  
Suppose that $\zeta\in D_d$ satisfies
\begin{itemize}
\item $\zeta(\s(\pi_1)) = \s(\pi_2)$, if $f\in \Diff^r_+(S^1)$, or
\item $\zeta(\s(\pi_1)) = I(\s(\pi_2))$, if $f\in \Diff^r_-(S^1)$,
\end{itemize}
where $I\co {\mathcal S}_d\to {\mathcal S}_d$ is the involution
that reverses the sign of the last $d$ coordinates. 

Then there exists a unique  $(\pi_1, \pi_2 , \zeta)$--ramified lift
of $f$.  We denote this lift by $\hat{f}(\pi_1, \pi_2, \zeta)$,
or by $\hat{f}(\pi,\zeta)$, if $\pi_1=\pi_2 = \pi$.

Furthermore, we have that
if $\zeta \in C_d$, then  $\hat{f}(\pi_1, \pi_2, \zeta)\in
\Diff^r_+(S^1)$.
\end{lemma}

\proof 
Suppose first that $f$ preserves orientation.
Since the restriction of $\pi_1$ to 
$\pi_1^{-1}(S^1 \setminus \{p_1\})$ and the restriction of $\pi_2$ to
$\pi_2^{-1}(S^1 \setminus \{p_2\})$ are both regular
$C^r$ covering maps of degree $d$, for any $\zeta\in D_d$ 
there is a unique $C^r$ diffeomorphism $\hat{f}_0\co  S^1\setminus\pi_1^{-1}\{p_1\}
\rightarrow S^1 \setminus \pi_2^{-1}\{p_2\}$
such that $\hat{f}_0(e) = \zeta(e)$
for all $e\in {\mathcal E}(\pi_1)$, and
the diagram (3) commutes on the restricted domains.
The condition ${\bf s}(\pi_2) = \zeta ({\bf s}(\pi_1))$
implies that $\hat{f}_0$ extends to a unique homeomorphism
$\hat{f}$ such that $\hat{f}(q) = \zeta(q)$, for all $q\in
\pi_1^{-1}(p_1)$
and such that the diagram in (3) commutes.
It remains to show that $\hat{f}$ is a $C^r$ diffeomorphism.

It suffices to show that $\hat{f}$ 
is a $C^r$ diffeomorphism at each $q \in \pi_1^{-1}(p_1)$.
By Lemma~\ref{l=local}, there are local coordinates near $q$ and
$\hat f(q)$, identifying both of these points with $0$, such that
$$(\hat{f}(x))^{s(\hat{f}(q))} \; = \; f(x^{s(q)})$$
for some integers $s(\hat f(q))$ and $s(q)$.
Since  ${\bf s}(\pi_2) = \zeta ({\bf s}(\pi_1))$, we have
$s(\zeta(q)) = s(q)$.  Let $j = s(q)= s(\hat{f}(q))$.  
Since $f$ is $C^r$ and has a fixed point at $0$,
$$f(x) = a_1 x + a_2 x^2  + \ldots + x^r + o(x^r)$$
and we can assume that the  coordinates have been 
chosen so that $a_1 > 0$.  So near $x = 0$,
\begin{eqnarray*}
\hat{f}(x) &=& (a_1 x^j + a_2 x^{2j} + \ldots + x^{rj} + o(x^{rj}))^\frac{1}{j} \\
	&=& x \: (a_1 + a_2 x^j + \ldots + x^{(r-1)j} + o(x^{(r-1)j}) )^\frac{1}{j},
\end{eqnarray*} where the root is chosen so that $\hat f^\prime(0) > 0$.
Since $a_1 > 0$ and $r\geq 2$, 
$\hat{f}$ is a $C^r$ diffeomorphism at $0$. Similarly, we see that
if $f$ is analytic, then $\hat{f}$ is analytic.  
Finally, we note that since $f$ is orientation preserving, if
$\zeta \in C_d$, then $\hat{f}$ must also be
orientation preserving.  

Now suppose that $f \in $ Diff$^r_-(S^1)$, and that $\zeta(\s(\pi_1))
= I(s(\pi_2))$.   Let $\overline \pi_1 = f \circ \pi_1$.  Setting $\hat{f}$
to be the $(\overline\pi_1,\pi_2, \zeta)$--lift of the identity map,
we obtain the desired conclusions.\eproof

\begin{lemma}\label{l=liftdiff2}  
Let $f_1$ and $f_2$ be $C^r$ diffeomorphisms of $S^1$, 
both with a fixed point at $p$, let $\pi\co  S^1 \rightarrow S^1$ be a
$C^r$ ramified covering map over $p$ with signature ${\bf s}$,  
and let $\zeta_1$, $\zeta_2 \in
D_d$. Suppose that $\zeta_1$ and $\zeta_2$ satisfy
$\zeta_i({\bf s}) = {\bf s}$ if $f_i\in \Diff^r_+(S^1)$,
and $\zeta_i({\bf s}) = I({\bf s})$ if $f_i\in \Diff^r_-(S^1)$.
Then  
$$\hat{f_2}(\pi, \zeta_2) \circ \hat{f_1}(\pi, \zeta_1) = \widehat{f_2 \circ f_1}
(\pi, \zeta_2 \circ \zeta_1).$$
\end{lemma}

\proof  The map $\hat{f_2}(\pi, \zeta_2) \circ \hat{f_1}(\pi,
\zeta_1)(q)$
satisfies:
\begin{enumerate}
\item  $\hat{f_2}(\pi, \zeta_2) \circ \hat{f_1}(\pi, \zeta_1)(q) = 
{\zeta_2 \circ \zeta_1}(q)$, for all $q\in\pi^{-1}(p)$,
\item  $\hat{f_2}(\pi, \zeta_2) \circ \hat{f_1}(\pi, \zeta_1)(e) = 
\zeta_2 \circ \zeta_1(e)$, for all $e\in {\mathcal E}(\pi)$,
and 
\item  the following diagram commutes:

\fig{
\object {a}  0,8 {S^1};
\object {b}  8,8 {S^1};
\object {c}  0,4  {S^1};
\object {d}  8,4  {S^1};
\object {e}  16,8  {S^1};
\object {f}  16,4  {S^1};
\arrow  {a}  {b}  {\hat{f_1}(\pi, \zeta_1)};
\arrow  {c}  {d}  {f_1};
\arrow  {a}  {c}  {\pi};
\arrow  {b}  {d}  {\pi};
\arrow  {b}  {e}  {\hat{f_2}(\pi, \zeta_2)};
\arrow  {d}  {f}  {f_2};
\arrow  {e}  {f}  {\pi};}
\end{enumerate}

By Lemma~\ref{l=liftdiff}, we must have 
$\; \hat{f_2}(\pi, \zeta_2) \circ \hat{f_1}(\pi, \zeta_1) = \widehat{f_2 \circ f_1}
(\pi, \zeta_2 \circ \zeta_1).$  \eproof

The following proposition is a $C^r$ version of Proposition~\ref{p=lift}.
\begin{prop}\label{p=cr}  
Suppose that $G$ is a group, and that
$\rho\co  G \rightarrow {\rm Diff }^r_+(S^1)$ is a representation
with global fixed point $p$.  Let $\pi\co S^1\to S^1$ be a $C^r$
ramified cover over $p$ with  signature vector ${\bf s}$. Then for 
every homomorphism $h\co G \rightarrow {\rm Stab}_{D_d}({\bf s})$, 
there is a unique
representation $$\hat{\rho} = \hat{\rho}(\pi, h): 
G \rightarrow {\rm Diff }^r(S^1)$$
such that, for all $\gamma \in G$,  $\hat \rho(\gamma)$ is the $(\pi, h(\gamma))$--
ramified lift of $\rho(\gamma)$.
If  $h$ takes values in Stab$_{C_d}({\bf s})$,
then $\hat{\rho}$ takes values in $\Diff^r_+(S^1)$.
\end{prop}

\proof  This follows immediately from the previous two lemmas. \eproof

The following lemma is a $C^r$ version of Lemma~\ref{l=conj}.
%\marginpar{We should restate this slightly to make it match up with the
%new statement of Lemma~\ref{l=conj}}

\begin{lemma}\label{l=crconj}
Let $G$ be a group, and let
$\rho\co  G \rightarrow {\rm Diff }^\omega(S^1)$ be a representation
with global fixed point $p$.  
 Let $\pi_1, \pi_2 \co  S^1 \rightarrow S^1$ be 
$C^r$ ramified covers over $p\in S^1$, with ${\bf s}(\pi_1) = \zeta({\bf s}(\pi_2))$,
for some $\zeta \in D_d$.

Then for
every homomorphism $h\co G \rightarrow {\rm Stab}_{D_d}({\bf s})$,
the representation $\tilde{\rho}(\pi_1, h)$ is conjugate to 
$\tilde{\rho}(\pi_2, \zeta h \zeta^{-1})$ 
in $\Diff^r(S^1)$, where $(\zeta h \zeta^{-1})(\gamma): = \zeta h(\gamma)\zeta^{-1}$.  
If $\rho$ takes values in $\Diff^\omega_+(S^1)$, if $\zeta \in C_d$, 
and if $h$ takes values in Stab$_{C_d}({\bf s})$, then $\tilde{\rho}(\pi_1, h)$ 
and $\tilde{\rho}(\pi_2, \zeta h \zeta^{-1})$ are conjugate in $\Diff^r_+(S^1)$.
\end{lemma}

\proof  This lemma follows from the diagram below, which commutes by Proposition~\ref{p=cr}
and Lemma~\ref{l=liftdiff}.
(Here $\widehat{id} =\widehat{id}(\pi_1,\pi_2,\zeta^{-1})$).

\fig{
\object {a}  0,10 {S^1};
\object {b}  12,10 {S^1};
\object {c}  0,2  {S^1};
\object {d}  12,2  {S^1};
\object {e}  10,20 {S^1};
\object {f}  22,20 {S^1};
\object {g}  10,12  {S^1};
\object {h}  22,12  {S^1};
\arrow  {a}  {b}  {\hat\rho(\pi_1,h)};
\arrow  {c}  {d}  {\rho};
\arrow  {a}  {c}  {\pi_1};
\arrow  {b}  {d}  {\pi_1};
\arrow  {e}  {f}  {\hat\rho(\pi_2,\zeta h \zeta^{-1})};
\arrow  {g}  {h}  {\rho};
\arrow  {e}  {g}  {\pi_2};
\arrow  {f}  {h}  {\pi_2};
\arrow  {a}  {e}  {\widehat{id}};
\arrow  {b}  {f}  {\widehat{id}};
\arrow  {c}  {g}  {id};
\arrow  {d}  {h}  {id};}
\eproof

Consider two lifts $\hat f(\pi_1,\zeta)$, $\hat f(\pi_2,\zeta)$
of the same diffeomorphism $f$ (or, more generally, of conjugate 
diffeomorphisms).
For purely topological reasons, if these lifts are conjugate by a
map with rotation number $0$, then
${\bf s}(\pi_1)$ and ${\bf s}(\pi_2)$ have the same length $2d$,
and the final $d$ entries in these vectors must agree. (More generally, if
the conjugacy has nonzero rotation number, then the final $d$ entries
of the first vector must lie in the $D_d$--orbit of the final $d$ entries 
of the second).  We now examine
the first $d$ entries of both vectors.  We show that, 
under appropriate regularity 
assumptions on $f$ and on the conjugacy, these entries must also agree, so
that ${\bf s}(\pi_1)= {\bf s}(\pi_2)$.  The next lemma is the key reason
for this.

\begin{lemma}\label{l=c2contr} Let $c\co [0,\infty)\to [0,\infty)$ be a $C^2$ contraction.  
Suppose that, for some integers $m,n>0$, the maps $v_1(x) = c(x^m)^{1/m}$
and $v_2(x) = c(x^n)^{1/n}$ are conjugate by a $C^1$ diffeomorphism 
$h\co [0,\infty)\to [0,\infty)$.  Then $m=n$.
\end{lemma}

\proof  Since $c$ is a $C^2$ contraction, the standard distortion estimate 
(see, eg \cite{Ko}) implies that for all $x,y\in [0,\infty)$, there exists
an $M\geq 1$, such that for all $k\geq 0$,
\begin{eqnarray}\label{e=distest}
\frac1M &\leq & \frac{(c^k)'(x)}{(c^k)'(y)} \leq M.
\end{eqnarray}
Assume without loss of generality that $n>m$ and suppose that there exists
a $C^1$ diffeomorphism $h\co [0,\infty)\to [0,\infty)$ such that
$hv_1(x) = v_2(h(x))$, for all $x\in[0,\infty)$.  Let $H(x) = h(x^{1/m})^n$.
Note that the $C^1$ function
$H\co [0,\infty)\to [0,\infty)$ has the following properties:
\begin{enumerate}
\item $H'(x) \geq 0$, for all $x\in[0,\infty)$, and $H'(x) = 0$ iff $x=0$;
\item for all $k\geq 0$, $H\circ c^k = c^k\circ H$.
\end{enumerate}
Then (2) implies that for every $x\in [0,\infty)$:
\begin{eqnarray*}
H'(x) &=& H'(c^k(x)) \frac{(c^k)'(x)}{(c^k)'(H(x))}
\end{eqnarray*}
for all $k\geq 0$.  But (\ref{e=distest}) implies that  
${(c^k)'(x)}/{(c^k)'(H(x))}$ is bounded independently of $k$, so that
$H'(x) = \lim_{k\to\infty}H'(c^k(x)) = 0$, contradicting property (1).
\eproof

\begin{cor}\label{c=sigdet} Let $G$ and $H$ be infinite subgroups of $\widehat{{\rm Aff}}^{{\bf s}_1}(\R)$ and   $\widehat{{\rm Aff}}^{{\bf s}_2}(\R)$, respectively, for some
${\bf s}_1, {\bf s}_2\in \overline{\mathcal S}$. If there exists 
$\alpha\in\Diff^1(S^1)$ such that $\alpha G \alpha^{-1} = H$, then ${\bf s}_1 = {\bf s}_2$.
%and $\alpha = \widehat{{\rm id}}(\zeta)$, for some $\zeta\in {\rm Stab}_{D_d}({\bf s}_1)$.

If $G < \widehat{{\rm Aff}}_+^{{\bf s}_1}(\R)$ and   
$H< \widehat{{\rm Aff}}_+^{{\bf s}_2}(\R)$, with ${\bf s}_1, {\bf s}_2\in \overline{\mathcal S}_+$, and  there exists $\alpha\in\Diff^1_+(S^1)$ 
such that $\alpha G \alpha^{-1} = H$, then ${\bf s}_1 = {\bf s}_2$.
%and $\alpha = \widehat{{\rm id}}(\zeta)$, for some 
%$\zeta\in {\rm Stab}_{C_d}({\bf s}_1)$.
\end{cor}

\proof  Let $G,H$ and $\alpha$ be given.  Note that ${\bf s}_1$ and
${\bf s}_2$ must have the same length $2d$, since the global finite invariant sets
of $G$ and $H$ must be isomorphic.  Let $g,h$ be elements of $G$ and $H$ with 
rotation number $0$ such that $h = \alpha g \alpha^{-1}$.  Since dilations
have twice as many fixed points in $\RP$ as translations, if $g$ is
a ramified lift of a translation, then so is $h$.  Assume that 
$g = \hat{S}(\pi_{{\bf s}_1}, id)$ and $h = \hat{T}(\pi_{{\bf s}_2}, id)$, where 
$S\co x\mapsto x+ s$ and $T\co x\mapsto x+t$ are translations with $s,t>0$.
 Let $q_1,\ldots, q_d$ and $\alpha(q_1), \ldots, \alpha(q_d)$ be the preimages
of $\infty$ under $\pi_{{\bf s}_1}$ and $\pi_{{\bf s}_2}$, respectively.  
In a neighborhood of $q_i$, the map $g$ is conjugate 
to $x\mapsto (S(x)^{m_i})^{1/m_i}$ and in a neighborhood of $\alpha(q_i)$, $h$ 
is conjugate 
to $x\mapsto (T(x)^{n_i})^{1/n_i}$, where $m_i = s(q_i)$ and $n_i = s(\alpha(q_i))$.
Since $S$ is a $C^2$ contraction in a neighborhood of $\infty$ and $T$ is conjugate to $S$,
it follows from Lemma~\ref{l=c2contr} that $m_i=n_i$ for $1\leq i\leq d$, which implies
that ${\bf s_1} = {\bf s}_2$.  

Suppose instead that $g$ is a ramified lift of a map in ${\rm Aff}(\R)$
conjugate to the dilation
$D\co  x\mapsto a x$, for some $a>1$.  Since $g$ must have $d$ fixed points
with derivative $a$, so must $h$, and so  
$h$ is also a ramified lift of a map in ${\rm Aff}(\R)$
conjugate to $D$.
Around $\infty$, the map $D$ is a $C^2$ contraction, and the same proof as above 
shows that ${\bf s}_1 = {\bf s}_2$.

The proof in the orientation-preserving case is analogous.
\eproof

\noindent
{\bf Proof of Proposition~\ref{p=c1conj}}\qua Let 
$\rho_n\co  \BS(1,n) \rightarrow \Diff^\omega(S^1)$ be
the standard representation.
Suppose that $\hat{\rho}_n(\pi_{\s_1}, h_1)$ and $\hat{\rho}_n(\pi_{\s_2}, h_2)\in 
\mathcal{V}$ are conjugate
by $\alpha\in \Diff^1(S^1)$, where ${\bf s}_1, {\bf s}_2\in {\mathcal S}$.
%:  $$\alpha \hat{\rho}_n(\pi_{\s_1}, h_1) \alpha^{-1} = 
%\hat{\rho}_n(\pi_{\s_2}, h_2).$$
It follows from Corollary~\ref{c=sigdet} that ${\bf s}_1= {\bf s_s}$.

We next show $h_1=h_2$.
Let $\gamma\in \BS(1,n)$ and let $k = \rho_n(\gamma)$.
Let $k_1 = \hat{k}(\pi_{\s_1}, h_1(\gamma))$ and 
$k_2 = \hat{k}(\pi_{\s_2}, h_2(\gamma))$
Then for all $q\in\pi_{\s_1}^{-1}(\infty)$,
we have:
$$\alpha h_1(\gamma)(q) = \alpha k_1(q) = k_2(\alpha(q)) = h_2(\gamma)(\alpha(q)),$$
and for all $e\in{\mathcal E}(\pi_1)$,
$$\alpha h_1(\gamma)(e) = \alpha k_1(e) = k_2(\alpha(e)) = h_2(\gamma)(\alpha(e)).$$
Since $\alpha (\pi_{\s_1}^{-1}(\infty)) = \pi_{\s_1}^{-1}(\infty)$ and
$\alpha ({\mathcal E}(\pi_1)) ={\mathcal E}(\pi_1)$, it follows that 
$h_1(\gamma) = h_2(\gamma)$.  So $\alpha h_1 = h_2 \alpha$.  Recall that each element 
$\widehat{\rho_n}(\pi_{\s}, h)$
of $\mathcal V$ is given by a  signature vector $\s \in \overline{\mathcal S}$ and a
representative $h$ of a
conjugacy class  in Hom$(\BS(1,n),$ Stab$_{D_d}(\s))$.  So $h_1 = h_2$.  
\eproof

\section{Proof of Theorems~\ref{t=main} and \ref{t=maincr}}\label{p=mainproofs}

The construction behind this proof is very simple.  We are given a $C^r$
representation $\rho$ of $\BS(1,n)$.  Using elementary arguments, we are
reduced to the case where $f=\rho(a)$ and $g=\rho(b)$ 
have a common finite invariant
set, the set of periodic orbits of $g$.  
Assume that the rotation numbers of $f$ and $g$ are both $0$.
Using the results from
Section~\ref{s=schwarz}, we obtain a local characterization of $f$ and
$g$ on the intervals between the common fixed points.  On each of these intervals,
$f$ is conjugate to the dilation $x \mapsto nx$ and $g$ is conjugate to the
translation $x \mapsto x+1$.  Gluing together the conjugating maps gives us a
$C^r$ ramified covering map over $\infty$.    Hence
$\rho$ is a $C^r$ ramified lift of the standard
representation. Proposition~\ref{p=ratramcovexists} implies
that there is a rational ramified cover with the same signature
as the the given $C^r$ ramified cover. 
Lemma~\ref{l=crconj} implies that $\rho$ is $C^r$
conjugate to a ramified lift of $\rho_n$ under the rational ramified cover.
It remains
to handle the case where the rotation numbers of $f$ and $g$ are not
$0$,
but this is fairly simple to do, since the elements of the 
standard representation
embed in analytic vector fields. We now give the complete proof.

Let $\rho\co  \BS(1,n)\to \Diff^r(S^1)$ be a
representation, where $r\in [2,\infty]$, or $r = \omega$.  If $r < \infty$, 
assume that
$\sigma(\rho) \leq \left(\frac{1}{n}\right)^{\frac{1}{r-1}}.$
If $r = \infty$, we assume that $\sigma(\rho) < 1.$

Let $f=\rho(a)$ and
$g=\rho(b)$, where $aba^{-1}=b^n$. Since $g$ is conjugate
to $g^n$, it follows that $\tau(g) = \pm \tau(g^n) = \pm n\tau(g)$,
where $\tau(h)$ denotes the rotation number of $h\in{\rm Homeo }(S^1)$.
Hence $g$ has rational rotation number.

\begin{lemma}\label{l=perpres} $f$ preserves the set of periodic points of $g$.
\end{lemma}

\proof  This follows from the relation $fg=g^nf$.  If $g^k(q) = q$, then
$g^{nk} (f(q)) = fg^k(q) = f(q)$.
So $f(q)$ is also periodic for $g$.
 \eproof

Suppose that $\tau(f)$ is irrational.  Then by Lemma~\ref{l=perpres},
the periodic points of $g$
are dense in $S^1$, which implies that $g^k = id$, for some $k\leq n+1$.
This implies that conclusion $(1)$ of Theorem~\ref{t=main} holds.  

Suppose, on the other hand, that $\tau(f)$ is rational.  Choose $l$
so that $g^l$ and $f^l$ are both orientation-preserving
and both have rotation number $0$.  %Note that $l\leq \max\{2n-2, n+1\}$. 
Then $f^l$ leaves ${\rm Fix}(g^l)$ invariant.  
Choose $p\in{\rm Fix}(g^l)$.  Any
accumulation of $\{f^{ln}(p)\}$ must be a fixed point for $f^l$ and for $g^l$. 
We have shown:

\begin{lemma} $f^l$ and $g^l$ have a common fixed point.
\end{lemma}

Note that the fixed points for $f^l$ are isolated;  if $f$ is not
analytic, then $\sigma(\rho) < 1$, which implies that the fixed points
for $f^l$ are hyperbolic.  Let $w_1 < w_2 < \ldots < w_k$ be the set of fixed points 
of $f^l$.  We will see that if $g^l$ is not the identity map, then the set of
fixed points for $g^l$ is exactly equal to the set of sinks for $f^l$.  

\begin{lemma}\label{l=sinkorid}  
If $g^l(w_i)=w_i$ and $(f^l)^\prime(w_i) > 1$, then $g^l = id$ on $[w_{i-1}, w_{i+1}]$.
\end{lemma}

\proof  Suppose that $(f^l)^\prime(w_i) = \lambda > 1$, and let $\alpha\co [w_i, w_{i+1})
\rightarrow [0,\infty)$ be a $C^1$ linearizing diffeomorphism such that
$\alpha f^l \alpha^{-1}(x) = \lambda x$ for all $x \in [0, \infty)$.  Let $F = 
\alpha f^l \alpha^{-1}$, and let $G = \alpha g^l \alpha^{-1}$.  If $g^l \neq id$ on 
$[w_i, w_{i+1})$, then there is a point $x_0 \in [0,\infty)$ such that
$G(x_0) \neq x_0$.  Let $x_0$ be any such point.  We may assume that $G^k(x_0) \rightarrow c$ as $k \rightarrow
\infty$, for some $c < \infty$, because this will be true for either $G$ or $G^{-1}$.
Since $GF^{-k} = F^{-k}G^{n^k}$ for all $k \in \N$, it follows that
\begin{eqnarray*}
G^\prime(F^{-k}(x_0)) &=& \frac{(F^{-k})^\prime(G^{n^k}(x_0))}{(F^{-k})^\prime(x_0)} \;
(G^{n^k})^\prime(x_0),
\end{eqnarray*}
for all $k \in \N$.  But since $G^\prime(0) = 1$ (by Lemma~\ref{l=derivs}),
this means that $(G^{n^k})^\prime(x_0)$ $\to 1$,
as $k\to \infty$ (or $k\to -\infty$), for every point $x_0$ that is not
fixed by $G$.  Since $G$ is not the identity, this is not
possible. Hence
$g=id$ on $[w_i, w_{i+1}]$. A similar argument shows that $g = id$ on
$[w_{i-1}, w_i]$,  \eproof

\begin{cor}\label{c=fixorid} If $g^l$ has a fixed point in the interval $(w_i, w_{i+1})$, then 
$g^l = id$ on $[w_i, w_{i+1}]$.  That is, $\partial{\rm Fix}(g^l)
\subseteq {\rm Fix}(f^l)$.
\end{cor}

\proof  Suppose that $g^l(p) = p$ for some $p \in (w_i, w_{i+1})$, and suppose that
$f^{kl}(p) \rightarrow w_{i}$ as $k \rightarrow -\infty$.  By Lemma~\ref{l=perpres},
$f^{lk}(p)$ is periodic for $g^l$ for all $k \in \Z$.  Since $g^l$ is
an orientation preserving circle diffeomorphism with a fixed point,
$f^{lk}(p)$ is a fixed point of $g^l$ for all $k$. By continuity, $w_i$ is a common
fixed point for $f^l$ and $g^l$.  Since $(f^l)^\prime(w_i) > 1$, 
Lemma~\ref{l=sinkorid} implies that $g^l = id$ on $[w_i, w_{i+1}]$.
Similarly, if $f^{kl}(p) \rightarrow w_{i+1}$ as $k \rightarrow -\infty$, then 
$g = id$ on $[w_i, w_{i+1}]$.
\eproof

This has the immediate corollary:

\begin{cor}\label{c=fixcomp} $f^l$ fixes every component of $S^1\setminus {\rm Fix}(g^l)$.
\end{cor}

\noindent
{\bf Remark}\qua Corollary~\ref{c=fixcomp} also follows from 
Theorem~\ref{t=navas}.  We have given a different proof here since we
will need Lemma~\ref{l=sinkorid} for the proof of Lemma~\ref{l=local}.

Let $-\infty\leq q_1< q_2 < \cdots < q_d <\infty $ be the elements of $\partial
{\rm Fix}(g^l)$.  
 
\begin{lemma}\label{l=local}  On each interval $(q_{i-1}, q_{i}]$, either
$g^l = id$, or there is a $C^r$ map $\alpha_i\co (q_{i-1}, q_{i}] \rightarrow 
(-\infty, \infty]$ such that
\begin{enumerate}
\item  $\alpha_i$ 
 conjugates $f^l$ to the map $x \mapsto n^l x$, and 
conjugates $g^l$ to the map $x \mapsto x+1$;
\item $\alpha_i\vert_{(q_{i-1}, q_{i})}$ is a $C^r$ diffeomorphism onto
 $(-\infty, \infty)$ 
\item For all $p$ in a neighborhood of $q_i$,  
$$\alpha_i(p) = o_i h(p)^s$$
where $h$ is a 
$C^r$ orientation-preserving diffeomorphism onto a neighborhood
of $\infty$, $1\leq s < r$, and $o_i \in \{\pm1\}$.  
\end{enumerate}
\end{lemma}

\proof  This follows from Proposition~\ref{p=schwarz}.  Note that we can
apply Proposition~\ref{p=schwarz} in this setting since we know that if $g^l
\neq id$ on $(q_{i-1}, q_i]$, then $(f^l)^\prime(q_i) \leq 1$ 
(by Lemma~\ref{l=sinkorid}).  By our assumptions on $\sigma(\rho)$,
if $2\leq r < \infty$, then $(f^l)^\prime(q_i) \leq (\frac{1}{n})^\frac{1}{r-1}$,
and if $r = \infty$, then  $(f^l)^\prime(q_i) <1$.  Therefore one of the assumptions
(A)--(C) of Proposition~\ref{p=schwarz} will hold.  \eproof

\begin{cor}\label{l=finfix} 
Either $g^l=id$, or $\partial {\rm Fix}(g^l) = {\rm Fix}(g^l) =
\{q_1,\ldots, q_d\}.$ 
\end{cor}

\proof Assume that
$\partial {\rm Fix}(g^l) = \{q_1,\ldots, q_d\} \neq {\rm Fix}(g^l)$,
but $g^l\neq id$.
Then there is an interval $[q_{i-1},q_i]$ on which $g^l=id$
but where $g^l \neq id$
on $[q_{i},q_{i+1}]$.  By Lemma~\ref{l=schwarz},  either $g^l$ or
$g^{-l}$ is $C^r$ conjugate in a neighborhood $[q_i, p)$ to the map
$x \mapsto x/(1-x^s)^\frac{1}{s}$ for some integer $1\leq s<r$.
But this map is not $r$--flat at $x=0$, so $g^l$ is not $C^r$ at $q_i$,
a contradiction.  
\eproof

\begin{cor}\label{c=extend} If $g^l\neq id$, then the map 
$\pi\co  S^1\to\RP$ defined by:
$$\pi (p) = \alpha_i(p),\,\quad\hbox{for } p\in (q_{i-1}, q_{i}]$$
is a $C^r$ ramified covering map over $\infty$,
$f^l$ is a $\pi$--ramified lift of 
$x \mapsto n^l x$, and $g^l$ is a $\pi$--ramified lift of the
map $x \mapsto x+1$.
\end{cor}

\proof  Let $q_i \in$ Fix$(g^l)$.  
Applying Lemma~\ref{l=local} to the interval $[q_i, q_{i+1})$,
we obtain a map $\overline\alpha_{i+1}\co  [q_i, q_{i+1}) \rightarrow [-\infty, \infty)$
which is a $C^r$ diffeomorphism on $(q_i, q_{i+1})$, and which is a power
of a $C^r$ diffeomorphism in a (right) neighborhood of $q_i$; 
$\overline\alpha_{i+1}(p) = h(p)^s$ for $p$ near $q_i$, for some $C^r$ diffeomorphism $h$ and some
integer $1 \leq s < r$.  Similarly, on the interval $(q_{i-1}, q_i]$
there is a map $\alpha_i\co  (q_{i-1}, q_i] \rightarrow 
(\infty, -\infty]$ which is a power of a diffeomorphism in a (left)
neighborhood of $q_i$; $\alpha_i(p) = h_*(p)^{s_*}$ near $q_i$.
We will show that $s=s_*$, and that the diffeomorphisms $h$ and $h_*$
glue together to give a $C^r$ diffeomorphism in a neighborhood of
$q_i$.  This will prove that the map
$$\pi_i (p) = 
\begin{cases} \alpha_i(p),\,\quad\hbox{for } p\in (q_{i-1}, q_{i}] \\
\overline{\alpha}_{i+1}(p) ,\,\quad\hbox{for } p\in [q_{i}, q_{i+1})
\end{cases}$$ is the restriction to $(q_{i-1}, q_{i+1})$ of a $C^r$ ramified
covering map over $\infty$.  By construction, the restrictions of $f^l$
and $g^l$ to $(q_{i-1}, q_{i+1})$ are $\pi_i$--ramified lifts of the maps 
$x \mapsto n^l x$ and $x \mapsto x+1$ respectively.

The diffeomorphism $1/h$ maps $q_i$ to $0$, and conjugates $g^l$ to the
map $x\mapsto x/(1 + x^s)^{1/s}$.  Similarly, $1/h_*$ conjugates $g^l$
to $x\mapsto x/(1 + x^{s_*})^{1/{s_*}}$.  Since $g$ is 
$C^r$, we must have $s=s_*$.   Both $1/h$ and $1/h_*$ are
linearizing maps for $f^l$ at $q_i$, and it is not hard to see that they
define a $C^r$ diffeomorphism $H$ in a neighborhood of $q_i$.  Therefore
$h$ and $h_*$ glue together to give a $C^r$ diffeomorphism $1/H$ in a neighborhood
of $q_i$.

It remains to show that $\pi_i = \pi_{i+1}$ on $(q_i, q_{i+1})$. 
Since the restriction of both of these maps
to $(q_i, q_{i+1})$ are  diffeomorphisms which linearize $f^l$, they are the 
same up to a constant multiple.  There is a unique point $x_0 \in  (q_{i-1}, q_{i})$ 
satisfying $f^l(x_0) =  g^{n^l-l}(x_0)$; -- this is the point $x_0 = g^l(y)$, where
$y$ is the unique fixed point for $f^l$ in $(q_i, q_{i+1})$. 
Both $\pi_i$ and $\pi_{i+1}$ send the point $x_0$ 
to the same point $1\in \R$.  So we have
$\pi_i = \pi_{i+1}$ on $(q_i, q_{i+1})$.
\eproof

It follows from Lemma~\ref{l=crconj} that the representation of $\BS(1, n^l)$ 
generated by $f^l$ and $g^l$ is $C^r$ conjugate to an element of $\mathcal{V}$.
In the remainder of this section, we will show that the diffeomorphisms $f$ and $g$
are $C^r$ ramified lifts of the generators of the standard action of $\BS(1,n)$ on $S^1$,
hence the representation they generate is also $C^r$ conjugate to an element of 
$\mathcal{V}$.  We begin with some lemmas about ramified lifts of flows on $S^1$.

\begin{lemma}\label{l=liftflow}  Let $\varphi\co  S^1 \rightarrow S^1$ be a
$C^r$ flow with a fixed point at $p$, and let $\pi\co  S^1 \rightarrow
S^1$ be a $C^r$ ramified covering map over $p$.  Let 
$F = \widehat{\varphi^1}(\pi, id)$ be the $\pi$--ramified lift of the time-$1$
map $\varphi^1$ with rotation number zero.  Then $F$ embeds as the time-$1$
map of a $C^r$ flow $F^t$ on $S^1$, and for all $t \in \R$, 
$F^t = \widehat{\varphi^t}(\pi, id)$.
\end{lemma}

\proof  By Lemma~\ref{l=liftdiff}, given any $t \in \R$ there is a unique $(\pi, id)$ -
ramified lift of $\varphi^t$, $F^t := \widehat{\varphi^t}(\pi, id)$.
Lemma~\ref{l=liftdiff2} imples that
$F^t \circ F^s = F^{s+t}  = F^s \circ F^t$ for all $s, t \in \R$.  

  Let $X$ be the $C^{r-1}$ vector field that generates $\varphi$, and
let $\hat X$ be the lift of $X$ under $\pi$. This  vector field is clearly 
$C^{r-1}$ on 
$S^1\setminus \pi^{-1}(p)$ and clearly generates the flow $F^t$ on $S^1$.
In a neighborhood of $q_i \in \pi^{-1}(p)$, $\hat X$ takes the form
$$\hat X (q) = d_{\pi (q)}\pi^{-1} X(\pi q),$$ and $\pi$ takes the form
$\pi(x) = x^s.$  A straightforward calculation shows that
the vector field $\hat X$ is $C^{r-1}$.  Similarly, $\hat X$ is analytic if $X$
and $\pi$ are. This completes the proof.
\eproof

\begin{lemma}\label{l=cflow}  Let $F\co  S^1 \rightarrow S^1$ be the time-$1$ map of a $C^r$ flow $F^t$,
where $r \geq 2$.  Suppose that $F$ is not $r$--flat, and $\tau(F) = 0$.  If $G$ is a $C^r$
orientation preserving
diffeomorphism such that $FG=GF$, and if $\tau(G) = 0$, then $G = F^t$ for some $t \in \R$.
\end{lemma}

\proof 
Since $\tau(F) = 0$ and $F$  is not $r$--flat, $F$ has a finite set of fixed points.
Let $q_1 < \ldots < q_d$ be the elements of Fix$(F)$.  If $FG=GF$, then $G$ permutes the
fixed points of $F$, and since $G$ is orientation preserving and has
rotation number zero,
$G([q_i,q_{i+1}]) = [q_i,q_{i+1}]$ for all $q_i \in {\rm Fix}(F)$.
By Lemma~\ref{l=kopell},
on $[q_i, q_{i+1})$, $G = F^{t_i}$ for some $t_i \in \R$, and on $(q_i, q_{i+1}]$, 
$G = F^{s_i}$ for some $s_i$.  Clearly, $t_i = s_i$.  So for $1 \leq i \leq d$,
$$G|_{[q_i, q_{i+1}]} = F^{t_i}, \hbox{ for some } t_i \in \R.$$
If $F^\prime(q_i) \neq 1$ for some $q_i \in {\rm Fix}(F)$,  then since $G$ is $C^1$
at $q_i$, it follows that $t_i = t_{i-1}$.  If $F^\prime(q_i) = 1$, then in
local coordinates, 
identifying $q_i$ with $0$,
$$F(x) = F^1(x) = x + ax^k + o(x^k)$$ for some $a \neq 0$ and $k \leq r$.
Therefore
$$G(x) = \begin{cases}   x + t_{i-1}ax^k + o(x^k), & \hbox{ for } x \in (q_{i-1}, q_i] \\
 x + t_{i}ax^k + o(x^k), &  \hbox{ for } x \in [q_i, q_{i+1}).
\end{cases}$$
Since $k \leq r$ and $G$ is $C^r$, $t_i = t_{i-1}$.
\eproof

\begin{cor}\label{c=cflow}  Let $\pi\co  S^1 \rightarrow S^1$ be a $C^r$ ramified covering
map over $\infty$, and let $F = \hat{k}(\pi, id)$ be the 
$(\pi, id)$--ramified lift
of 
$k \in {\rm Aff}(\R)$, $k\neq id$.  Let ${\bf s}(\pi) = (s_1, \ldots, s_d,
o_1, \ldots, o_d)$,
where $s_i \leq r-1$ for $1 \leq i \leq d$.  By Lemma~\ref{l=liftflow}, $F$
embeds as the time-$1$ map of a $C^r$ flow $F^t$.  If $H\co  S^1 \rightarrow
S^1$ is a $C^r$ orientation preserving 
diffeomorphism such that $FH=HF$, and if $\tau(H) = 0$,
then $H = F^t$ for some $t \in \R$.
\end{cor}

\proof  By Lemma~\ref{l=cflow}, it is enough to show that $F$ is not $r$--flat.
In coordinates identifying a fixed point with $0$,
$$ F(x) = \frac{x}{(b+ax^{s})^\frac{1}{s}},$$ 
where either $a\neq 1$ or $b\neq 0$,
which is clearly not $r$--flat, if $s < r$.
\eproof

\begin{prop}\label{p=root}  
Let $F\in \Diff^r(S^1)$ be a diffeomorphism such that
$F^l$ is orientation-preserving and $\tau(F^l) = 0$, for some $l>0$. 
Suppose that $F^l= \widehat{k^l}(\pi, id)$ is a $C^r$ ramified lift of
$k^l \neq id$, 
where $k\in {\rm Aff}_+(\R)$, and suppose that
${\bf s}(\pi) = (s(q_1), \ldots, s(q_d), o_1, \ldots, o_d)$, where
$s(q_i)
 \leq r-1$ for $1 \leq i \leq d$.
Then either $F$ is a $\pi$--ramified lift of  $k$ or $F$ is a $\pi$--ramified lift of $-k$.
\end{prop}

\proof  
Let $\zeta \in D_d$
be such that $\zeta(q) = F(q)$ for all $q \in \pi^{-1}(\infty)$, and
$\zeta(e) = F(e)$ for all $e \in {\mathcal E}(\pi)$.

\begin{lemma}\label{l=stabhash}  $\zeta \in {\rm Stab }^{\#}_{D_d}(\s(\pi))$.
\end{lemma}

\proof  Given any $q \in \pi^{-1}(\infty)$, there is an interval $[q,p)$
and $C^r$ diffeomorphisms $h_1\co [q,p)\rightarrow [0,\infty)$ and
$h_2\co [F(q),F(p))\rightarrow [0,\infty)$ such that
$$h_1F^lh_1^{-1}(x) = [k^l(x^{s})]^\frac{1}{s}, \; \; \hbox{ and } \; \;
h_2F^lh_2^{-1}(x) = [k^l(x^{t})]^\frac{1}{t},$$ where $s = s(q)$ and
$t=s(F(q))$.
We can assume that $k^l$ is a contraction on $[0,\infty)$.  (If not, then
use $k^{-l}$ and $F^{-l}$).  Since $F^l|_{[q,p)}$ is conjugate by $F$ to 
$F^l|_{[F(q),F(p))}$,  Lemma~\ref{l=c2contr} implies that $s(q) =
s(F(q))$, and therefore $\zeta \in {\rm Stab}^{\#}_{D_d}(\s(\pi))$.
\eproof

By Lemma~\ref{l=stabhash}, either $\zeta(\s(\pi)) = \s(\pi)$, or
$\zeta(\s(\pi)) = I(\s(\pi))$.  If $\zeta(\s(\pi)) = \s(\pi)$, then let
$\alpha\co  S^1 \rightarrow S^1$ be the $(\pi, \zeta)$--ramified lift of the identity
map: $\alpha = \widehat{id}(\pi, \zeta)$.
By Lemma~\ref{l=liftdiff2}, $\alpha$
commutes with $F^l$.  So $F \alpha^{-1}$ commutes with $F^l$, and by 
construction, $F \alpha^{-1}$ fixes every 
interval $(q_i, q_{i+1}) \subseteq \pi^{-1}(\RP \setminus \{\infty\})$.
By Lemma~\ref{l=cflow}, $F^l$ embeds as the time-$1$ map of a $C^r$ flow, $F^t$,
and $F \alpha^{-1} = F^{t_0} = \widehat{\phi^{t_0}}(\pi, id)$ for some $t_0 \in \R$,
where $\phi$ is an analytic flow with $\phi^1 = k^l$. Therefore 
$F = \widehat{\phi^{t_0}}(\pi, id) \circ \widehat{id}(\pi, \zeta) = 
\widehat{\phi^{t_0}}(\pi,\zeta)$ (using Lemma~\ref{l=liftdiff2}).  It follows that
$t_0 = 1/l$, and therefore $F$ is the $(\pi, \zeta)$--ramified lift of the map $k$.

If $\zeta(\s(\pi)) = I(\s(\pi))$, then we let $\alpha = \widehat{-id}(\pi,\zeta)$,
the $(\pi, \zeta)$--ramified lift of $-id\co  x \rightarrow -x$.  
As above, $F \alpha^{-1}$ is the $(\pi, \zeta)$--ramified lift of 
$k$, and therefore $F = \psi^t \alpha = \widehat{-k}(\pi, \zeta)$.
\eproof

\begin{cor}  $f$ and $g$ are $C^r$
 ramified lifts under $\pi$ of the generators of the standard action of $\BS(1,n)$ on $S^1$.
\end{cor}

\proof The standard representation 
$\rho_{n^l}\co  \BS(1, n^l) \rightarrow {\rm Diff }^\omega(\RP)$
is analytically conjugate to the representation
$\kappa\co  \BS(1, n^l) \rightarrow {\rm Diff }^\omega (\RP)$ with generators
$\kappa(a^l)\co  x \mapsto n^l x$ and $\kappa(b^l)\co  x \mapsto x + l$.  So there
is a $C^r$ ramified
covering map  $\pi\co  S^1 \rightarrow S^1$ over $p$ such that 
$f^l = \widehat{\kappa(a^l)}(\pi, id)$ and $g^l =
\widehat{\kappa(b^l)}(\pi, id)$.  By Proposition~\ref{p=root},
either $f$ is  a ramified lift of $\rho_n(a)\co  x \mapsto nx$,
or $f$ is a ramified lift of $-\rho_n(a)\co  x \mapsto -nx$.  Similarly, $g$ is
either a ramified lift of $\rho_n(b)\co  x \mapsto x+1$, or a ramified lift
of $-\rho_n(b)\co  x \mapsto -x-1$.  Since $f$ and $g$ satisfy the relation
$fgf^{-1} = g^n$, the maps that they are lifted from must also satisfy this relation.
Given this requirement,
the only possibility is that $f$ is a $\pi$--ramified lift of $\rho_n(a)$
and $g$ is a $\pi$--ramified lift of $\rho_n(b)$.
\eproof

Since the generators $\rho(a) = f$ and $\rho(b) = g$ of the
representation $\rho$ are ramified lifts under $\pi$ of the generators
$\rho_n(a)$ and $\rho_n(b)$, respectively, of $\rho_n$, it follows that,
for every $\gamma\in \BS(1, n)$, there exists a unique $h(\gamma)\in D_d$ (or
in $C_d$ if $\rho$ is orientation-preserving) such that:
$$\rho(\gamma) = \widehat{\rho_n(\gamma)}(\pi, h(\gamma)).$$ 
Since $\rho(\gamma_1\gamma_2) = \rho(\gamma_1)\rho(\gamma_2)$, it
follows that $h\co  \BS(1,n)\to D_d$ ($C_d$) is a homomorphism. 
Finally, note that $h$ must take values in ${\rm Stab}_{D_d}({\bf s})$
(or ${\rm Stab}_{C_d}({\bf s})$, if $\rho$ is orientation-preserving).

This concludes the proof of Theorems~\ref{t=main} and
\ref{t=maincr}.\eproof

Finally, we sketch the proof of Theorem~\ref{t=navas}.

\medskip

\noindent
{\bf Sketch of proof of Theorem~\ref{t=navas}}\qua Let $\rho$ be a $C^r$
representation of $\BS(1,n)$, with $r\geq 2$, let $f=\rho(a)$ and
$g=\rho(b)$.  We may assume that $f$
has rational rotation number.  By taking powers of the elements of
$\BS(1,n)$, we may assume that both
$f$ and $g$ have rotation number $0$. Assume that $g$ is not the identity map.

Let $J$ be a component of the complement of ${\rm Fix}(g)$.  Using a distortion estimate and 
the group relation one shows that $J$ must be fixed by $f$, as
follows.  Otherwise, the $f$--orbit of $J$ must accumulate at both ends on a
fixed point of $f$.  The standard $C^2$ distortion estimate shows that there
is an $M>1$ such that for all $x,y$ in same component of the $f$--orbit of $J$, and
for all $k\in\Z$,
$$\frac{1}{M} < | \frac{(f^k)'(x)}{(f^k)'(y)}| < M.$$  
But, for all $k\in N$, we have that $f^k g f^{-k} = g^{n^k}$.  Hence, for all $p \in J$,
we have:
$$(g^{n^k})'(p) = g'(y)  \frac{(f^k)' (g(y))} {(f^k)' (y)},$$
where $y = f^{-k}(p)$.  Note that $y$ and $g(y)$ lie in the same component 
$f^{-k}(J)$, and $g'(y)$ is uniformly bounded.  This implies that for all $p \in
J$ and all $k\in\N$, $(g^{n^k})'(p)$ is bounded, so that $g=id$ on $J$, a contradiction.

So $f$ fixes each component of the complement of ${\rm Fix}(g)$.  Let $J$ be such a
component.  Since $g$ has no fixed points on $J$, $g$ embeds in a $C^1$ flow $g^t$,
defined on $J$ minus one of its endpoints, that is $C^r$ in the interior
of $J$ (see, eg \cite{ZL}).  
Furthermore, for all $t$, $fg^tf^{-1} = g^{nt}$ (this follows from
Kopell's lemma).  Fixing some point $p$ in the interior of $J$, this flow 
defines a $C^r$ diffeomorphism between the real line and the interior of $J$, sending 
$t \in \R$ to $g^t(p)\in J$.  Conjugating by this diffeomorphism, $g^t$ is sent to a 
translation by $t$, and $f$ is sent to a diffeomorphism $F$ satisfying $F(x+t) = F(x)+tn$,
for all $t,x \in \R$.  But this means that $F'(x) = n$ for all $x \in
\R$. Up to an affine change of coordinates, $g$ is conjugate on $J$ to
$x\mapsto x+1$ and $f$ is conjugate to $x\mapsto nx$. \eproof

\section{Proof of Proposition~\ref{p=deform}}
Let $\rho\co \BS(1,n)\to \Diff^\omega(S^1)$ be a $\pi$--ramified lift of the
standard representation $\rho_n$ with $\sigma(\rho) =\left( \frac{1}{n}
\right)^\frac{1}{r-1}$, for some $r\geq 2$.  Let ${\mathcal Q}$
be the
set of all points $q\in \pi^{-1}(\infty)$ satisfying $s(q) = r-1$;
this set is nonempty since  $\sigma(\rho) =\left( \frac{1}{n}
\right)^\frac{1}{r-1}$.  Lemma~\ref{l=normal} implies that,
 in a neighborhood of $q$, $\pi\co  x \mapsto  x^{r-1}$,
in the appropriate coordinates identifying $q$ with $0$.

For $t\in (-1,1)$ we deform $\pi$  to
obtain a $C^{r-1+t^2}$ map $\pi_t\co S^1\to S^1$ with the following properties:
\begin{itemize}
\item $\pi_0 = \pi$ and $\pi^{-1}(\infty) = \pi_t^{-1}(\infty)$, for all $t$;
\item $\pi_t\vert_{S^1\setminus \pi_t^{-1}(\infty)}$ is a $C^{\infty}$ covering
map onto its image;
\item about each $q\in \pi^{-1}(\infty)\setminus{\mathcal Q}$, $\pi_t$
is locally equal to $\pi$;
\item about each $q \in {\mathcal Q}$, $\pi_t$ is locally 
$x\mapsto  x^{r-1+t^2}$, in the same charts identifying $q$ with $0$
described above.
\end{itemize}

A slight modification of the proof of
Proposition~\ref{p=cr} also shows that $\rho_n$ has
a lift to a $C^r$ representation $\rho_t\co \BS(1,n)\to \Diff^r(S^1)$ so
that the following diagram commutes, for all $\gamma\in \BS(1,n)$:

\fig{
\object {a}  0,8 {S^1};
\object {b}  8,8 {S^1};
\object {c}  0,4  {S^1};
\object {d}  8,4  {S^1};
\arrow  {a}  {b}  {\rho_t(\gamma)};
\arrow  {c}  {d}  {\rho_n(\gamma)};
\arrow  {a}  {c}  {\pi_t};
\arrow  {b}  {d}  {\pi_t};}
(One merely needs to check that the integer $j$ in the proof
of Lemma~\ref{l=liftdiff} can be replaced by the real number
$r-1+t^2$).

Notice that $\rho_t$ has the property that $\sigma(\rho_t) =
\left(\frac{1}{n}\right)^\frac{1}{r-1+t^2}$, so that $\rho_s$ is not $C^1$
conjugate to $\rho_t$ unles $s=t$.  One can further modify this
construction by replacing the points of ${\mathcal Q}$ by intervals
of length $\eps_t$, extending $\rho_t(b)$ isometrically
across these intervals, and extending $\rho_t(a)$
in an arbitrary $C^r$ fashion to these intervals.
 Since $\rho_t(b)$
is $r$--flat (by Lemma~\ref{l=derivs}) on ${\mathcal Q}$ for $t\neq 0$ and
$r-1$ flat for $t=0$, 
the representation $\rho_t$ is $C^r$
and varies $C^{r-1}$ continuously in 
$t$ if we choose $\eps_t\to 0$ as $t\to 0$.  
In this way, one can create uncountably many deformations of $\rho$.
(Note that, in essence, we have deformed $\pi$ to obtain a ``broken
$C^r$ ramified cover'' \`a la Theorem~\ref{t=navas}).
\eproof

 \section{Proof of Theorem~\ref{t=solv}}
Let $r \in \{\infty, \omega\}$, and let
$G < \Diff^r(S^1)$ be a solvable group without infinitely flat
elements.  Suppose that $G^m:= \{g^m: g \in G\}$ is not abelian, for
any $m \in \Z$. We begin by showing that the group $G^2 < {\rm Diff}_+^r(S^1)$ 
has a finite set of points that is globally invariant.

\begin{lemma}\label{l=G2}  $G^2$ contains a non-trivial normal abelian subgroup
$N$, such that $N$ contains an element of infinite order. There is an integer
$d>0$, and a finite set $\{q_1,\ldots, q_d\}$,
with $q_1<q_2<\cdots < q_d$, such that:
\begin{enumerate}
\item for all $f\in G^2$, $\tau(f^d) = 0$ and $f\{q_1,\ldots, q_d\} =
\{q_1,\ldots, q_d\}$;
\item for all $g\in N$, either $g^d = id$ or ${\rm Fix}(g^d) =
\{q_1,\ldots, q_d\}$.
\end{enumerate} 
\end{lemma}

\proof Note that $G^2$ is a solvable group, and every diffeomorphism
in $G^2$ is orientation preserving.  Let 
$$G^2= G_0 > G_1 > \ldots > G_n > G_{n+1} = \{id\}$$
be the derived series for $G^2$, and let $N = G_n$ be the terminal
subgroup in this series.  Recall that $N$ is a normal abelian subgroup 
of $G^2$.  We first show that $N$ contains an element of infinite order.
We will use the following result of Ghys (\cite{Gh2}  Proposition 6.17):

\begin{lemma}\label{l=rotno1}
If $H \subset {\rm Homeo}_+(S^1)$ is solvable, then the rotation
number $\tau\co  H$ $ \rightarrow \R / \Z$ is a homomorphism.
\end{lemma}

Suppose that every diffeomorphism in $N$ has finite order.  Since $G^2 
\neq N$ (because $G^2$ is not abelian), $G_{n-1}$ cannot be abelian -- 
if it were, then $G_n$ would be trivial.  Suppose that $f,h \in
G_{n-1}$.  By Lemma~\ref{l=rotno1}, $\tau(fhf^{-1}h^{-1}) = 0$. 
But $fhf^{-1}h^{-1}$ is orientation
preserving and has finite order, since $fhf^{-1}h^{-1} \in G_n$.
Therefore $fhf^{-1}h^{-1} = id$, and $G_{n-1}$ is abelian, a
contradiction.
So $N$ contains a diffeomorphism with infinite order.

If $\tau(g)$ is irrational, for some orientation-preserving $g\in N$, then the elements
of $N$ are simultaneously conjugate to rotations.  But, since
$N$ is normal in $G^2$, this implies that
the elements of $G^2$ are simultanously conjugate to rotations, which 
implies that $G^2$ is abelian, a contradiction.

Hence $\tau(g)\in \Q/\Z$, for every $g\in N$.  Note that every $g \in N$ either
has finite order, or a finite set of periodic points: 
if Fix$(g^l)$ is infinite, for
some integer $l \neq 0$, then there is a point $q \in {\rm Fix}(g^l)$ that
is an accumulation point for a sequence $\{q_i \} \subset$ Fix$(g^l)$.  But this
implies that $g^l$ is infinitely flat at $q$, and therefore $g^l =
id$.  

Hence there exists $g \in N$ with infinite order and a finite fixed set,
${\rm Fix}(g)=\{q_1,\ldots, q_d\}$.  
%Then  $g^d$ is orientation-preserving,
%$\tau(g^d) = 0$ and ${\rm Fix}(g^d) = \{q_1,\ldots, q_d\}$.
If $h\in N$ is another element of $N$, then, since $h$ commutes with
$g$, it follows that $h(\{q_1,\ldots, q_d\})  = \{q_1,\ldots, q_d\}$,
and so $\tau(h^d) = 0$.  If the set of fixed points for $h^d$ is
infinite, then $h^d = id$, and if ${\rm Fix}(h^d)$ is finite, then 
${\rm Fix}(h^d)  = \{q_1,\ldots, q_d\}$.

Finally, let $f\in G^2 \setminus N$ and pick $g\in N$ satisfying
${\rm Fix}(g)  = \{q_1,\ldots, q_d\}$.  Then there exists
a $\overline g\in N$ such that $fgf^{-1} = \overline g$.  This implies
that $f({\rm Fix}(g)) = {\rm Fix}(\overline g)$; that is,
$f(\{q_1,\ldots, q_d\})  = \{q_1,\ldots, q_d\}$.  It follows that
$\tau(f^d) = 0$.  This completes the proof.  \eproof

Let $\{q_1,\ldots, q_d\}$ be given by the previous lemma,
labelled so that $-\infty\leq q_1 < q_2 <\cdots < q_d <\infty$, and
let $l=2d$. 
We will begin by working with the group $G^l$.
Note that every $g \in G^l$ is orientation-preserving, has zero rotation 
number, and 
fixes every point in the set $\{q_1, \ldots, q_d\}$.
Throughout this section, 
we will be working on the intervals $(q_i, q_{i+1})$, 
where we adopt the convention that $q_{d+1} = q_1$.

Let $M$ be a normal abelian subgroup of $G^l$ which contains an
element of infinite order.   For the rest of the proof, fix a
diffeomorphism $g \in M$ which has infinite order.

\begin{lemma} Let ${\mathcal C}(g) = \{f\in G^l\,\mid\, gf
= fg\}$.  Then ${\mathcal C}(g) \neq G^l$.
\end{lemma}

\proof  A proof  of this lemma is essentially contained in 
\cite{FaSh}. This lemma is implied by the following theorem, which is classical. 

\begin{theorem}[H{\"o}lder's Theorem]\label{t=h}  If a group of homomorphisms acts 
freely on $\R$, then it is abelian.
\end{theorem}
If $f \in \mathcal{C}(g)$, then Fix$(f)$ = Fix$(g)$.  So on every interval $(q_i, q_{i+1})$,
$1 \leq i \leq d$, no element of $\mathcal{C}(g)$ has a fixed point.  By Theorem~\ref{t=h},
the restriction of the action of $\mathcal{C}(g)$ to each interval $(q_i, q_{i+1})$ is
abelian.  Since $f(q_i) = q_i$ for all $f \in \mathcal{C}(g)$ and for all $q_i \in$
Fix$(g)$,  $\mathcal{C}(g)$ is an abelian subset of $G^l$.  But $G^l$ is not abelian, so
${\mathcal C}(g) \neq G^l$.  \eproof

\begin{lemma}\label{l=local2} 
Let $f\in G^l\setminus \mathcal{C}(g)$.  
Then for every interval $(q_{i-1}, q_{i}]$
there is a positive real number $\lambda_i$ and a $C^r$ map
$\alpha_i\co (q_{i-1},q_{i}]\to \RP$ with the following properties:
\begin{enumerate}
\item $\alpha_i g (p) = \alpha_i(p) + 1$ 
and $\alpha_i f (p) = \lambda_i \alpha_i(p)$;
\item $\alpha_i\vert_{(q_{i-1}, q_{i})}$ is a $C^r$ diffeomorphism onto
$(-\infty,\infty)$;
\item  there is an orientation-preserving $C^r$ diffeomorphism
$h_i$  from a neighborhood of $q_{i}$ to a neighborhood of $\infty$ 
and integers $s_i \in \{1,\ldots, r-1\}$, 
$o_i\in\{\pm 1\}$,
such that, for all $p$ in this neighborhood:
$$\alpha_i(p) = o_i  h_i (p)^{s_i}$$
\end{enumerate}
The same conclusions hold, with the same $\lambda_i$, $o_i$ 
and $\alpha_i\vert_{(q_{i-1}, q_{i})}$, but different local
(orientation-reversing) diffeomorphism $h_i^\ast$ and integer $s_i^\ast$, 
when $q_{i}$ is replaced by $q_{i-1}$ and
the interval $(q_{i-1}, q_{i}]$ is replaced by $[q_{i-1}, q_{i})$.

\end{lemma}

\remark To ensure that the conditions $\lambda_i >0$ (as opposed to
$\lambda_i\neq 0$)  hold in Lemma~\ref{l=local2}, it is necessary that
we chose $l$ to be even.

\proof We use the following fact, proved by Takens:

\begin{theorem}[\cite{Ta}, Theorem 4] Let $h\co [0,1)\to [0,1)$ be a $C^\infty$ diffeomorphism
with unique fixed point $0\in [0,1)$.  If $h$ is not infinitely flat, then there exists a unique $C^\infty$ vector field $X$ on $[0,1)$ such that
$h = h^1$, where $h^t$ is the flow generated by $X$.
\end{theorem}

For each $1 \leq i \leq d$, Let $g_i^t\co (q_{i-1},q_{i}]\to 
(q_{i-1},q_{i}]$ be the flow 
given by this theorem with $g_i^1 = g\vert_{(q_{i-1},q_{i}]}$.
If $f \in G^l \setminus \mathcal{C}(g)$,
then since $g \in M$, we have $fgf^{-1} \in M$, and therefore 
$fgf^{-1} \in \mathcal{C}(g)$.
By Lemma~\ref{l=kopell}, for $1 \leq i \leq d$, we must have
$$fgf^{-1} = g_i^{\lambda_i}$$ on $(q_{i-1}, q_{i}]$, 
for some $\lambda_i \in 
\R \setminus\{0\}$, $\lambda_i \neq 1$.
Note that, because $l$ is even, $\lambda_i$ must be positive, for all $i$.
So assumption (D) of Section~\ref{s=schwarz} holds in the interval
$(q_{i-1}, q_{i}]$ for each $q_i \in \{q_1, \ldots, q_d\}$.  

The same reasoning can
be applied to $[q_{i-1},q_{i})$, using a possibly different flow
$\tilde{g}_i^t$ and constant $\mu_i$ with 
$$fgf^{-1} = \tilde{g}_i^{\mu_i}.$$
Since $g_i^{\lambda_i}$ and $\tilde{g}_i^{\mu_{i}}$ coincide on
$(q_i, q_{i+1})$, it is not hard to see that we must have $\lambda_i =
\mu_{i}$. Now the result follows from
Proposition~\ref{p=schwarz}, as in the proof of Lemma~\ref{l=local}
and Corollary~\ref{c=extend}.
\eproof

\begin{cor} For every $f\in G^l\setminus {\mathcal C}(g)$,
there is a positive real number $\lambda = \lambda(f) \neq 1$ such that $\lambda_i =
\lambda$, for all $1\leq i\leq d$, where $\lambda_i$ is given by 
Lemma~\ref{l=local2}.  For every $i$, $s_i = s_{i+1}^\ast$, where addition is
mod $d$.
\end{cor}

\proof   Let $1 \leq i \leq d$. As in the proof of
Corollary~\ref{c=extend}, we have that 
$g$ is conjugate in a left neighborhood of
$q_i$ to  $x\mapsto x/(1 + x^{s_i})^{1/{s_i}}$
and $g$ is conjugate in a right neighborhood of $q_i$ to
$x\mapsto x/(1 + x^{s_{i+1}^\ast})^{1/{s_{i+1}^\ast}}$.  Since $g$ is 
$C^\infty$, we must have $s_i=s_{i+1}^\ast$. But then
$f$ is conjugate in a left neighborhood of
$q_i$ to  $x\mapsto x/\lambda_i^{1/s_i}$
and $f$ is conjugate in a right neighborhood of $q_i$ to
$x\mapsto x/(\lambda_{i+1})^{1/{s_{i+1}^\ast}}$.  It follows
that $\lambda_i = \lambda_{i+1}$.  Set $\lambda$ to be
this common value. \eproof

The proof of the next corollary is identical to the proof of Corollary~\ref{c=extend}.

\begin{cor}\label{l=a} For every $f\in G^l\setminus {\mathcal C}(g)$,
the map $\pi\co  S^1 \to\RP$ defined by:
$$\pi (p) = \alpha_i(p),\,\quad\hbox{for } p\in (q_{i-1}, q_{i}]$$
is a $C^r$ ramified cover with 
signature ${\bf s} = (s_1, \ldots s_d, o_1, \ldots, o_d)$. The diffeomorphism
$f$ is a $\pi$--ramified lift of the map $x \mapsto \lambda(f)x$, and
$g$ is a $\pi$--ramified lift of the map $x \mapsto x+1$.
\end{cor}

\begin{cor}\label{c=existsvarphi}
 $g$ embeds in a unique $C^r$ flow $g^t$,
with $g = g^1$.  The elements
of $\mathcal{C}(g)$ belong in the flow for $g$ and, for each $f \in G^l \setminus \mathcal{C}(g)$,
lie in the ramified lift
under $\pi_f$ of the translation group $\{x\mapsto x+\beta\,\mid\,
\beta\in \R\}$. That is, for any $h\in \mathcal{C}(g)$, there exist real numbers
$\beta$, $t$ such that $h = g^t$ is a $\pi_f$--ramified lift of the map
$x \mapsto x + \beta$.
\end{cor}

\proof  This corollary follows directly from Corollary~\ref{l=a}, Lemma~\ref{l=liftflow} and\break
Lemma~\ref{l=cflow}.  \eproof

\begin{lemma}\label{l=existsgamma} For any $f_1, f_2\in G^l \setminus \mathcal{C}(g)$,  
there exists a real number $\gamma$ such that $f_2$ is a $\pi_{f_1}$--ramified lift of 
the map $x\mapsto \lambda(f_2) x + \gamma$.
\end{lemma}

\proof The proof is expressed in a series of commutative diagrams.

\begin{lemma}\label{l=existsalpha} There exists an $\alpha\in \R$ such that 
$g^\alpha$ is the $(\pi_{f_1}, \pi_{f_2}, id)$--ramified lift of the
identity map.
\end{lemma}

\proof The following
diagram shows that if $\widehat {id}$ is the $(\pi_{f_1}, \pi_{f_2},
id)$--ramified lift of the identity map on $\RP$, 
then $\widehat {id} \circ g = g \circ \widehat id$:

\fig{
\object {a}  0,10 {S^1};
\object {b}  12,10 {S^1};
\object {c}  0,2  {\RP};
\object {d}  12,2  {\RP};
\object {e}  10,20 {S^1};
\object {f}  22,20 {S^1};
\object {g}  10,12  {\RP};
\object {h}  22,12  {\RP};
\arrow  {a}  {b}  {g};
\arrow  {c}  {d}  {x\mapsto x+1};
\arrow  {a}  {c}  {\pi_{f_1}};
\arrow  {b}  {d}  {\pi_{f_1}};
\arrow  {e}  {f}  {g};
\arrow  {g}  {h}  {x\mapsto x+1};
\arrow  {e}  {g}  {\pi_{f_2}};
\arrow  {f}  {h}  {\pi_{f_2}};
\arrow  {a}  {e}  {\widehat{id}};
\arrow  {b}  {f}  {\widehat{id}};
\arrow  {c}  {g}  {id};
\arrow  {d}  {h}  {id};}

Since $g$ embeds in a flow $g^t$ that is a ramified lift of 
an affine flow, 
it follows from Corollary~\ref{c=existsvarphi} that there exists an $\alpha$
such that $\widehat {id} = g^\alpha$.\eproof

\begin{lemma}\label{l=existsgamma0} 
For every $t_0\in\R$ there exists $\gamma\in \R$ such that
$g^{t_0}$ is the $(\pi_{f_1}, \pi_{f_2}, id)$--ramified lift
of the map $x\mapsto x - \gamma$.
\end{lemma}

\proof Let $\alpha$ be given by the previous lemma. Let
$\gamma$ be the real number 
such that $g^{t_0-\alpha}$  is a $\pi_{f_1}$ - ramified lift
of $x\mapsto x - \gamma$.  The proof follows 
from the following diagram:

\fig{
\object {a}  0,8 {S^1};
\object {b}  8,8 {S^1};
\object {c}  0,4 {\RP};
\object {d}  8,4  {\RP};
\object {e}  16,8  {S^1};
\object {f}  16,4  {\RP};
\arrow  {a}  {b}  {g^{t_0-\alpha}};
\arrow  {c}  {d}  {x\mapsto x - \gamma};
\arrow  {a}  {c}  {\pi_{f_1}};
\arrow  {b}  {d}  {\pi_{f_1}};
\arrow  {b}  {e}  {\widehat{id} = g^\alpha};
\arrow  {d}  {f}  {id};
\arrow  {e}  {f}  {\pi_{f_2}};}
%\arrow  {a}  {e}  {\widehat id = g^\alpha};
%\arrow  {c}  {f}  {id};}

The composition of the maps on the 
top row is $g^\alpha \circ g^{t_0-\alpha} = g^{t_0}$.
The composition of the maps on the bottom row
is $x \mapsto x - \gamma$. \eproof

\begin{lemma}  For all $t \in \R$, 
$f_2 g^t f_2^{-1} = g^{\lambda(f_2)t}$.  
\end{lemma}

\proof  The proof follows from the following diagram:

\fig{
\object {a}  0,8 {S^1};
\object {b}  10,8 {S^1};
\object {c}  0,4 {\RP};
\object {d}  10,4  {\RP};
\object {e}  20,8  {S^1};
\object {f}  20,4  {\RP};
\object {g}  30,8  {S^1};
\object {h}  30,4  {\RP};
\arrow  {a}  {b}  {f_2^{-1}};
\arrow  {c}  {d}  {x\mapsto \frac{1}{\lambda(f_2)}x};
\arrow  {a}  {c}  {\pi_{f_2}};
\arrow  {b}  {d}  {\pi_{f_2}};
\arrow  {g}  {h}  {\pi_{f_2}};
\arrow  {b}  {e}  {g^t};
\arrow  {d}  {f}  {x \mapsto x + t};
\arrow  {e}  {f}  {\pi_{f_2}};
\arrow  {e}  {g}  {f_2};
\arrow  {f}  {h}  {x \mapsto \lambda(f_2) x};
}
The composition of the maps on the bottom row gives $x \mapsto x + \lambda(f_2) t$.
By uniqueness (Lemma~\ref{l=liftflow}), $f_2 g^t f_2^{-1} = g^{\lambda(f_2)t}$
for all $t \in \R$.
\eproof

Let $\alpha$
be given by Lemma~\ref{l=existsalpha}, and let $\gamma$ be given by 
Lemma~\ref{l=existsgamma0}, with $t_0 = \lambda(f_2)\alpha$. 
From the following diagram:

\fig{
\object {a}  0,10 {S^1};
\object {b}  12,10 {S^1};
\object {c}  0,2  {\RP};
\object {d}  12,2  {\RP};
\object {e}  10,20 {S^1};
\object {f}  22,20 {S^1};
\object {g}  10,12  {\RP};
\object {h}  22,12  {\RP};
\arrow  {a}  {b}  {\widehat{id} = g^\alpha};
\arrow  {c}  {d}  {id};
\arrow  {a}  {c}  {\pi_{f_1}};
\arrow  {b}  {d}  {\pi_{f_2}};
\arrow  {e}  {f}  {g^{\lambda\alpha}};
\arrow  {g}  {h}  {x\mapsto x-\gamma};
\arrow  {e}  {g}  {\pi_{f_1}};
\arrow  {f}  {h}  {\pi_{f_2}};
\arrow  {a}  {e}  {f_2};
\arrow  {b}  {f}  {f_2};
\arrow  {c}  {g}  {F};
\arrow  {d}  {h}  \fliplabel{x\mapsto \lambda(f_2) x};}

it follows that $F(x) = \lambda(f_2) x + \gamma$, completing the proof
of Lemma~\ref{l=existsgamma}.\eproof

\begin{prop} Fix $f\in G^l\setminus \mathcal{C}(g)$.  Then for each $h\in G$, there
exists $F \in {\rm Aff}(\R)$ such that $h$ is a $\pi_f$--ramified lift of $F$.
\end{prop}

\proof By Corollary~\ref{c=existsvarphi} and Lemma~\ref{l=existsgamma},
we have that for each $h\in G$, there exists $k\in {\rm Aff}_+(\R)$,
so that $h^l$ is a $\pi_f$--ramified lift of $k^l$.  
Therefore, by Proposition~\ref{p=root},
$h$ is a $\pi_f$--ramified lift of either $k$ or $-k$.\eproof

By Lemma~\ref{l=crconj}, this completes the proof of 
Theorem~\ref{t=solv}.  \eproof

\addcontentsline{toc}{subsection}{Acknowledgements}
\sh{Acknowledgements}
Many useful conversations with Gautam Bharali, Christian Bonatti,
Keith Burns, Matthew Emerton,  Benson Farb, 
Giovanni Forni, John Franks, Ralf Spatzier, Jared Wunsch and 
Eric Zaslow are gratefully acknowledged.  We thank
Andr{\'e}s Navas and Etienne Ghys 
for very useful comments on earlier versions 
of this paper, and for pointing  out several references to us.
Ghys supplied the simple proof of Proposition~\ref{p=ramcovexists},
and Curt McMullen supplied the simple proof of Lemma~\ref{l=normal}.
Finally, we thank Benson Farb for reminding us that the BS groups are often 
interesting.  The second author was supported by an NSF grant.

\end{document}